\documentclass[%
 notitlepage,superscriptaddress,
 amsmath,amssymb,
 aps
 amsmath,amssymb,
aps,pra,onecolumn
]{revtex4-2}

\usepackage{graphicx}
\usepackage{dcolumn}
\usepackage{bm}
\usepackage{mathtools}

\usepackage{hyperref}

\begin{document}

\preprint{APS/123-QED}

\title{Improving the accuracy of meshless methods via resolving power optimisation using multiple kernels}

\author{H. M. Broadley}
\email{henry.broadley@manchester.ac.uk}

\author{J. R. C. King}%

\affiliation{%
 Department of Mechanical and Aerospace Engineering, The University of Manchester}%

\author{S. J. Lind}

\affiliation{School of Engineering, Cardiff University}%

\date{\today}

\begin{abstract}
Meshless methods are commonly used to determine numerical solutions to partial differential equations (PDEs) for problems involving free surfaces and/or complex geometries, approximating spatial derivatives at collocation points via local kernels with a finite size. Despite their common use in turbulent flow simulations, the accuracy of meshless methods has typically been assessed using their convergence characteristics resulting from the polynomial consistency of approximations to operators, with little to no attention paid to the resolving power of the approximation. Here we provide a framework for the optimisation of resolving power by exploiting the non-uniqueness of kernels to provide improvements to numerical approximations of spatial derivatives. We first demonstrate that, unlike in finite-difference approximations, the resolving power of meshless methods is dependent not only on the magnitude of the wavenumber, but also its orientation, before using linear combinations of kernels to maximise resolving power over a range of wavenumbers. The new approach shows improved accuracy in convergence tests and has little impact on stability of time-dependent problems for a range of Eulerian meshless methods. Solutions to a variety of PDE systems are computed, with significant gains in accuracy for no extra computational cost per timestep in Eulerian frameworks. The improved resolution characteristics provided by the optimisation procedure presented herein enable accurate simulation of systems of PDEs whose solution contains short spatial scales such as flow fields with homogeneous isotropic turbulence.

\end{abstract}

\maketitle

\section{Introduction}\addvspace{10pt}
\label{sec::Introduction}
When performing direct numerical simulations of systems of partial differential equations (PDEs) where the solution has a fine structure (e.g. turbulent flows in fluid dynamics), resolving how the solution varies across the smallest wavelength admitted by the discretisation is often essential to obtaining an accurate solution. In many cases, optimising the resolution of all length scales permitted by a given discretisation (i.e. the resolving power of a numerical scheme), is preferred to improving the polynomial consistency of numerical approximations. Resolving power is acknowledged to be a key metric in the finite-difference community, where compact (i.e. implicit) schemes with spectral-like resolving power (see \cite{Lele}) are commonly preferred to higher-order discretisations (both explicit and implicit) with worse resolving power.

When solving PDEs in complex geometries or with a moving surface, meshless methods, which can discretise the domain trivially, are a natural choice --- for a broad overview we refer the reader to \cite{Li_Review, Garg_Review}. Throughout this manuscript we refer to collocation points as nodes, define our local stencil at each node $i$ to be the set of nodes $\mathcal{N}_i$. At each node we must choose a kernel/basis function/set of basis functions (depending on the meshless method in question) which enables local weights to be calculated. These weights are subsequently used to approximate derivative operators. For ease of exposition we refer to the functions which produce weights as kernels throughout.

The most widely known meshless method is Smoothed Particle Hydrodynamics (SPH). SPH is an interpolation scheme originally developed for astrophysical applications \cite{gingold,Lucy} but is now also a commonly used method within the wider fluid dynamics community, particularly for free surface flows \cite{Monaghan_Review}. Spatial derivatives are approximated using weighted sums of properties at neighbouring nodal points, with the weights determined by a smoothing kernel. Modern SPH methods typically use polynomial kernels with compact support. For a discrete set of nodal points SPH is formally zeroth order accurate for gradient operators and divergent for the Laplacian, though convergent-like behaviour is obtained for poorly resolved domains as higher order errors typically dominant.

For many years, improvements to SPH approximation have focused on improving the polynomial consistency of approximation \cite{Bonet&Lok,Sibilla}, ignoring resolving power as a method of quantifying error. Recent advances in SPH have exploited the non-uniqueness of weights by using multiple kernels to approximate operators, providing a route to optimisation of a chosen metric. In \cite{Cabezon&Garcia} a linear combination of kernels was used to minimise the pairing instability found in SPH simulations, while \cite{MK_Arxiv} used multiple kernels to minimise two different cost functions. In \cite{MK_Arxiv}, the improved kernels were used to solve the Navier--Stokes equations, with $L_1$-error reductions up to $40\%$ observed in several test cases. Due to the modified SPH methods (geometric density average force SPH (GDSPH) and integral SPH
(ISPH)) used in \cite{MK_Arxiv} relatively large stencil sizes were typically required for this study compared to standard SPH implementations, and the choice of cost function was seen to have a large impact on the accuracy of the approximations. 

The Generalised Finite-Difference Method (GFDM) is an extension of finite-differences to scattered nodal distributions. When approximating differential operators, weights are determined by solving a linear system to ensure polynomial consistency up to a desired order. However, GFDM is typically limited to low order due to ill-conditioning of the linear system and poor accuracy. In contrast, the Generalised Moving Least Squares (GMLS) method of \cite{TraskGMLS} uses polynomial reconstruction to provide high-order approximations to differential operators. GMLS has recently been extended to a compact scheme (CMLS) \cite{TraskCMLS}. Therein, two different approaches to obtaining implicit approximations were considered: (i) so-called Weinan type derivatives \cite{Weinan} which require information about the PDE to obtain the compact approximation; (ii) operators similar to those of \cite{Lele} (called Lele-type derivatives in \cite{TraskCMLS}). For both types of compact schemes, clear improvements in error obtained when comparing GMLS and CMLS schemes of the same polynomial accuracy. 

Radial Basis Function (RBF) methods, first developed in the field of cartography \cite{RBF_Cartography}, can also be used to interpolate functions and approximate derivatives on scattered point clouds. Early RBF methods were global in nature, requiring information from every point in the domain to obtain weights. These required the solution of large, dense linear systems which were often poorly-conditioned and costly to solve, negating the exponential convergence and spectral accuracy provided. To mitigate these issues, a local approach known as Radial Basis Function-Finite Differences (RBF-FDs), developed independently by \cite{shu},\cite{Cecil}, and \cite{Wright}, is commonly preferred when solving PDE systems. In RBF-FDs only a local stencil is used, and to improve convergence characteristics polynomial consistency can be prescribed to any given order (as usual, higher order convergence requires larger stencils/more neighbours). An overview of the theory underpinning RBFs and RBF-FDs as well as various applications can be found in \cite{Fornberg&Flyer2015}.

Recent work by \cite{Hybridkernels_RBFFD} on RBF-FDs by has focused on `hybrid' kernels, where a linear sum of two RBFs is used. This approach therefore differs from a multiple kernel approximation as only one linear system is solved (and therefore only one set of weights produced), as opposed to obtaining the final set of weights from a linear combination of two previously calculated sets.

The Local Anisotropic Basis Function Method (LABFM) is another high-order meshless method which has been developed in recent years \cite{king_2020}. LABFM approximates operators to arbitrary polynomial consistency by using a weighted sum of anisotropic basis functions (ABFs), with weights determined by solving a local linear system. Further improvements to LABFM are detailed in \cite{king_2022} where using orthogonal polynomials as ABFs and a stencil reduction algorithm was demonstrated to increase accuracy and reduce computational cost whilst ensuring the linear system remained well-conditioned. LABFM has been applied to weakly compressible non-Newtonian flows \cite{king_2024_ve}, and fully compressible combustion flows \cite{king_2024_Combustion}, yielding high-order solutions in complex geometries.

 Despite the widespread use and ongoing development of meshless methods, the notion of resolving power has received sparse attention from the community, in large part due to the larger number of neighbours required to compute derivatives compared to finite-differences. Consequently, the focus has typically been on the trade-off between computational cost and polynomial consistency. In addition, the choice of kernel which generates the (non-unique) weights is, although informed by experience, often ad-hoc and not rigorously justified. Herein we present an approach by which one may use the notion of resolving power to improve the accuracy of a meshless method by using weights generated by multiple kernels. 

 The structure of this paper is as follows. In \S\ref{sec::resolving_power_formulation} we outline the notion of resolving power for disordered nodal distributions. In \S\ref{sec::linear_combinations_of_weights} we then show how using linear combination of kernels can be used to improve resolving power of meshless methods, and in \S\ref{sec::Convergance_Studies} perform convergence studies for three different methods, comparing the accuracy of approximations to the single-kernel formulation. Several PDE systems are then solved  in \S\ref{sec:PDE_solns} using one method (LABFM), with comparisons for solutions using single and multi-kernel approaches. Details of the numerical implementations of each meshless method used can be found in Appendix \ref{sec::Numerical_implementations}.

 In the results which follow, where we present comparisons with single-kernel approaches, unless otherwise stated, the single-kernel approaches use the following smoothing kernel/basis functions: SPH - Wendland C2 kernel; RBF-FDs - Gaussian RBF; LABFM - Hermite polynomials with Wendland C2 kernel. These have been chosen due to their widespread use in the relevant literature, though we expect similar results for other choices. The second kernels used in the multi-kernel approach of \S\ref{sec::linear_combinations_of_weights} were found to perform slightly worse in convergence studies than those chosen for the single-kernel approach. Details of the second kernel used for each numerical method can be found in Appendix \ref{sec::Numerical_implementations}.

\section{Resolving Power of Meshless Methods}
\label{sec::resolving_power_formulation}
The use of resolving power, a form of Fourier analysis, as a metric to determine the efficacy of collocated numerical methods dates back to the 1960s~\cite{Roberts&Weiss}. A thorough discussion how one may use Fourier analysis to characterise errors of numerical methods can be found in \cite{Vichnevetsky&Bowles}. The essence of this analysis is to determine how well a numerical method approximates the desired derivatives of waves for a range of wavenumbers up to the Nyquist wavenumber (the maximal wavenumber admitted by the discretisation). For a function with a given wavenumber $k_x$, the numerical method will calculate a derivative which has an effective wavenumber $k_{eff}$. If the numerical scheme exactly approximates derivatives (e.g. as in spectral schemes) then $k_{eff}=k_x$ for all values of $k_x$. Such an analysis is typically performed for one-dimensional waves (e.g. for derivatives with respect to $x$ a resolving power analysis considers the derivative of the function $e^{ik_xx}$) and directional dependence is only considered in the context of reducing numerical anisotropy rather than minimising resolving power errors of individual operators. However, most simulations of PDE systems are performed in two or three dimensions, therefore if computed derivatives have a dependence on wavenumbers associated with other dimensions this must be included in a resolving power analysis. In what follows, we use consider a two-dimensional domain, though extension to higher dimensions follows naturally.

Let us consider a domain $\Omega\subset \mathbb{R}^2$. To ease the Fourier analysis we assume that $\Omega=[-\pi L_x,\pi L_x]\times[-\pi L_y,\pi L_y]$ is rectangular and that the system we are solving is subject to periodic boundary conditions on $\partial \Omega$, although non-periodic boundaries can be dealt with relatively straightforwardly --- see \cite{Lele} for examples. 

We begin by defining a function $\phi:\Omega \to\mathbb{R}$ which can be expressed as a two-dimensional Fourier series
\begin{equation}
\label{eq::fourier_series}
    \phi(x,y)=\sum_{k_x\in \mathbb{Z},}\sum_{k_y\in \mathbb{Z}} c_{k_x,k_y}e^{ i(k_xx+k_yy)},
\end{equation}
where the $c_{k_x,k_y}$ are complex constants and $k_x,k_y$ are (appropriately scaled) wavenumbers. Let us now discretise $\Omega$ using $N$ nodal points, making no assumptions about the distribution of nodes. Explicit collocated numerical methods approximate differential operators $\mathcal{L}$ acting on the function $\phi$ at a point $i$ using weighted sums of function properties at nodes within the computational stencil of the point $(x_i,y_i)$, therefore can be written in the form
 \begin{equation}
        \mathcal{L}(\phi)\Big|_i=\sum_{j\in\mathcal{N}_i} \phi_{ji} w_{ji}^L.\label{eq::generic_operator}
    \end{equation}
where $\phi_{ji}\coloneqq \phi_j-\phi_i$, $\mathcal{N}_i$ is the set of nodes inside the computational stencil of point $i$ and $w_{ji}^L$ are the set of weights. This equation describes finite-differences, the antisymmetric formulation of SPH, RBF-FDs, and LABFM (see Appendix \ref{sec::Numerical_implementations}).

When the domain $\Omega$ is discretised by a finite number of nodal points this limits the maximum wavenumber (or equivalently the shortest wavelength) present in the numerical approximation of the derivatives of $\phi$. This maximum wavenumber is the Nyquist wavenumber $k_{Ny}$. This wavenumber may vary in the different dimensions of the domain. For example, if our numerical method is uniform finite-differences but we discretise our domain with different nodal spacings in the $x$ and $y$ directions, then the Nyquist wavenumber will be different in each direction. As the focus of this paper is on meshless methods we shall assume that the distribution of points is locally isotropic, though the average nodal spacing may vary through the domain. The local average nodal spacing for the point $i$ is forthwith denoted $s_i$, and the (local) Nyquist wavenumber $k_{Ny}=\pi/s_i$ is the same for all modes in the Fourier series. We note that the resolving power analysis performed below is resolution independent, though of course increasing resolution increases the Nyquist wavenumber.

We are now ready to consider the resolving power of two-dimensional collocated numerical schemes which may be expressed in the form of \eqref{eq::generic_operator}. As an illustrative example, we consider the resolving power of the first derivative with respect to $x$ at point $i$ in our domain (implicit approximations will be considered in a future study). This operator is approximated through the expression
\begin{equation}
\label{eq::grad_approximation}
    \frac{\partial \phi}{\partial x}\Big|_i=\sum_{j\in\mathcal{N}_i} \phi_{ji} w_{ji}^x.
\end{equation}
Upon substitution of the Fourier series \eqref{eq::fourier_series} into \eqref{eq::grad_approximation} one finds that for each pair of wavenumbers $(k_x,k_y)=\boldsymbol{k}$ the effective wavenumber $k_{eff}$ of the scheme at point $i$ is
\begin{equation}
\label{eq::grad_resolving_power}
    k_{eff}=\sum_{j\in\mathcal{N}_i}\sin (k_x x_{ji}+k_y y_{ji})w_{ji}^x+i\sum_{j\in\mathcal{N}_i}\left[1-\cos(k_xx_{ji}+k_yy_{ji})\right] w_{ji}^x.
\end{equation}
The two terms on the right-hand side of this expression relate to dispersion and dissipation of the numerical scheme. If $\Re\{k_{eff}\}\neq k_x$ then the numerical scheme has a dispersion error, and if $\Im\{k_{eff}\}\neq 0$ then there is a dissipation error. Some numerical methods have properties which simplify the above expression. If we use a Cartesian mesh which is aligned with the coordinated axes and use a computational stencil which is a line (as in finite-differences) then $y_{ji}=0$ for all ${j\in\mathcal{N}_i}$. The effective wavenumber of the scheme in this special case is independent of $k_y$, and a one-dimensional Fourier analysis is justified, as in the analysis of (for example),\cite{Vichnevetsky&Bowles}, \cite{Lele}, and \cite{Pencil_code}. We also note that uniform central finite-differences have the property 
\begin{equation}
    \sum_{j\in\mathcal{N}_i}\left[1-\cos(k_xx_{ji}+k_yy_{ji})\right] w_{ji}^x\equiv0
\end{equation}
and thus have no dissipative error.

For general meshless methods however, \eqref{eq::grad_resolving_power} describes the effective wavenumber of the approximation at point $i$. We note that we have
\begin{gather}
\label{eq::RP_real_properties}
    \Re\{k_{eff}(k_x,k_y)\}=-\Re\{k_{eff}(-k_x,-k_y)\},\qquad \Re\{k_{eff}(-k_x,k_y)\}=-\Re\{k_{eff}(k_x,-k_y)\},\\
    \label{eq::RP_im_properties}
    \Im\{k_{eff}(k_x,k_y)\}=\Im\{k_{eff}(-k_x,-k_y)\},\qquad \Im\{k_{eff}(-k_x,k_y)\}=\Im\{k_{eff}(k_x,-k_y)\},
\end{gather}
therefore we can determine the resolving power by considering a half-disc in $(k_x,k_y)$ space such that $k_x^2+k_y^2\leq k_{Ny}^2$.

By the nature of meshless methods, both the size of the stencil and the location of nodal points within each stencil may vary significantly on a node-by-node basis, consequently resolving power differs from node to node. To obtain a robust measure of the resolving power of a given method across a nodal distribution we compare the average effective wavenumber of all nodes, although we note the procedure outlined in $\S$\ref{sec::linear_combinations_of_weights} is implemented on a stencil-by-stencil basis. The resolving power of a given numerical method is of course dependent on the initial nodal distribution in Eulerian methods, and time-dependent for Lagrangian methods, though for nodal distributions with a sufficiently large number of nodes we see negligible differences in the (Eulerian) results presented herein.

When evaluating the resolving power of a numerical method one typically compares the difference between the real and effective wavenumbers for $\lvert\lvert\boldsymbol{k}\rvert\rvert\in[0,k_{Ny}]$, which correspond to lines through the origin in $(k_x,k_y)$ space. To enable comparison between different such lines on a given discretisation we introduce the scaled real and effective wavenumbers, $\hat{k}$ and $\hat{k}_{eff}$, defined through
\begin{equation}
\label{eq:hatk_defn}
\hat{k}=\frac{||\boldsymbol{k}||}{k_{Ny}}k_x,\qquad
    \hat{k}_{eff}(k_x,k_y)=\frac{||\boldsymbol{k}||}{k_{Ny}}k_{eff}(k_x,k_y).
\end{equation}

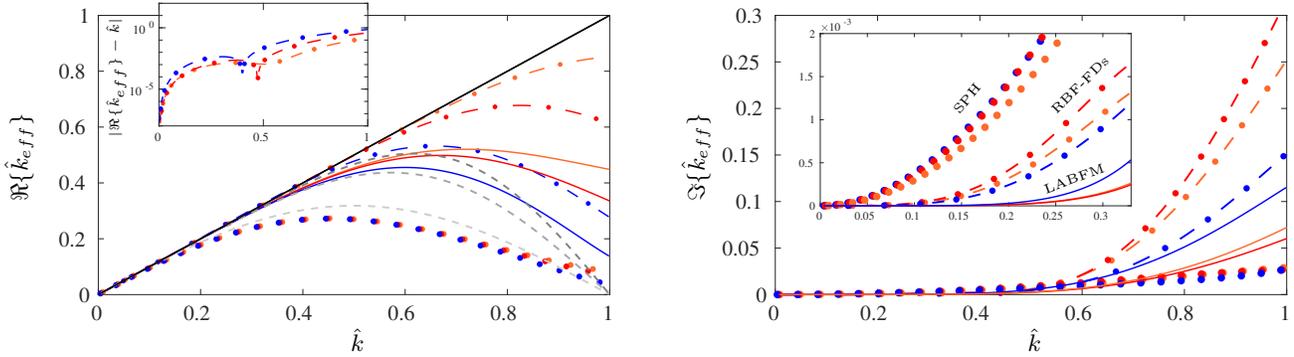
\begin{figure}
    \begin{center}
    \setlength{\unitlength}{1cm}
    \begin{picture}(18,5)(0,0)
    \put(0,0){\includegraphics[width=0.49\linewidth]{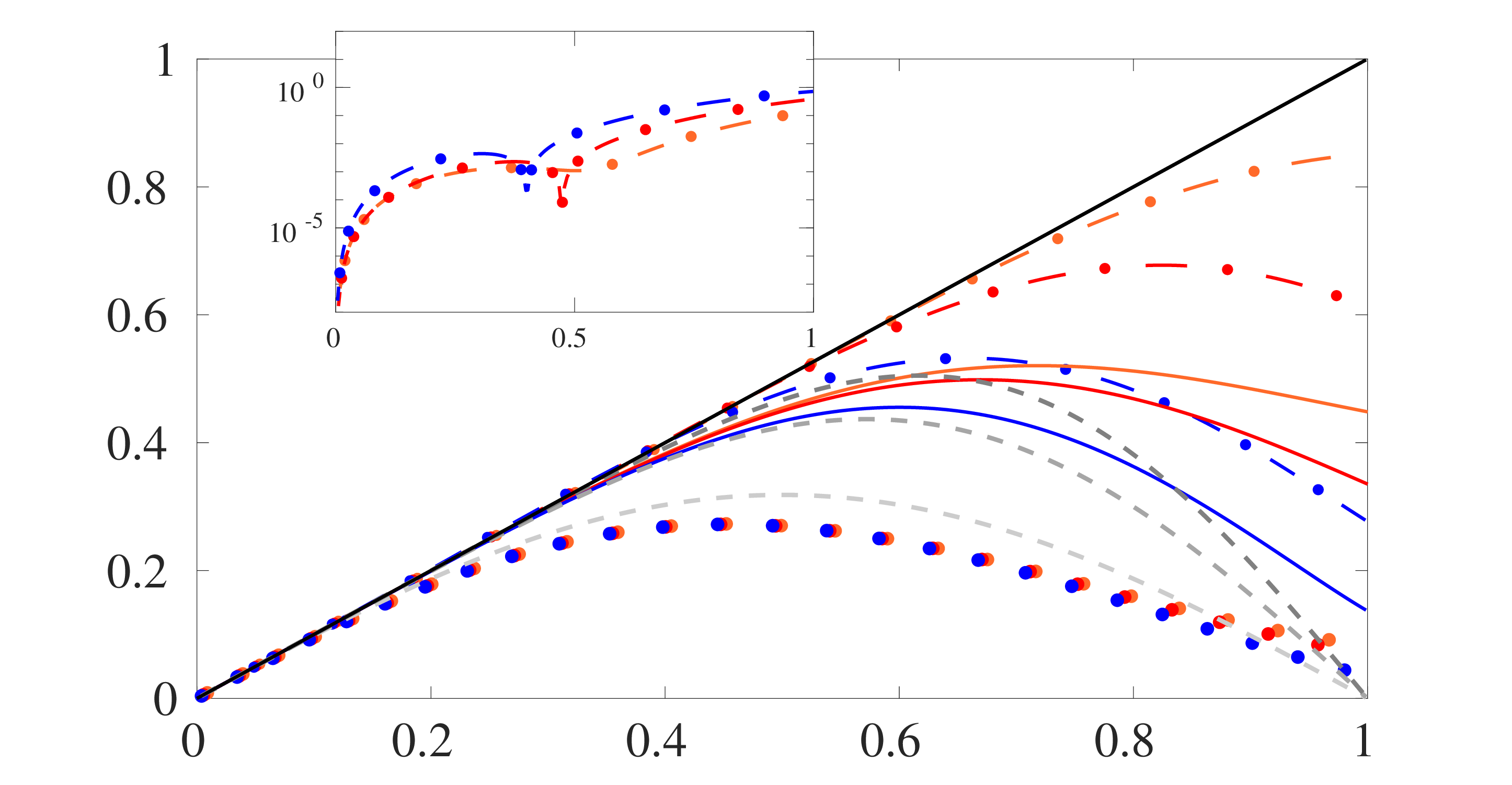}}
    \put(4.5,-0.25){$\hat{k}$}
    \put(9,00){\includegraphics[width=0.49\linewidth]{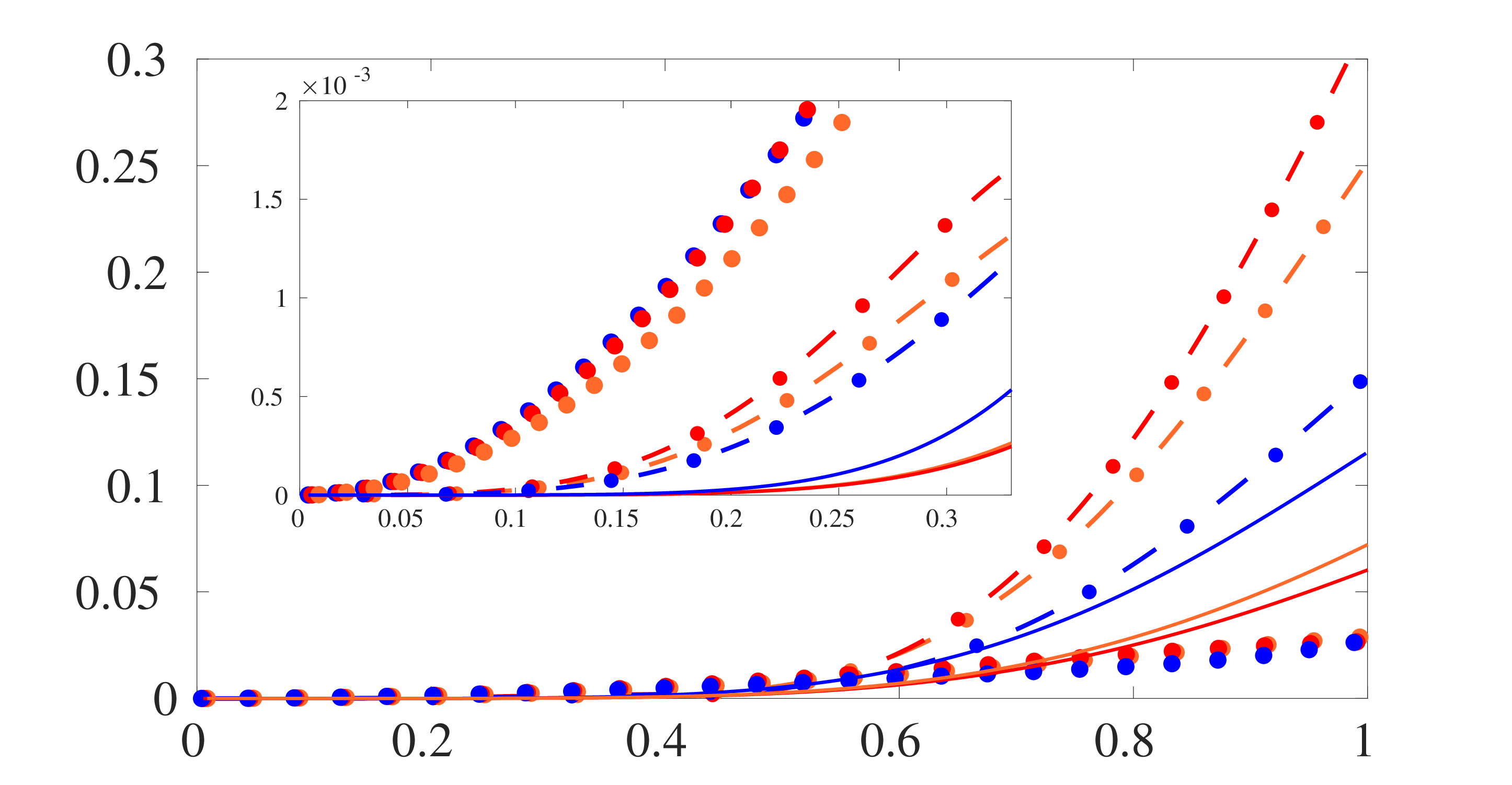}}
    \put(13.5,-0.25){$\hat{k}$}
  
  \put(12.5,2.9){\rotatebox{50}{\tiny{SPH}}}
   \put(13.7,1.85){\rotatebox{20}{\tiny{LABFM}}}
   \put(13.8,2.9){\rotatebox{45}{\tiny{RBF-FDs}}}
    \put(-0.1,1.75){\rotatebox{90}{$\Re\{\hat{k}_{eff}\}$}}
    \put(8.9,1.75){\rotatebox{90}{$\Im\{\hat{k}_{eff}\}$}}
    \put(1.25,2.6){\rotatebox{90}{{\tiny{$\lvert \Re\{\hat{k}_{eff}\}-\hat{k}\rvert$}}}}

    \end{picture}
    \end{center}
    \caption{Resolving power of the gradient operator for various meshless methods. Left panel $\Re\{k_{eff}\}$, right panel $\Im\{k_{eff}\}$. Orange lines denote $l=2k$, red $l=k$, blue $l=0$, black spectral accuracy. $\cdots$ SPH with Wendland C2 kernel $h/s=1.3$, $\cdot-\cdot$ $m=2$ RBF-FDs with Gaussian RBF $\lvert\mathcal{N}_i\rvert=20$, --- $m=4$ LABFM using Hermite polynomials and Wendland C2 kernel. Grey-scale dashed lines denote $2^{\mathrm{nd}}$, $4^{\mathrm{th}}$, and $6^{\mathrm{th}}$ order FDs for comparison (increasing resolving power corresponds to higher order). Left inset shows semi-log plot of $\hat{k}$ against $\lvert \hat{k}-\Re\{\hat{k}_{eff}\}\rvert$, for RBF-FDs. Right inset shows zoomed in graph of $\Im\{\hat{k}_{eff}\}$ for the range $0\leq \hat{k}\leq 1/3$.}
    \label{fig:rp_grad_examples}
\end{figure}

To give an example of the analysis carried out above, Fig. \ref{fig:rp_grad_examples} shows the real (left panel) and imaginary (right panel) parts of $\hat{k}_{eff}$ demonstrating resolving power along different lines in $(k_x,k_y)$ space of various collocated meshless numerical methods. There are substantial differences between the different methods. Considering the real part of $\hat{k}_{eff}$, RBF-FDs (dash-dot lines) are closest to spectral resolving power for high wavenumbers, however it can be seen in the inset that $\hat{k}_{eff}>\hat{k}$ across a wide range of wavenumbers which can cause issues in resolving acoustic waves, leading to catastrophic breakdown of the numerical solution. In addition, the imaginary part of the wavenumber is significantly larger than the other meshless methods, indicating that the dissipation error is significantly smaller. The resolving power of LABFM at $m=4$ (solid lines) remains close to spectral up to $\hat{k}=0.4$ (similar to $6^{\mathrm{th}}$ order FDs), with $\hat{k}_{eff}<\hat{k}$ across all wavenumbers. For both LABFM and RBF-FDs there is a substantial dependence on the angle of the line, though the sign of $k_y$ was noted to have only a small influence. In SPH (dotted lines), the real part of $\hat{k}_{eff}$ deviates from spectral accuracy very quickly and the imaginary part (dissipation error) is also significantly larger over the first third of wavenumbers. For high wavenumbers SPH has a smaller dissipation error than the other methods, although when solving PDE systems de-aliasing usually prevents accurate resolution of wavenumber $\hat{k}>\frac{2}{3}$ to ensure stability.

\begin{figure}
    \begin{center}
    \setlength{\unitlength}{1cm}
    \begin{picture}(18,5)(0,0)
    \put(0,0){\includegraphics[width=0.49\linewidth]{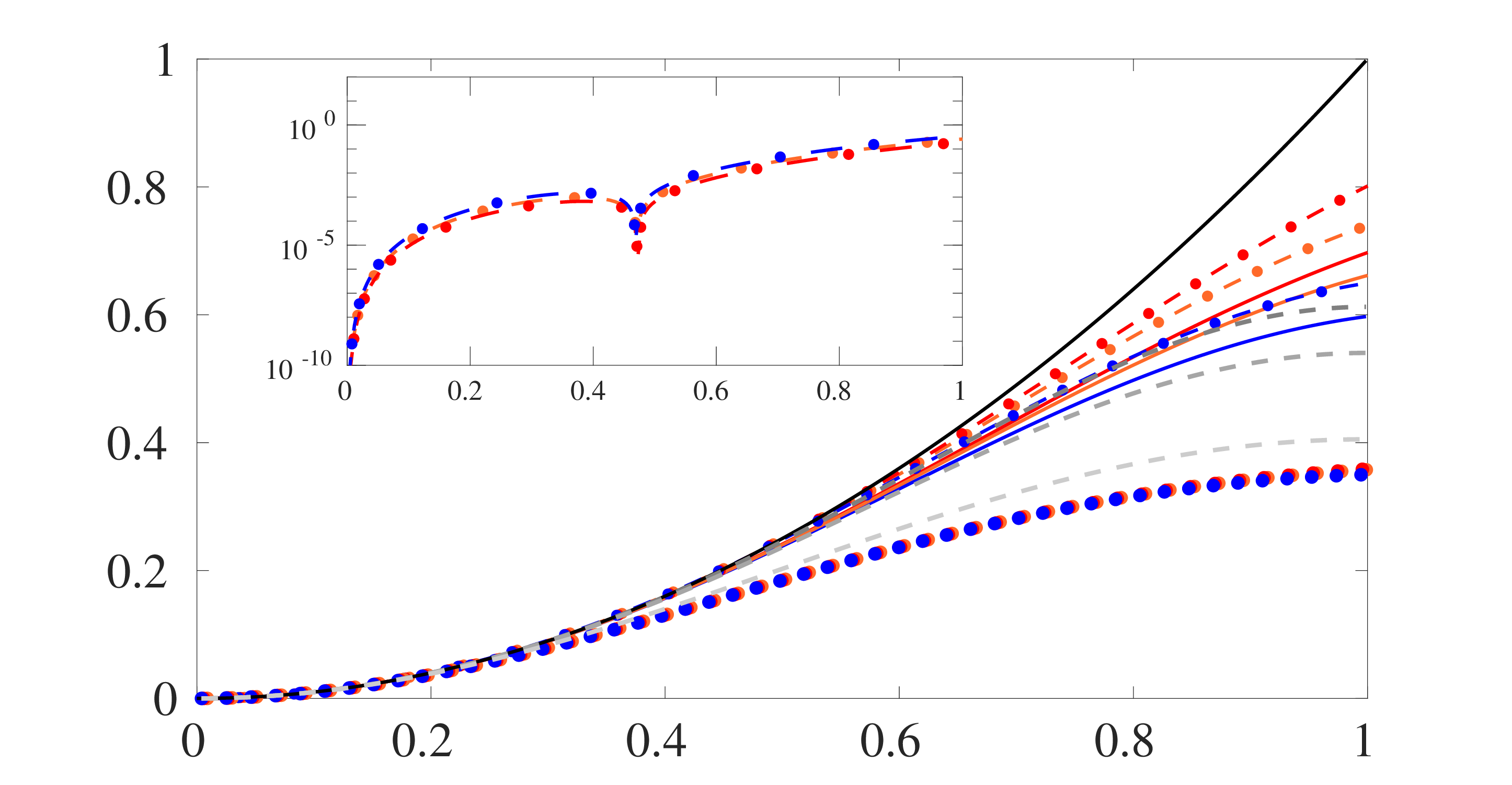}}
    \put(9,0){\includegraphics[width=0.49\linewidth]{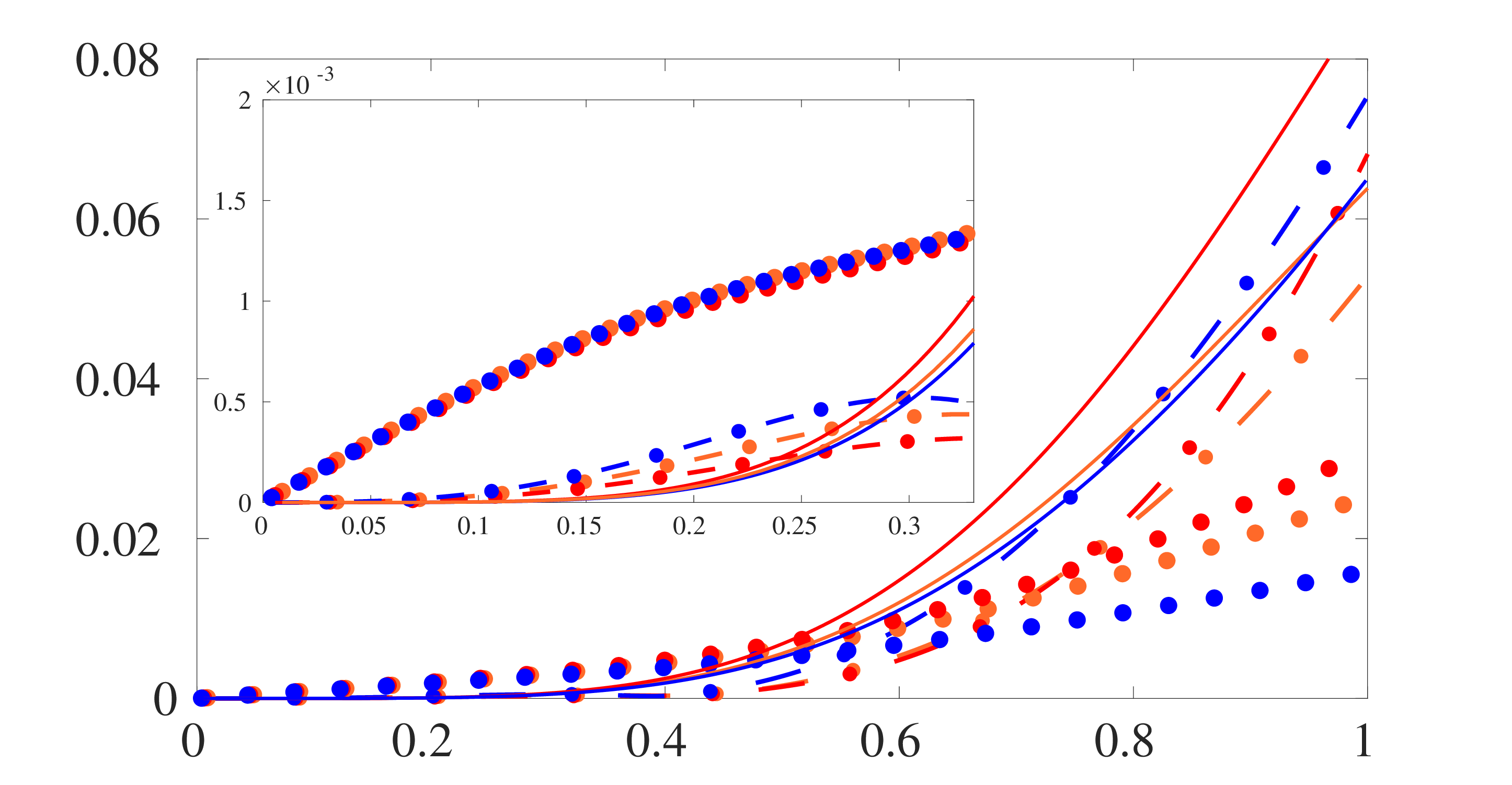}}
 \put(1.25,2.4){\rotatebox{90}{{\tiny{$\lvert \Re\{\hat{q}_{eff}^2\}-\hat{q}^2\rvert$}}}}
    \put(13.5,-0.25){$\hat{q}$}
    \put(-0.1,1.75){\rotatebox{90}{$\Re\{\hat{q}^2_{eff}\}$}}
    \put(8.9,1.75){\rotatebox{90}{$\Im\{\hat{q}^2_{eff}\}$}}

    \put(4.5,-0.25){$\hat{q}$}

    \end{picture}
    \end{center}
    \caption{Resolving power graphs of Laplace operator for various meshless methods. Left panel $\Re\{q_{eff}\}$, right panel $\Im\{q_{eff}\}$. Orange lines denote $l=2k$, red $l=k$, blue $l=0$, black spectral accuracy. $\cdots$ SPH with Wendland C2 kernel $h/s=1.3$, $\cdot-\cdot$ $m=2$ RBF-FDs with Gaussian RBF $\lvert\mathcal{N}_i\rvert=20$, --- $m=4$ LABFM using Hermite polynomials and Wendland C2 kernel. Grey-scale dashed lines denote $2^{\mathrm{nd}}$, $4^{\mathrm{th}}$, and $6^{\mathrm{th}}$ order FDs for comparison (increasing resolving power corresponds to higher order). Left inset shows semi-log plot of $\hat{q}$ against $\lvert \hat{kq}^2-\Re\{\hat{q}_{eff}^2\}\rvert$, for RBF-FDs. Right inset shows zoomed in graph of $\Im\{\hat{q}_{eff}^2\}$ for the range $0\leq \hat{q}\leq 1/3$.}
    \label{fig:rp_lap_examples}
\end{figure}

The analysis carried above can be easily extended to higher derivatives. For example, the resolving power of the Laplacian operator can be expressed through
\begin{equation}
    q_{eff}^2=\sum_{j\in\mathcal{N}_i}(1-\cos(k_xx_{ji}+k_yy_{ji}))w_{ji}^{L}-i\sum_{j\in\mathcal{N}_i}\sin(k_xx_{ji}+k_yy_{ji})w_{ji}^{L}
\end{equation}
where $q_{eff}^2=k_x^2+k_y^2$ for spectral methods. Here the real part of $q_{eff}^2$ relates to the degree of dissipation, and the imaginary part is dispersion error. Instead of\eqref{eq::RP_real_properties} \& \eqref{eq::RP_im_properties} we have
\begin{gather}
\label{eq::RPlap_real_properties}
    \Re\{q_{eff}^2(k_x,k_y)\}=\Re\{q_{eff}^2(-k_x,-k_y)\},\qquad \Re\{q_{eff}^2(-k_x,k_y)\}=\Re\{q_{eff}^2(k_x,-k_y)\},\\
    \label{eq::RPlap_im_properties}
    \Im\{q_{eff}^2(k_x,k_y)\}=-\Im\{q_{eff}^2(-k_x,-k_y)\},\qquad \Im\{q_{eff}^2(-k_x,k_y)\}=-\Im\{q_{eff}^2(k_x,-k_y)\},
\end{gather}
in the present case. Analogous to the case of the first derivative approximation we define a scaled real and effective wavenumber
\begin{equation}
    \hat{q}=\frac{||\boldsymbol{k}||}{k_{Ny}},\qquad \hat{q}_{eff}^2(k_x,k_y)=\frac{q_{eff}^2}{ k_{Ny}}.
\end{equation}
Fig. \ref{fig:rp_lap_examples} shows the resolving power for the Laplacian operator for the same meshless methods as in Fig. \ref{fig:rp_grad_examples}. The resolving power of RBF-FDs and LABFM is similar, though $q_{eff}^2>k_x^2+k_y^2$ across a wide range of wavenumbers along all lines in $(k_x,k_y)$ space for RBF-FDs. SPH again has significantly worse resolving power with little difference between the lines in $(k_x,k_y)$ space.

Now that we have a framework for evaluating resolving power for numerical methods on arbitrary discrestisations, we discuss how to improve accuracy of approximations.

\section{Linear Combinations of Weights}
\label{sec::linear_combinations_of_weights}
In many numerical methods, we have more points inside the stencil than terms in the Taylor series that need to be eliminated to realise the desired polynomial consistency. Consequently, the weights $w_{ji}^d$ are non-unique. For example, in SPH a different choice of kernel would result in different weights for a given particle distribution. If the kernels are non-unique for a given stencil, a natural question to ask is whether combinations of kernels may be used to provide some benefit - for example increased accuracy. In compact finite-difference schemes, the maximal polynomial consistency possible for a given stencil is often eschewed - resulting in a non-unique kernel - which may then be optimised to maximise resolving power. Here we investigate how one may optimise the choice of weights in meshless methods to maximise resolving power.

\subsection{Theoretical Framework}
\label{subsec::theoretical_framework}
To illustrate the theoretical idea behind this approach, consider a node $i$ which uses a computational stencil within which there are $\lvert\mathcal{N}_i\rvert$ neighbours with average nodal spacing $s_i$ to approximate derivatives. As a representative example, let us approximate the partial derivative of a given function $\phi$ (defined at all discrete points within the stencil) with respect to $x$. Assuming that there are more points in the stencil than equations we require for the approximation to have polynomial consistency of order $m$ (i.e. $\lvert\mathcal{N}_i\rvert>(m^2+3m)/2)$ in 2D) we can approximate $\partial \phi/\partial x$ to polynomial consistency $m$ using two different kernels producing two different sets of weights, $\hat{w}^x_{ji}$ and $\overline{w}^x_{ji}$ respectively. We may write
\begin{subequations}
\label{eq::different_basis_approximations}
    \begin{gather}
    \label{eq::different_basis_approximations1}
        \frac{\partial \phi}{\partial x}\Big|_i=\sum_{j\in\mathcal{N}_i} \phi_{ji} \hat{w}^{x}_{ji}+\mathcal{O}(s_i^{m+1}),\\
        \label{eq::different_basis_approximations2}
        \frac{\partial \phi}{\partial x}\Big|_i=\sum_{j\in\mathcal{N}_i} \phi_{ji} \overline{w}^{x}_{ji}+\mathcal{O}(s_i^{m+1})
    \end{gather}
\end{subequations}
Adding a linear combination of the two approximations together we obtain
\begin{equation}
    \label{eq::multi-basis_approximation_arb_c}
    (\hat{c}_i+\overline{c}_i)\frac{\partial \phi}{\partial x}\Big|_i=\sum_{j\in\mathcal{N}_i} \phi_{ji}(\hat{c}_i\hat{w}_{ji}^x+\overline{c}_i\overline{w}_{ji}^x)+\mathcal{O}(s_i^{m+1}),
\end{equation}
where $\hat{c}_i$, $\overline{c}_i$ are (at this point) arbitrary constants. Requiring that $\overline{c}_i=1-\hat{c}_i$ gives
\begin{equation}
    \label{eq::multi-basis_approximation}
    \frac{\partial \phi}{\partial x}\Big|_i=\sum_{j\in\mathcal{N}_i} \phi_{ji}\left[\hat{c}_i\hat{w}_{ji}^x+(1-\hat{c}_i)\overline{w}_{ji}^x\right]+\mathcal{O}(s_i^{m+1}).
\end{equation}
This approach gives an approximation of the operator to the same accuracy as both of those in \eqref{eq::different_basis_approximations}, however we now have a parameter, $\hat{c}_i$ which we are free to choose. In an Eulerian framework, if we choose $\hat{c}_i$ in a pre-processing step, then in a time-dependent simulation there is {\em{no}} additional cost per time-step of using a linear combination of weights. For Lagrangian methods, there is a computational cost associated with calculating two sets of weights and choosing $\hat{c}_i$, however for simulations of most PDE systems this is likely to only result in a small increase in overall cost of each time-step; calculation of all derivatives and neighbour finding algorithms are typically much more computationally expensive. 

This approach may be applied to any method with derivative operators of the form~\eqref{eq::generic_operator}, including SPH, LABFM and RBF-FDs. As we shall now use operators with more than one kernel/basis, we term this new approach a `multi-kernel' (MK) implementation, and refer to the usual approximations as single-kernel (SK) implementations.

\subsection{Optimisation of coefficients}
\label{subsec::Optimisation_of_coefficients}
A natural question now arises: how should the extra parameter $\hat{c}_i$ be used to improve approximations? We outline one such approach based on improving resolving power.

We first consider a gradient operator. We remind the reader that for first derivatives we wish (as far as possible) to have
\begin{equation}
    \Re\{k_{eff}\}=k_x,\qquad \Im\{k_{eff}\}=0
\end{equation}
over as large a range of wavenumbers $(k_x,k_y)$ as possible. We now define a scalar measure of the resolving power, $E_G$, as a double integral wavenumber space $\boldsymbol{k}$
\begin{equation}
\label{eq::rp_integral_minimiser}
\begin{split}
 E_G(\hat{c_i})&=\int_{\lvert\boldsymbol{k}\rvert\leq k_M}\left(k_x-\sum_{j\in\mathcal{N}_i}\sin(k_x x_{ji}+k_yy_{ji})\left[\hat{c}_i\hat{w}_{ji}^x+(1-\hat{c}_i)\overline{w}_{ji}^x\right]\right)^2 \\&
    +\left(\sum_{j\in\mathcal{N}_i}(1-\cos(k_xx_{ji}+k_yy_{ji}))\left[\hat{c}_i\hat{w}_{ji}^x+(1-\hat{c}_i)\overline{w}_{ji}^x\right]\right)^2\mathrm{d}\boldsymbol{k},
        \end{split}
\end{equation}
where $k_M\leq k_{Ny}$ is an upper limit of the wavenumbers we wish to optimise over. Note that \eqref{eq::rp_integral_minimiser} is quadratic in $\hat{c}_i$, therefore it is trivial to derive an analytic expression for the (unique) choice of $\hat{c}_i$ which minimises $E_G$. This value can then be calculated by approximating the appropriate integrals using any numerical method (e.g. the midpoint rule).

In minimising $E_G$ in the manner described above there is no preference for signs of the dispersive or dissipative errors, although the addition of weights to penalise different error types is a trivial extension. We note that if $k_{eff}>k_x$ for a range of wavenumbers then one may think of the numerical discretisation as `stiffening' the propagating medium. In simulations of time-dependent PDE systems (especially fluid dynamics simulations) such a property increases the likelihood of the numerical solution diverging due to fast travelling waves (e.g. acoustic waves). For this reason, it is in general preferable that $k_{eff}<k_x$ $\forall{k}_{x}$. To ensure this holds a gradient descent method using a penalty function could be implemented, though in the interests of clarity of exposition we have eschewed this approach in this study.

A similar approach can be taken to analyse the error in wavenumber space of the Laplacian operator. Here we wish
\begin{equation}
    \Re\{q_{eff}^2\}=k_x^2+k_y^2,\qquad \Im\{q_{eff}^2\}=0,
\end{equation}
and thus the error can now be expressed as 
\begin{equation}
\label{eq::rp_integral_minimiser_lap}
\begin{split}
 E_L(\hat{c_i})&=\int_{\lvert\boldsymbol{k}\rvert\leq k_M} \left(k_x^2+k_y^2-\sum_{j\in\mathcal{N}_i}(1-\cos(k_xx_{ji}+k_yy_{ji}))\left[\hat{c}_i\hat{w}_{ji}^L+(1-\hat{c}_i)\overline{w}_{ji}^L\right]\right)^2 \\&
    +\left(\sum_{j\in\mathcal{N}_i}\sin(k_xx_{ji}+k_yy_{ji})\left[\hat{c}_i\hat{w}_{ji}^L+(1-\hat{c}_i)\overline{w}_{ji}^L\right]\right)^2\mathrm{d}\boldsymbol{k}.
        \end{split}
\end{equation}
Unless otherwise stated, we shall minimise $E_G, E_L$ over wavenumbers in the range $k_M= 0.3 k_{Ny}$ for RBF-FDs and LABFM, and $k_M= 0.2 k_{Ny}$ for SPH operators. Larger choices of $k_M$ tend to lead to values of $\hat{c}_i$ where $\Re\{k_{eff}\}>k$ and/or $\Re\{q_{eff}^2\}>k_x^2+k_y^2$ for small wavenumbers (where absolute error is comparatively small). The upper limit of the integral is therefore chosen to maximise the range of wavenumbers optimised whilst retaining accuracy at small wavenumbers.

\section{Convergence Studies and Stability Tests}
\label{sec::Convergance_Studies}
We now wish to compare the properties of approximations generated using the new multi-kernel approach to those of the traditional single kernel implementations. For all the numerical methods and geometries used herein we distribute the domain with an unstructured point cloud in the same way. Node distributions are generated as in \cite{king_2022}, using propagating front algorithm of \cite{Fornberg&Flyer_disc} before performing ten iterations of a shifting algorithm to avoid clustering of particles. All nodal distributions are on average uniformly spaced ($s_i=s$ is constant) throughout the domain.

\subsection{Convergence}
\label{subsec::convergence}

We now demonstrate the power of using linear combinations of weights to improve meshless approximations to operators. We discretise a doubly-periodic square domain $(x,y) \in [0,1]\times [0,1]$ and follow \cite{king_2022} in defining a test function
\begin{equation}
    \label{eq::test_function}
    \phi(x,y)=\sin(2\pi y)\frac{4}{\pi}\sum_{k=1}^8\frac{\sin\left(2(2k-1)\pi(x-1/4)\right)}{2k-1},
\end{equation}
which has sharp gradients at $x=0.25, 0.75$ and is the first eight terms of the Fourier series of a top-hat function. To assess the convergence of the numerical operator we use the $L_2$ norm defined through
\begin{equation}
    \label{eq::l2_norm_defn}
    L_2-\mathrm{norm}(\cdot )= \frac{\left(\sum_{i=1}^N \left[L_{\mathrm{a}}(\cdot)^d_i-L_{n}(\cdot)_i \right]^2\right)^{1/2}}{\left(\sum_{i=1}^N \left[L_{\mathrm{a}}(\cdot)_i\right]^2\right)^{1/2}}
\end{equation}
where $L_a(\cdot)_i$ is the analytic operator in question, and $L_n(\cdot)_i$ is the numerical approximation. To measure the improvement provided by using linear combinations of weights as detailed in \S\ref{subsec::Optimisation_of_coefficients} we consider the ratio $\mathcal{R}$ between the $L_2-$norm of the solution produced by standard (SK) implementations of meshless methods to the new (MK) approaches outlined herein
\begin{equation}
    \label{eq::l2_norm_ratio}
    \mathcal{R}\coloneqq \frac{L_2-\mathrm{norm}(\mathrm{MK})}{L_2-\mathrm{norm}(\mathrm{SK})}.
\end{equation}

To demonstrate the potential of the new approach outlined herein for this convergence study we exploit the fact that the maximum magnitude of the wavenumber of $\phi$ is $||\boldsymbol{k}||=\pi\sqrt{904}$. Defining $k_M$ to be this maximal magnitude of wavenumber, in this case we limit the optimisation procedure to only consider wavenumbers such that $||\boldsymbol{k}||\leq k_M$.

\begin{figure}
    \begin{center}
    \setlength{\unitlength}{1cm}
    \begin{picture}(18,5)(0,0)
    \put(0,0){\includegraphics[width=0.49\linewidth]{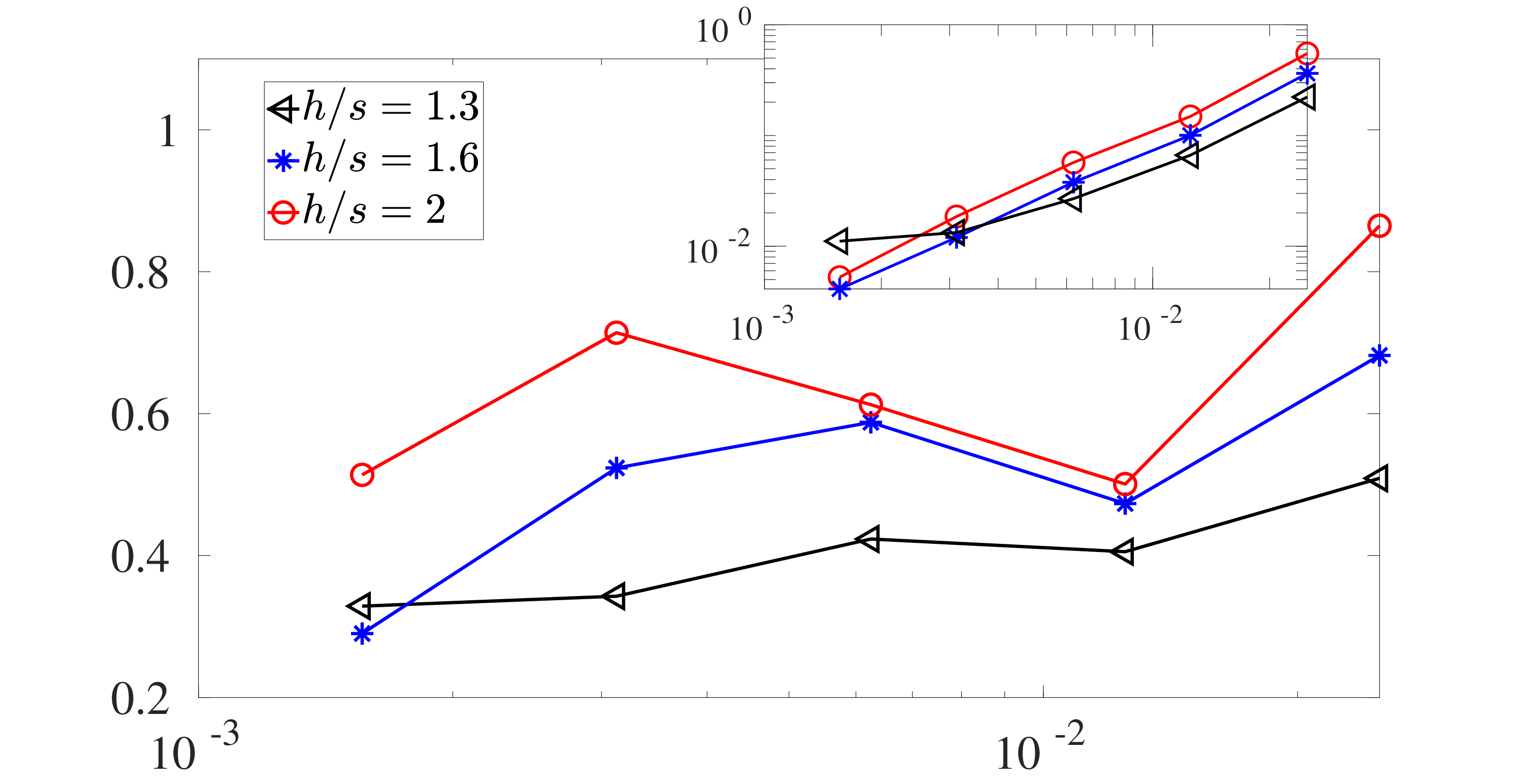}}
    \put(9,0){\includegraphics[width=0.49\linewidth]{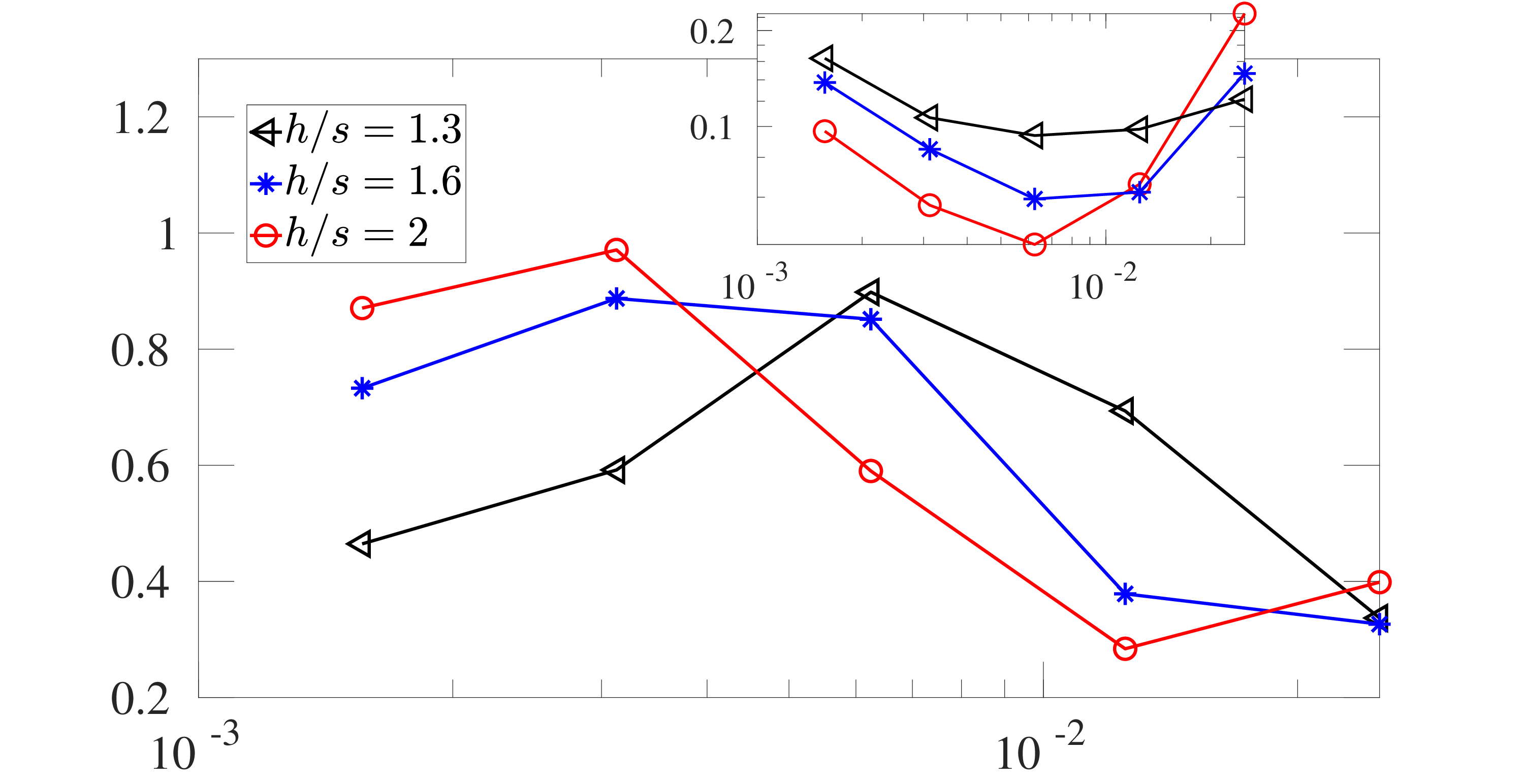}}
    \put(-0.1,2.25){\rotatebox{90}{$\mathcal{R}$}}
   \put(8.9,2.25){\rotatebox{90}{$\mathcal{R}$}}
    \put(4.5,-0.25){$s$}
    \put(13.5,-0.25){$s$}
\end{picture}
    \end{center}
    \caption{Ratio $\mathcal{R}$ of errors arising from SK and MK approach of \S\ref{subsec::Optimisation_of_coefficients} for SPH. Gradient (left) and Laplacian (right) approximations of \eqref{eq::test_function}. Red $+$, -- $h/s=2.0$; blue $*$, $--$ $h/s=1.6$; black $\triangleleft$, $\cdot-\cdot$ $h/s=1.3$. Insets show convergence of $L_2-$ norm of new approximations.}
    \label{fig:grad_conv_sph}
    \end{figure}
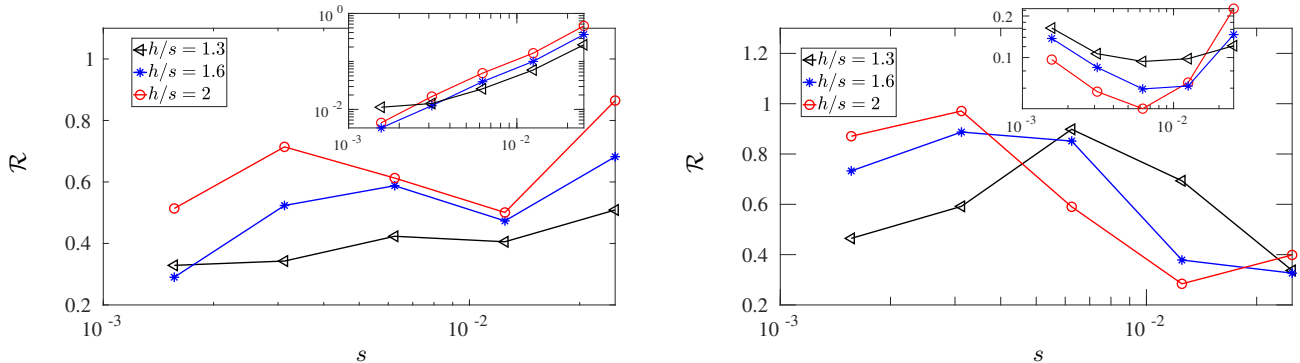

Fig. \ref{fig:grad_conv_sph} shows the variation with resolution of the ratio of errors $\mathcal{R}$ defined in \eqref{eq::l2_norm_ratio} of approximations to the gradient and Laplace operators of \eqref{eq::test_function}. When using SPH approximations, we use a stencil size of $2h$, where $h$ is the kernel smoothing length \cite{Li_Review}. For SPH approximations, the optimisation of resolving power consistently reduces error over all resolutions for the gradient operator, with an improvement in accuracy up to $70\%$ ($\mathcal{R}=0.3$). In particular, the smaller stencils are most improved by the new approach, and we note for readers unfamiliar with SPH that simulations are typically run with a stencil size of $2 h=2\times(1.3s)$. The accuracy of the Laplacian operator is also improved by maximising resolving power using two kernels, with the largest decrease seen on the coarsest discretisations, though the trend is less clear. 

RBF-FD approximations behave in a somewhat different manner. For both gradient and Laplacian operators comparatively little improvement is seen on the coarsest discretisation, whereas an order of magnitude reduction in error ($\mathcal{R}<0.1$) is obtained on highly resolved nodal distributions. Relatively significant improvements are still observed at moderate resolutions so the new MK approach is not limited to improving the accuracy of solutions of well-resolved problems.
\begin{figure}
    \begin{center}
    \setlength{\unitlength}{1cm}
    \begin{picture}(18,5)(0,0)
    \put(0,0){\includegraphics[width=0.49\linewidth]{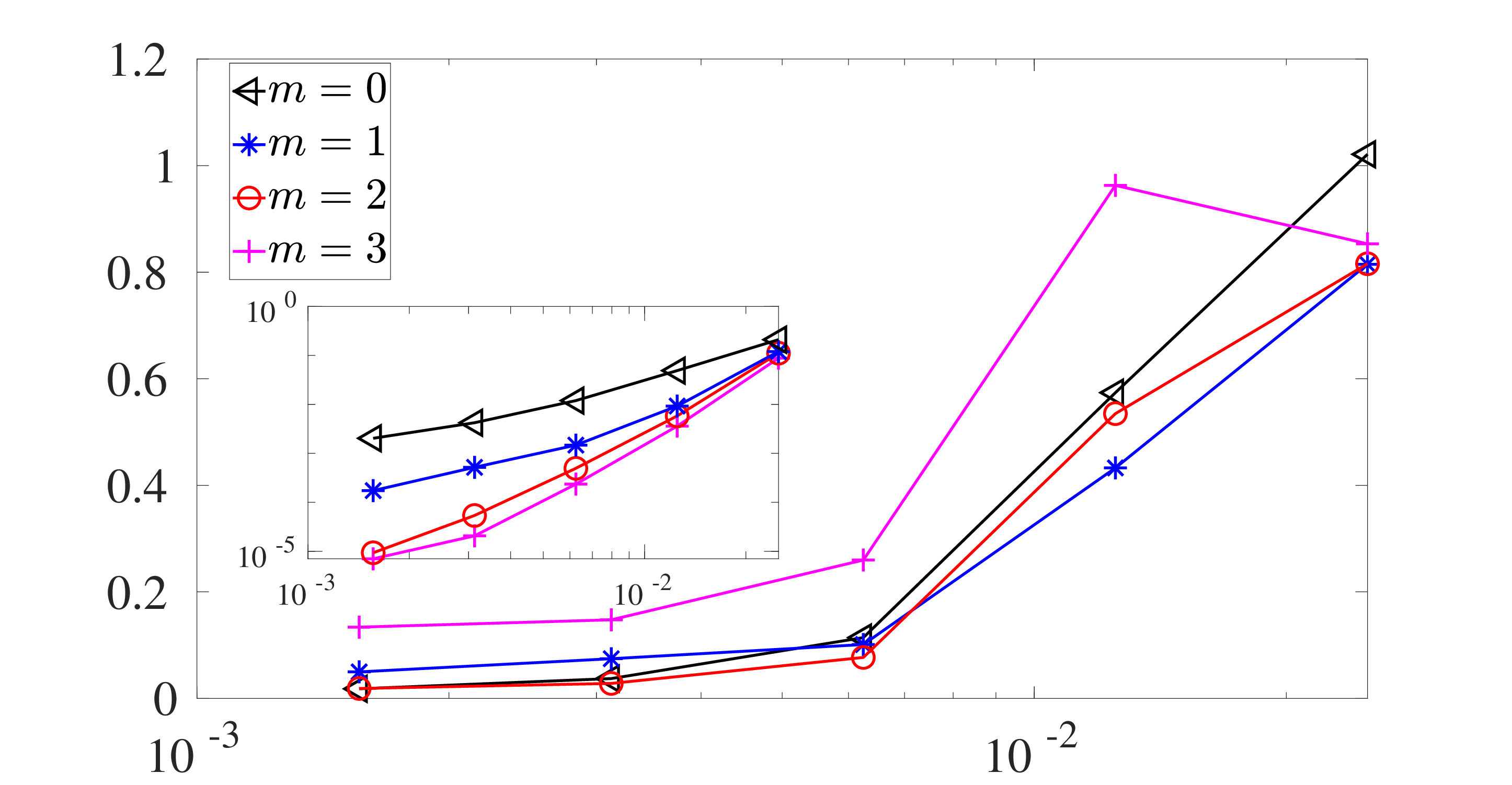}}
    \put(9,0){\includegraphics[width=0.49\linewidth]{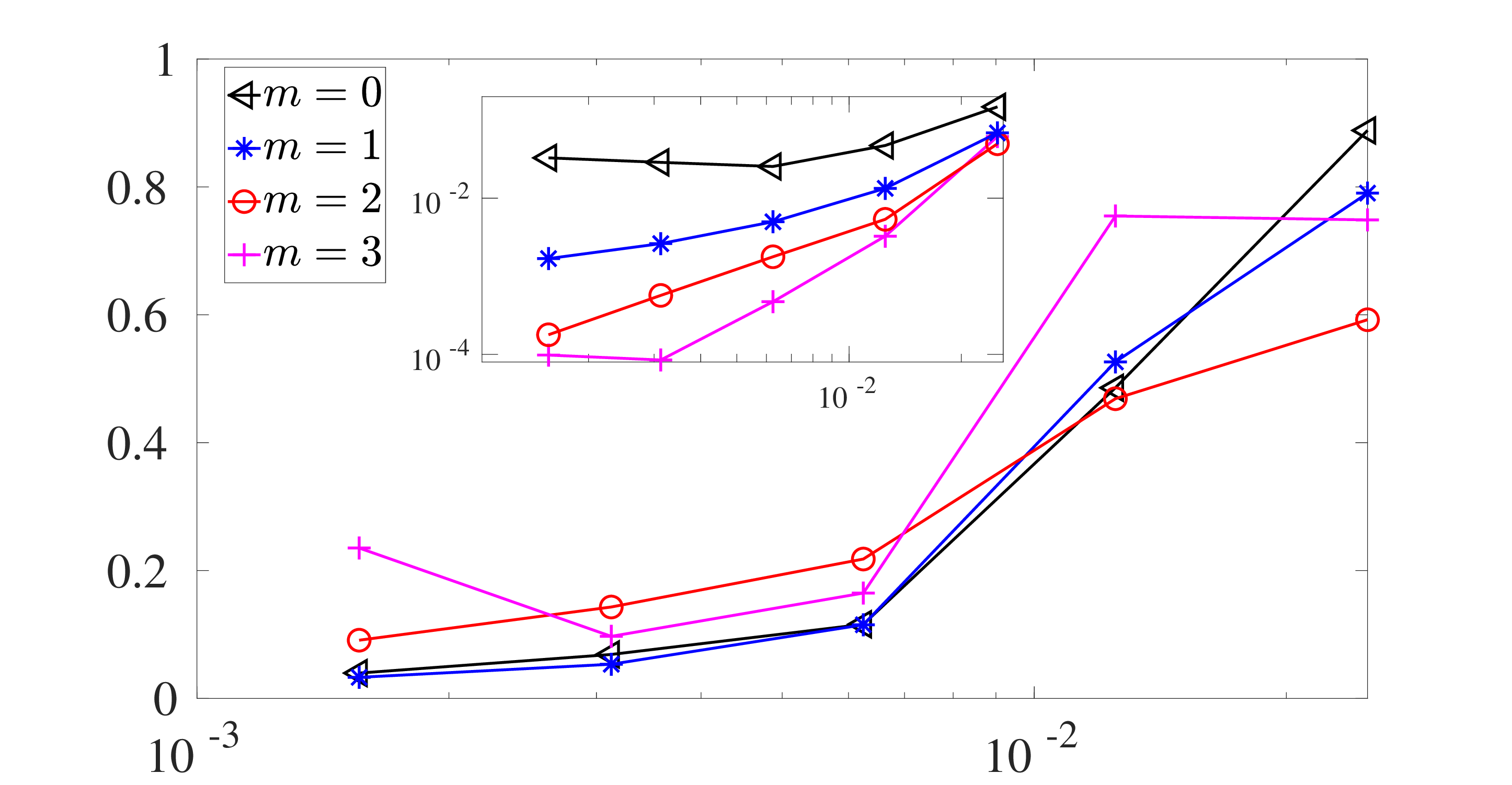}}
    \put(-0.1,2.25){\rotatebox{90}{$\mathcal{R}$}}
    \put(8.9,2.25){\rotatebox{90}{$\mathcal{R}$}}
    \put(4.5,-0.25){$s$}
    \put(13.5,-0.25){$s$}
\end{picture}
    \end{center}
    \caption{Ratio $\mathcal{R}$ of errors arising from standard implementation and new approach of \S\ref{subsec::Optimisation_of_coefficients} for RBF-FDs. Gradient (left) and Laplacian (right) approximations of \eqref{eq::test_function}. Black $\triangleleft$, $\cdot-\cdot$ $m=0$, $\lvert\mathcal{N}_i\rvert=10$; pink $+$, $\cdots$ $m=1$, $\lvert\mathcal{N}_i\rvert=15$; blue $*$, $--$ $m=2$, $\lvert\mathcal{N}_i\rvert=20$; red $\circ$ , --$m=3$, $\lvert\mathcal{N}_i\rvert=25$. Insets show convergence of $L_2-$ norm of new approximations.}
    \label{fig:grad_conv_rbffd}
    \end{figure}
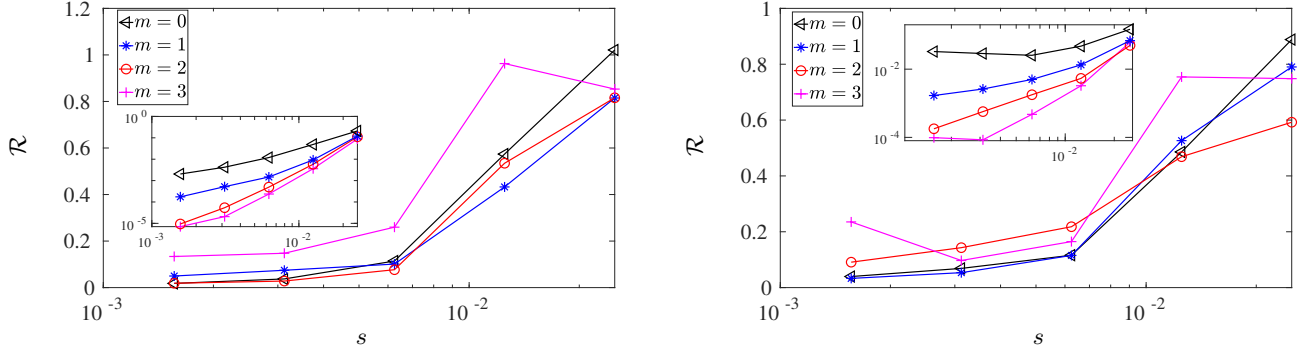

Approximations using LABFM are considered at higher order than RBF-FDs. Significant improvements in gradient approximations are observed at all orders, though applying the multi-basis formulation to Laplacian operators results in smaller gains. For gradient operators, the relative improvement is most significant for the $m=8$ scheme, though in all cases the improvement remains relatively constant ($20-50\%$ depending on $m$) for sufficiently refined meshes. For the Laplacian operator, although at some resolutions only an error reduction of $\sim10\%$ is attained, we stress that for a time-dependent problem no extra computational steps are required to calculate derivatives thus even this small gain is achieved for little to no cost.
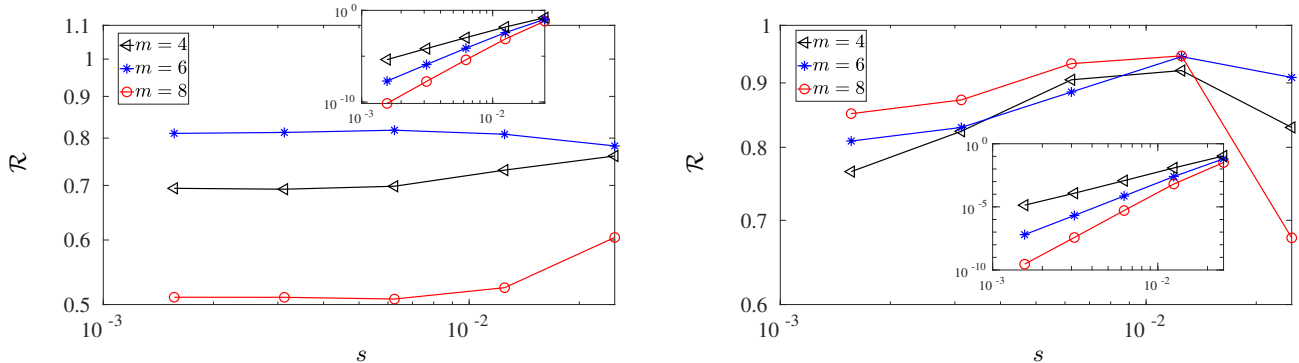
\begin{figure}
    \begin{center}
    \setlength{\unitlength}{1cm}
    \begin{picture}(18,5)(0,0)
    \put(0,0){\includegraphics[width=0.49\linewidth]{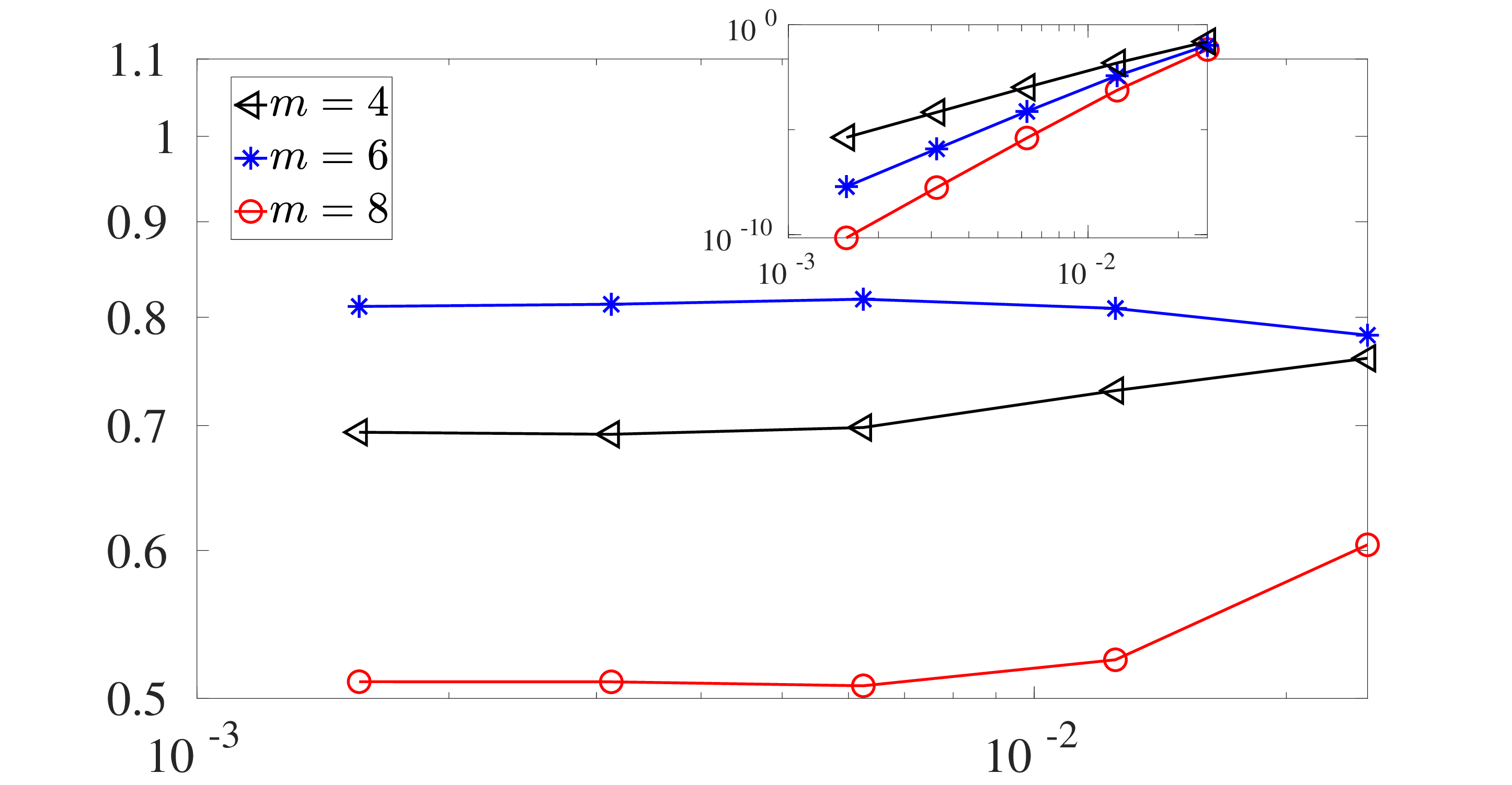}}
    \put(9,0){\includegraphics[width=0.49\linewidth]{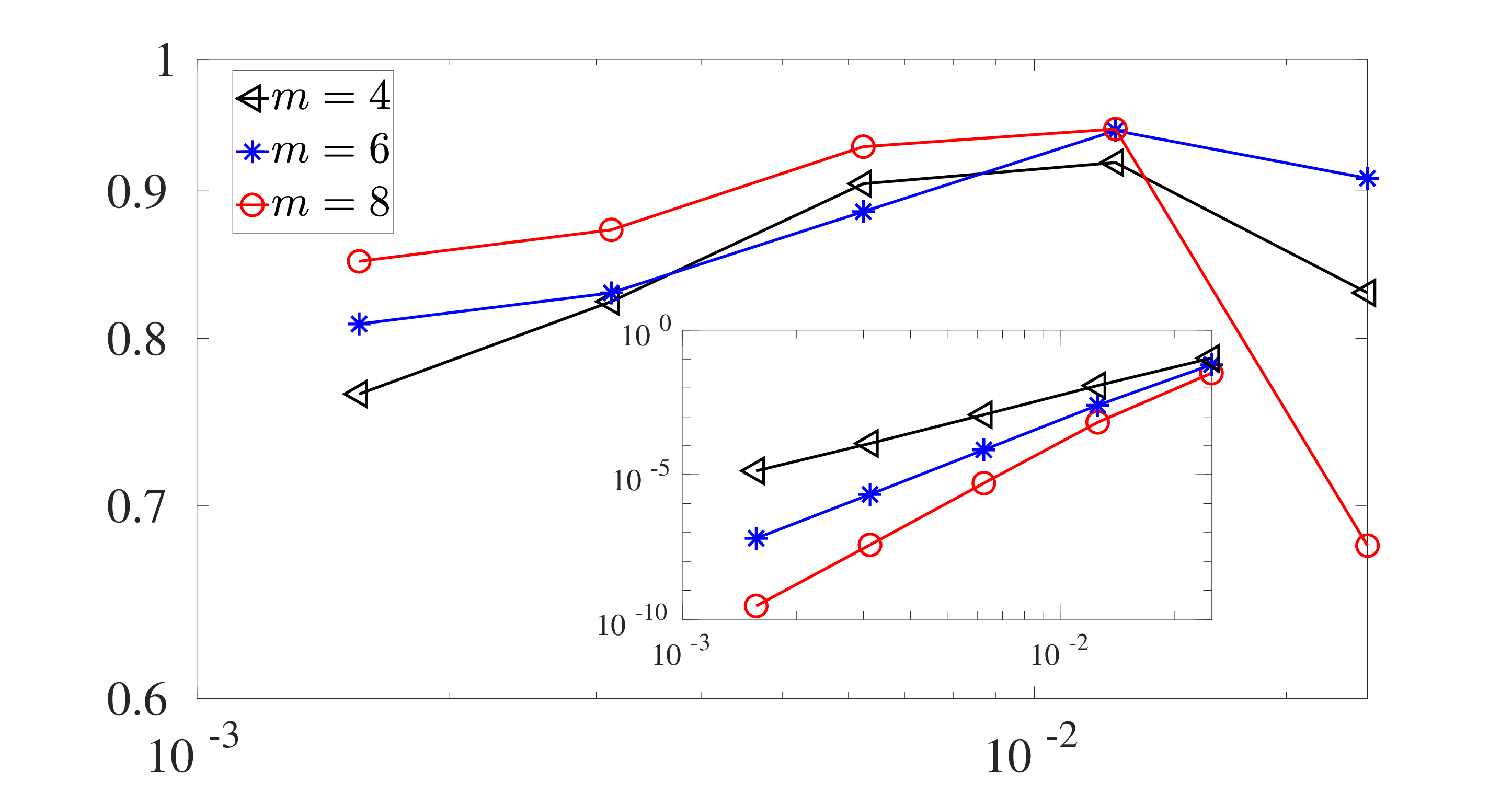}}
    \put(-0.1,2.25){\rotatebox{90}{$\mathcal{R}$}}
    \put(8.9,2.25){\rotatebox{90}{$\mathcal{R}$}}
    \put(4.5,-0.25){$s$}
    \put(13.5,-0.25){$s$}
\end{picture}
    \end{center}
    \caption{Ratio $\mathcal{R}$ of errors arising from standard implementation and new approach of \S\ref{subsec::Optimisation_of_coefficients} for LABFM. Gradient (left) and Laplacian (right) approximations of \eqref{eq::test_function}. Black $\triangleleft$ $m=4$; blue $*$ $m=6$; red $\circ$ $m=8$. Insets show convergence of $L_2-$ norm of new approximations.}
    \label{grad_conv_LABFM}
    \end{figure}

\subsection{Resolving Power}
\label{subsec::RP_results}
Now that clear improvements in approximations have been achieved using the new framework, we discuss the improvements in resolving power. As mentioned in \S\ref{sec::resolving_power_formulation}, resolving power has typically only been considered in one-dimension even for meshless methods. However, it has been shown herein that for (2D) meshless methods we must perform a multi-dimensional resolving power analysis. We now demonstrate the effect of using linear combinations of weights on the resolving power of each method. For brevity, in this subsection we shall consider only one order of polynomial consistency and stencil size per method, matching this to those of Fig. \ref{fig:rp_grad_examples}.

In order to quantify the improvement in resolving power attained using the optimisation procedure of \S\ref{subsec::Optimisation_of_coefficients} we introduce the percentage error reduction of resolving power of gradient operators $\varepsilon_G$ defined through 
\begin{gather}
    \varepsilon_G(k_x,k_y)=\left(1-\frac{\sqrt{(\Re\{\hat{k}_{eff}^{MK}\}-\hat{k})^2+(\Im\{\hat{k}_{eff}^{MK}\})^2}}{\sqrt{(\Re\{\hat{k}_{eff}^{SK}\}-\hat{k})^2+(\Im\{\hat{k}_{eff}^{SK}\})^2}}\right)\times100 \%
\end{gather}
where $\hat{k}, \hat{k}_{eff}$ are as defined in \eqref{eq:hatk_defn}, and superscripts $MK$ denote effective wavenumbers arising from the multi-kernel approach of \S\ref{sec::linear_combinations_of_weights} and $SK$ the usual single kernel approaches.

\begin{figure}
    \begin{center}
    \setlength{\unitlength}{1cm}
    \begin{picture}(18,5)(0,0)
    \put(0,0){\includegraphics[width=0.49\linewidth]{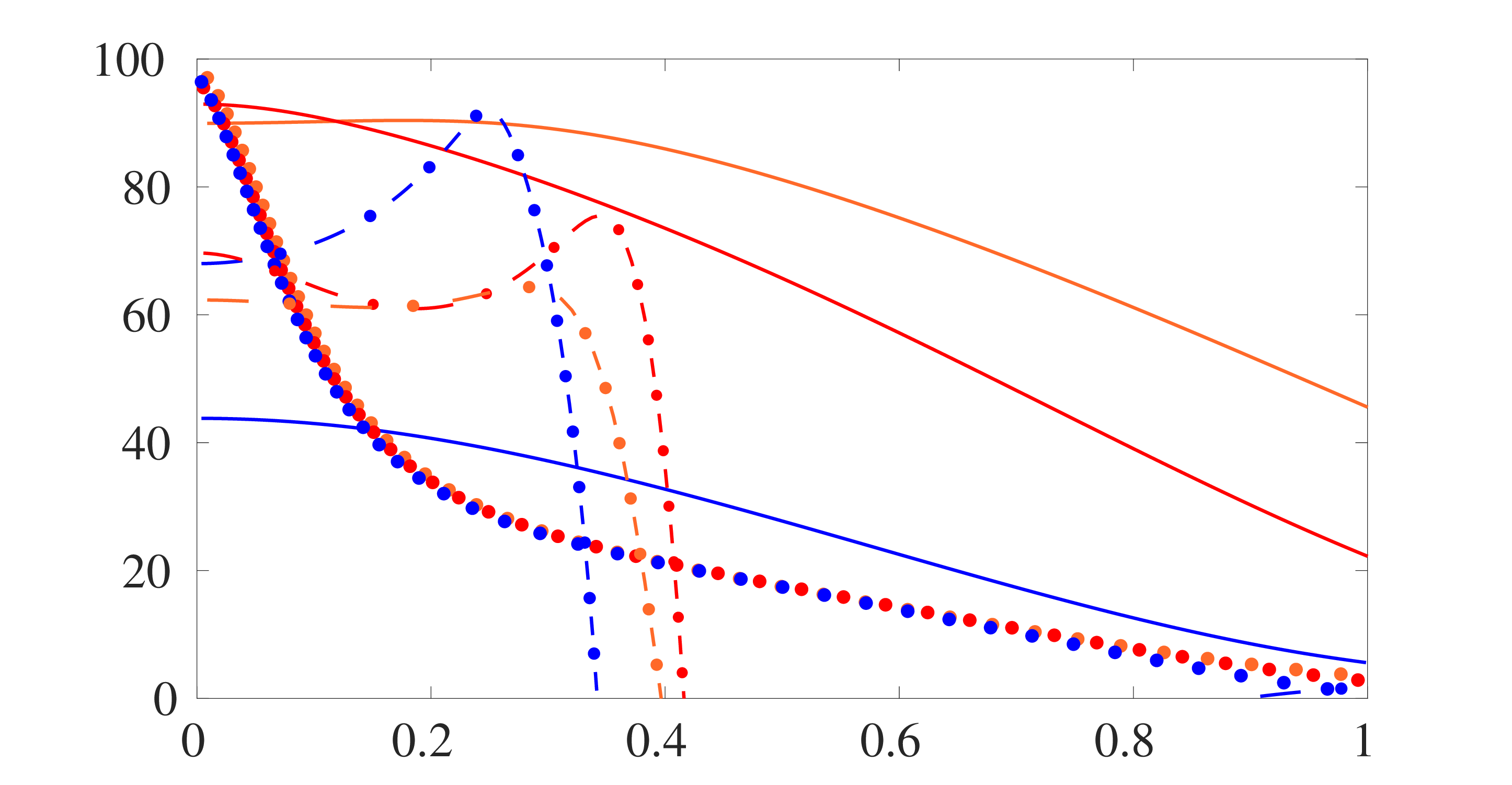}}
    \put(9,0){\includegraphics[width=0.49\linewidth]{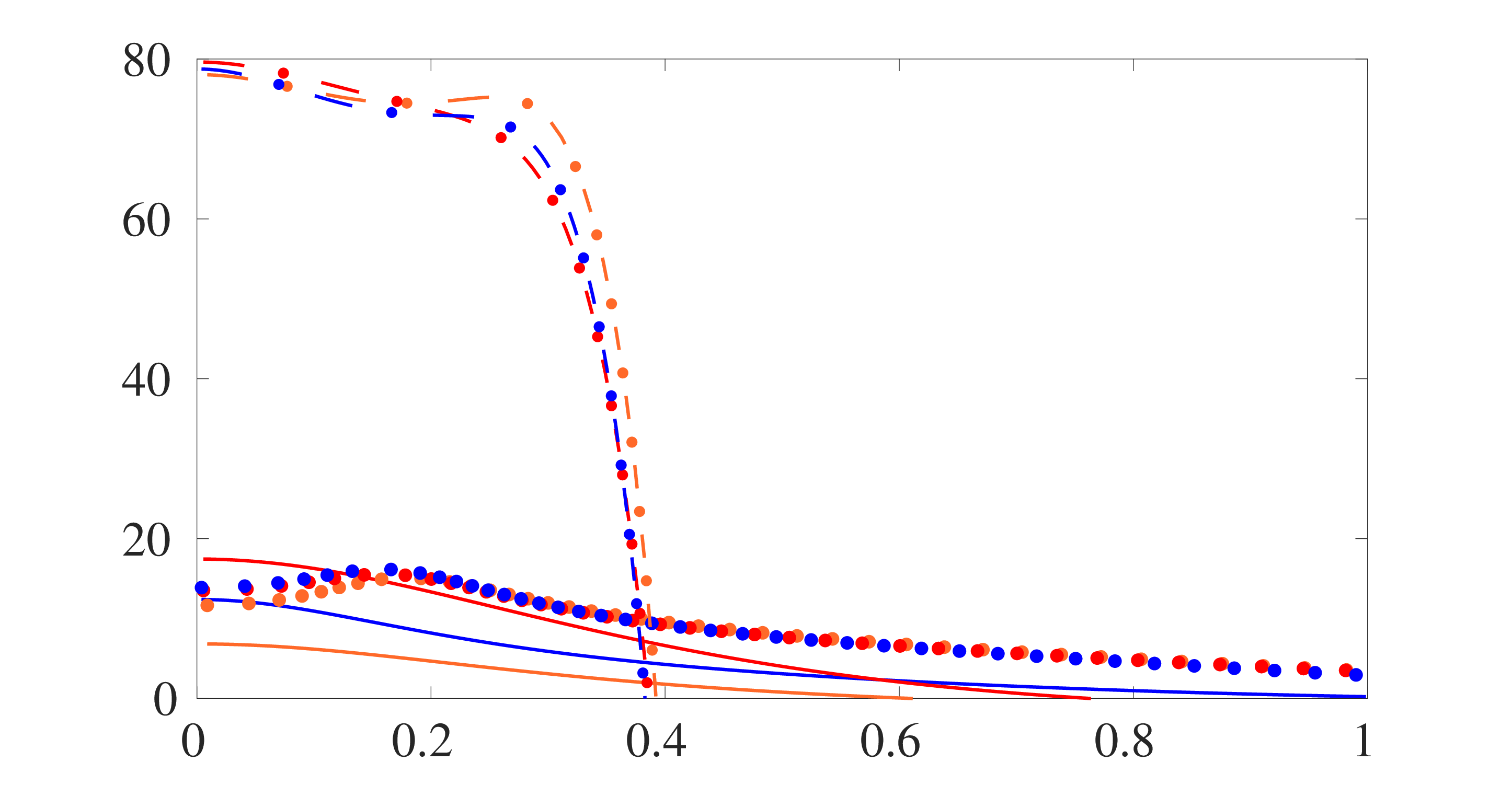}}
 
    \put(13.5,-0.25){$\hat{q}$}
    \put(-0.1,2){\rotatebox{90}{$\varepsilon_G$}}
    \put(8.9,2){\rotatebox{90}{$\varepsilon_L$}}

    \put(4.5,-0.25){$\hat{k}$}

    \end{picture}
    \end{center}
    \caption{Improvements in resolving power for multi-kernel approximations to operators for various meshless methods. Left panel gradient operators, right panel Laplace operators. Solid lines denote LABFM with $m=4$, dotted lines SPH, dot-dashed lines RBF-FDs with $m=2$. Blue $k_y=0$, red $k_y=k_x$, orange $k_y=2k_x$.}
    
    \label{fig:rp_grad_rp_approach}
\end{figure}
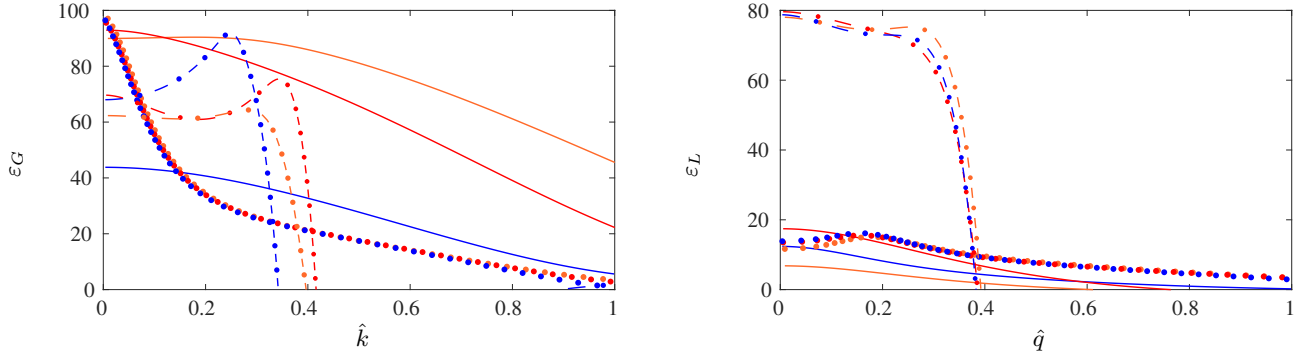

Fig. \ref{fig:rp_grad_rp_approach} (left panel) compares the resolving power of gradient operators using the new approaches to the classical approaches along different lines (all emanating from the origin) in $(k_x,k_y)$ space. For SPH operators (dots) there is consistent improvement of the effective wavenumber across all of $(k_x,k_y)$ space, with $\varepsilon_G$ having little dependence on the line being considered. The improvement is substantial for small wavenumbers, and even up to wavenumbers half the Nyquist wavenumbers there is still a $20\%$ improvement in resolving power.

For LABFM (solid lines), the improvement is dependent on the line in wavenumber space considered. Along $k_y=k_x,2k_x$ up to $90\%$ improvement in resolving power is attained at low wavenumbers and these large improvements persist up to $\hat{k}=1$. Modes with $k_y=0$ see comparatively less improvement at all wavenumbers, however the multi-kernel scheme is still closer to spectral across the whole range of $\hat{k}$.

For RBF-FDs (dash-dot lines) the picture is more complicated. Use of the multi-kernel approach helps to reduce/eliminate regions where $\hat{k}_{eff}>\hat{k}$ which occur at low wavenumbers. However, this generally comes at a cost of improvement at higher wavenumbers. Significant improvements are therefore observed at low wavenumbers (up to $\hat{k}=0.4$) and we note that due to the effects of dissipation/filtering this is where most of the energy of the wavenumber spectrum of a solution to a PDE system will likely sit. In addition shrinking the region in wavenumber space where $\hat{k}_{eff}>\hat{k}$ will stabilise simulations of fast-moving waves for reasons discussed in \S\ref{sec::resolving_power_formulation}.

We now turn our attention to the improvement in the resolving power of approximations to the Laplacian operators as shown in the right panel of Fig. \ref{fig:rp_grad_rp_approach}. We define a parameter $\varepsilon_L$ which measures percentage error reduction in a similar manner to $\varepsilon_G$, replacing real wavenumber $\hat{k}$ with $\hat{q}^2$ and effective wavenumbers $\hat{k}_{eff}$ with $\hat{q}_{eff}^2$.  For SPH, we see a consistent improvement in effective wavenumber with the new approach, however it is relatively small compared to the gradient error with a maximum of $20\%$ error reduction. This is consistent with Fig. \ref{fig:grad_conv_sph} for which smaller improvements were obtained using the multi-kernel approach compared to the gradient operator.

Similarly small reductions are observed for LABFM approximations. Maximum improvements are attained in the limit $\hat{k}\rightarrow 0$, with only $\sim10\%$ improvement in resolving power at $\hat{q}=0.5$.

For RBF-FDs we see similar behaviour to the gradient approximation. The multi-kernel approach has a large impact on error reduction for small wavenumbers, however this sacrifices some accuracy at higher wavenumbers. Although it is desirable to reduce error across the whole wavenumber space, given the special properties of RBF-FDs operators (i.e. larger values of effective wavenumber than real wavenumber) it is more important to reduce error for small wavenumbers, as this enables the large reductions in error (especially for well resolved approximations) seen in Fig. \ref{fig:grad_conv_rbffd}.

These results explain the improvement in approximations to operators seen in the previous subsection. For both SPH and LABFM the gains in resolving power are far greater for gradients compared to the Laplacian, hence the greater reduction in $L_2-$error for the former. The improvement in the resolving power of RBF-FDs for small $\hat{k}$ led to a near order of magnitude improvement in approximations on well refined discretisations.

\subsection{Stability Tests}
\label{subsec::stability}
Many solutions of systems of PDEs are time-dependent. It is therefore crucial that numerical approximations of derivatives are as stable as possible. To analyse the stability of each operator using a given method we construct a global derivative operator $\boldsymbol{A}^d$ following \cite{king_2020}, such that $L\left(\bm{\phi}\right)=\bm{A}\bm{\phi}$, where $\bm{\phi}$ is a vector containing all $\phi_{i}$. The eigenvalues $\boldsymbol{\lambda}$ of this matrix determine the stability of the operator. For gradient operators, we desire $\Re\{\lambda_i\}=0,\ \forall i$ (no growth or decay of modes in time), however for stability of purely convective problems we only require $\Re\{\lambda_i\}\leq0,\ \forall i$.

\begin{figure}[t]
    \begin{center}
    \setlength{\unitlength}{1cm}
    \begin{picture}(18,3)(0,0)
    \put(0,0){\includegraphics[width=0.32\linewidth]{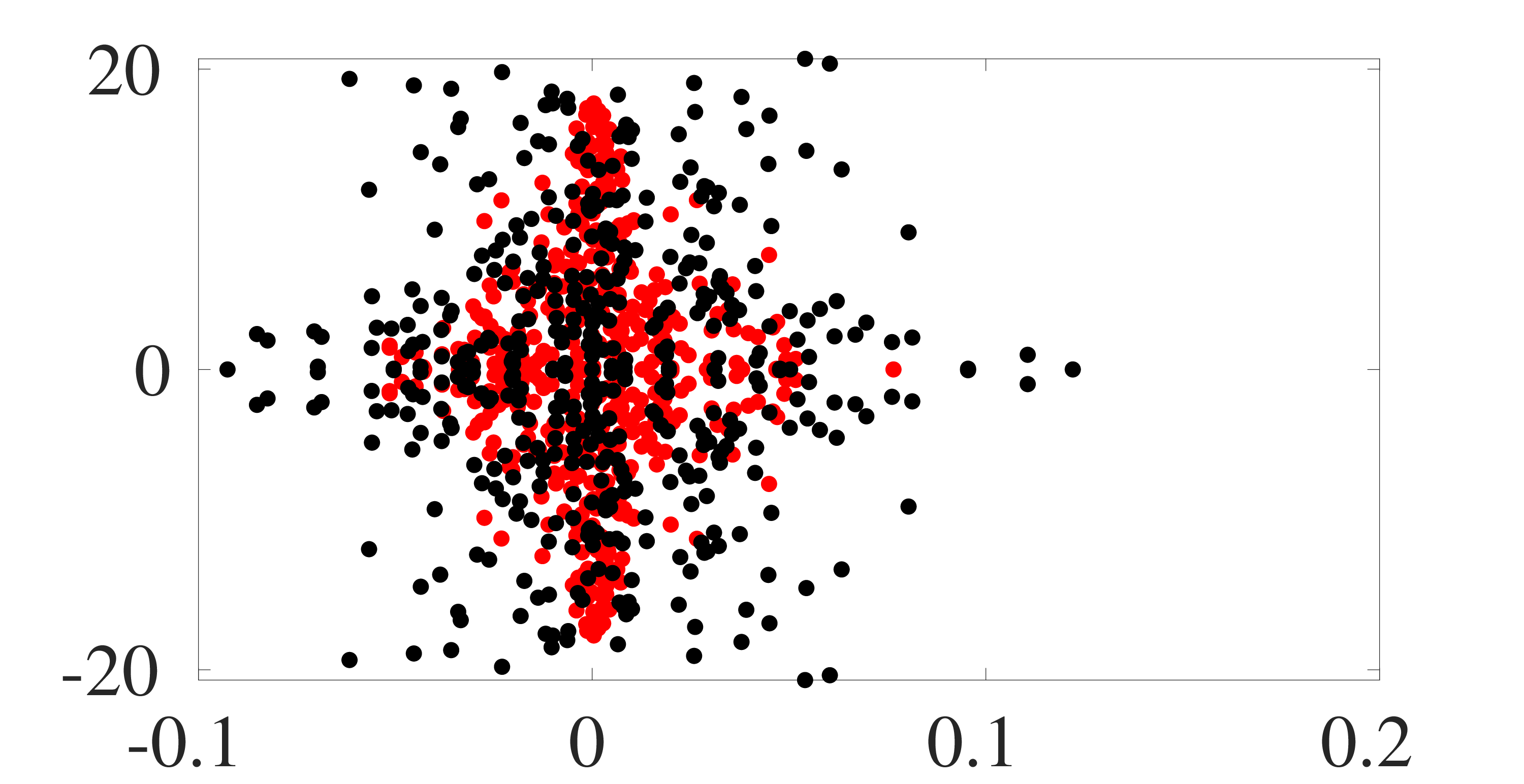}}
    \put(6,0){\includegraphics[width=0.32\linewidth]{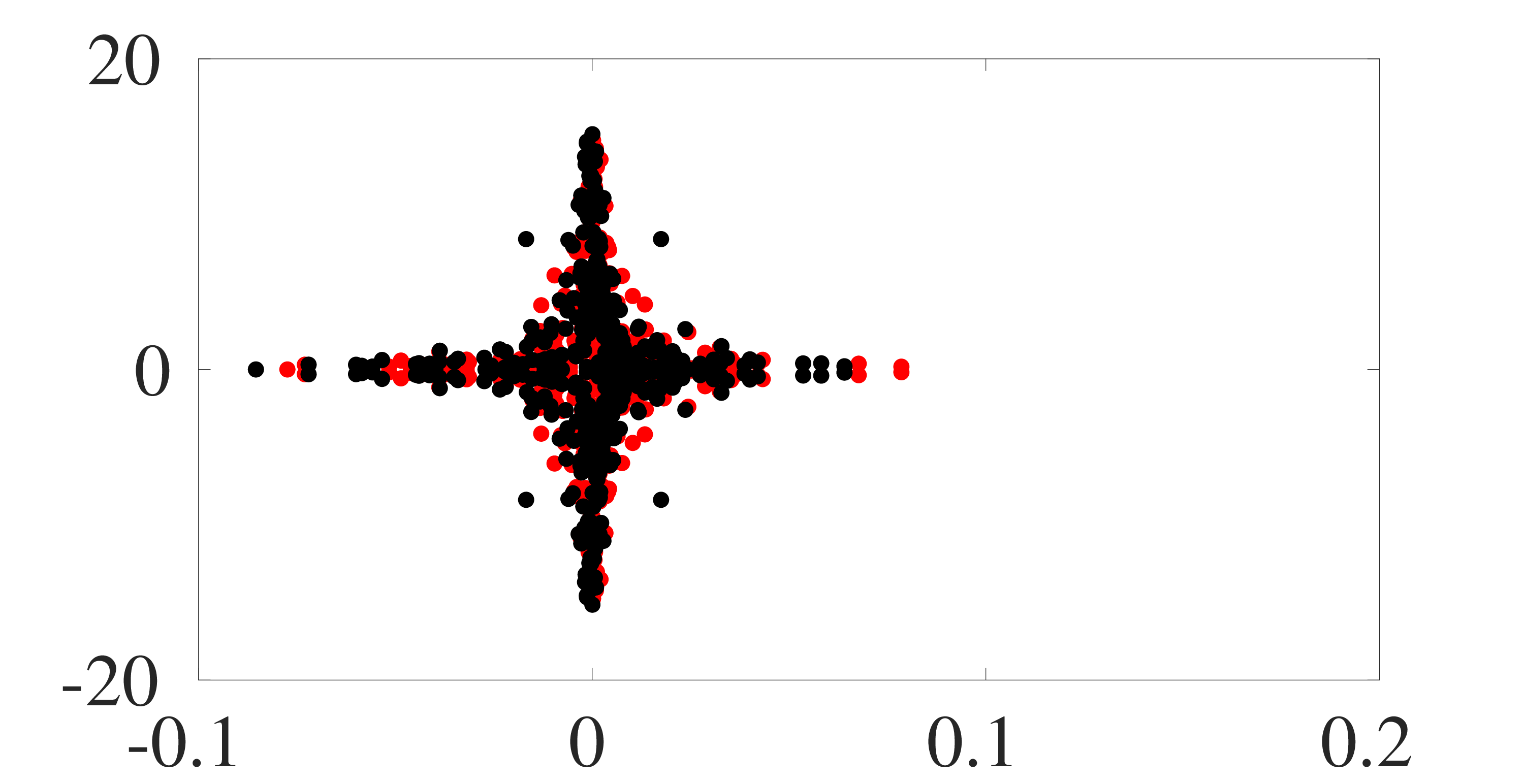}}
    \put(12,0){\includegraphics[width=0.32\linewidth]{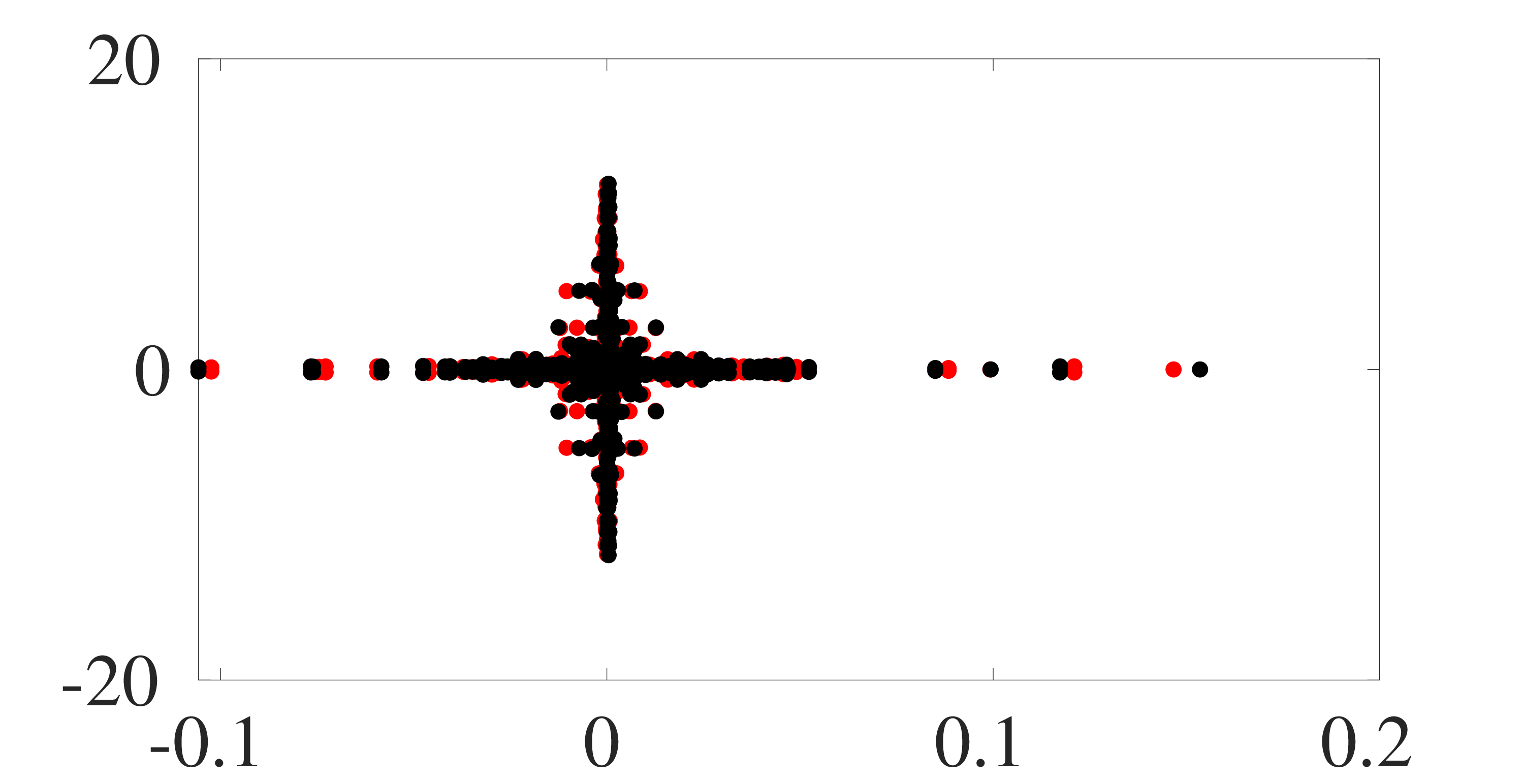}}
    \put(2.5,-0.45){$\Re\{\lambda\}$}
    \put(8.5,-0.45){$\Re\{\lambda\}$}
    \put(14.5,-0.45){$\Re\{\lambda\}$}
     \put(-0.25,1.15){\rotatebox{90}{$\Im\{\lambda\}$}}
     \put(5.75,1.15){\rotatebox{90}{$\Im\{\lambda\}$}}
     \put(11.75,1.15){\rotatebox{90}{$\Im\{\lambda\}$}}
    \end{picture}
    \end{center}
    \caption{Stability of gradient operator for SPH approximations on a typical grid with $441$ points. Red dots correspond to Wendland C2 kernel, black dots to those produced using a combination of two kernels as outlined in \S\ref{subsec::Optimisation_of_coefficients}. Left panel: $h/s=1.3$; Centre Panel: $h/s=1.6$; Right Panel: $h/s=2.0$.}
    \label{fig:sph_grad_stab}
\end{figure}
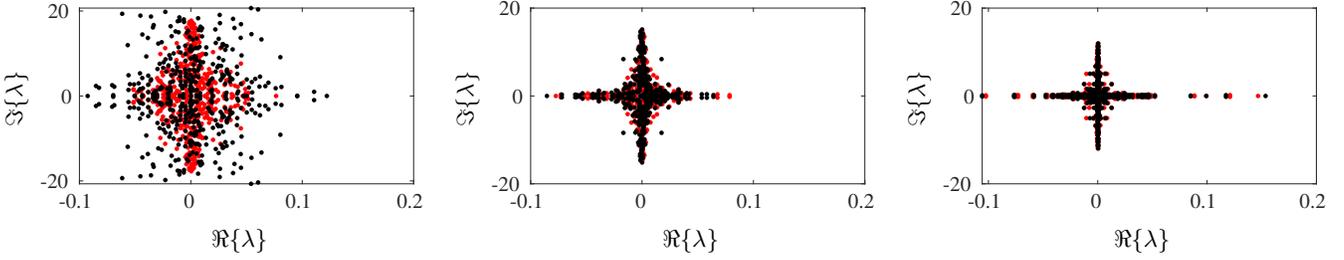

\begin{figure}[t]
    \begin{center}
    \setlength{\unitlength}{1cm}
    \begin{picture}(18,2)(0,0)
    \put(0,0){\includegraphics[width=0.24\linewidth]{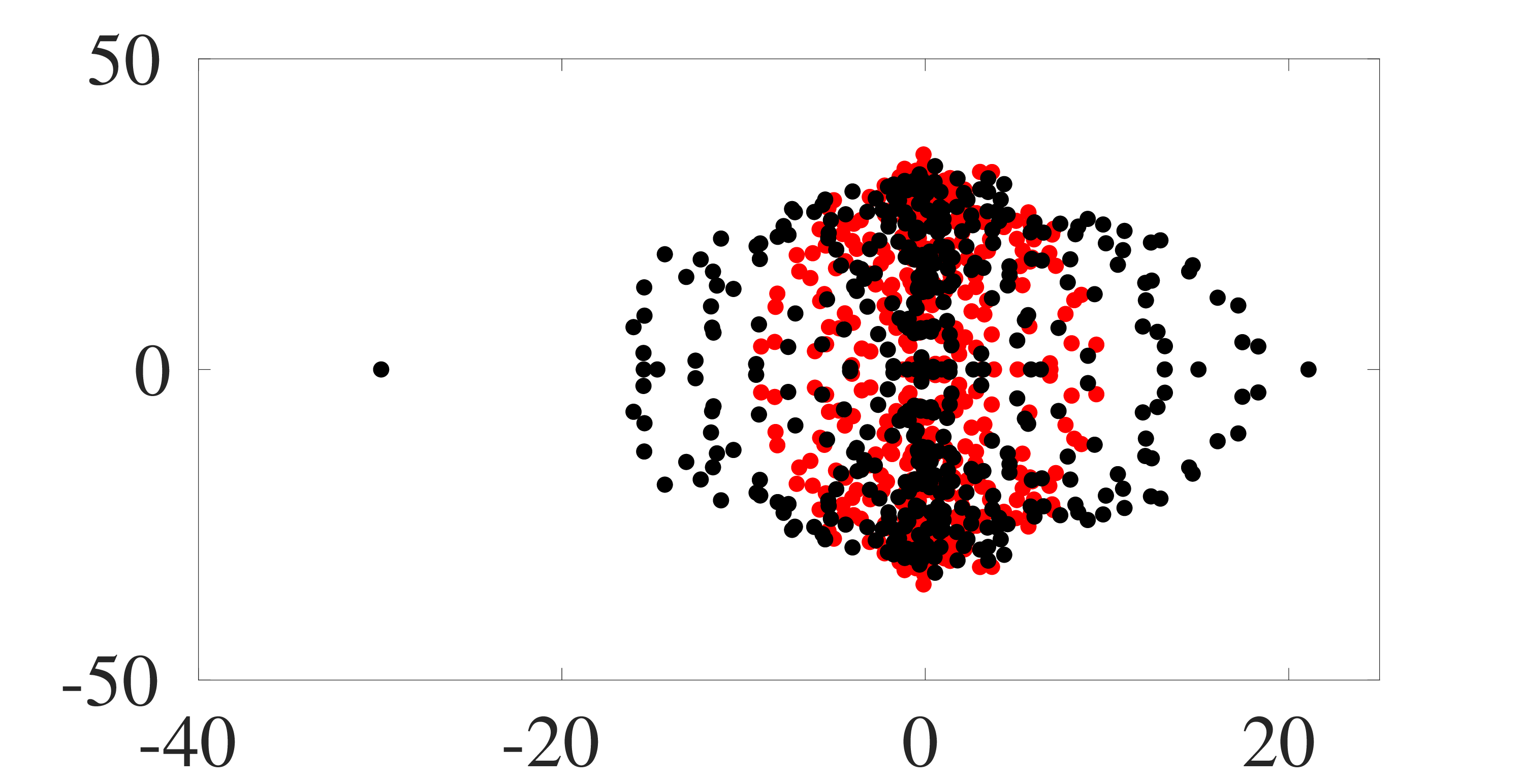}}
    \put(4.5,0){\includegraphics[width=0.24\linewidth]{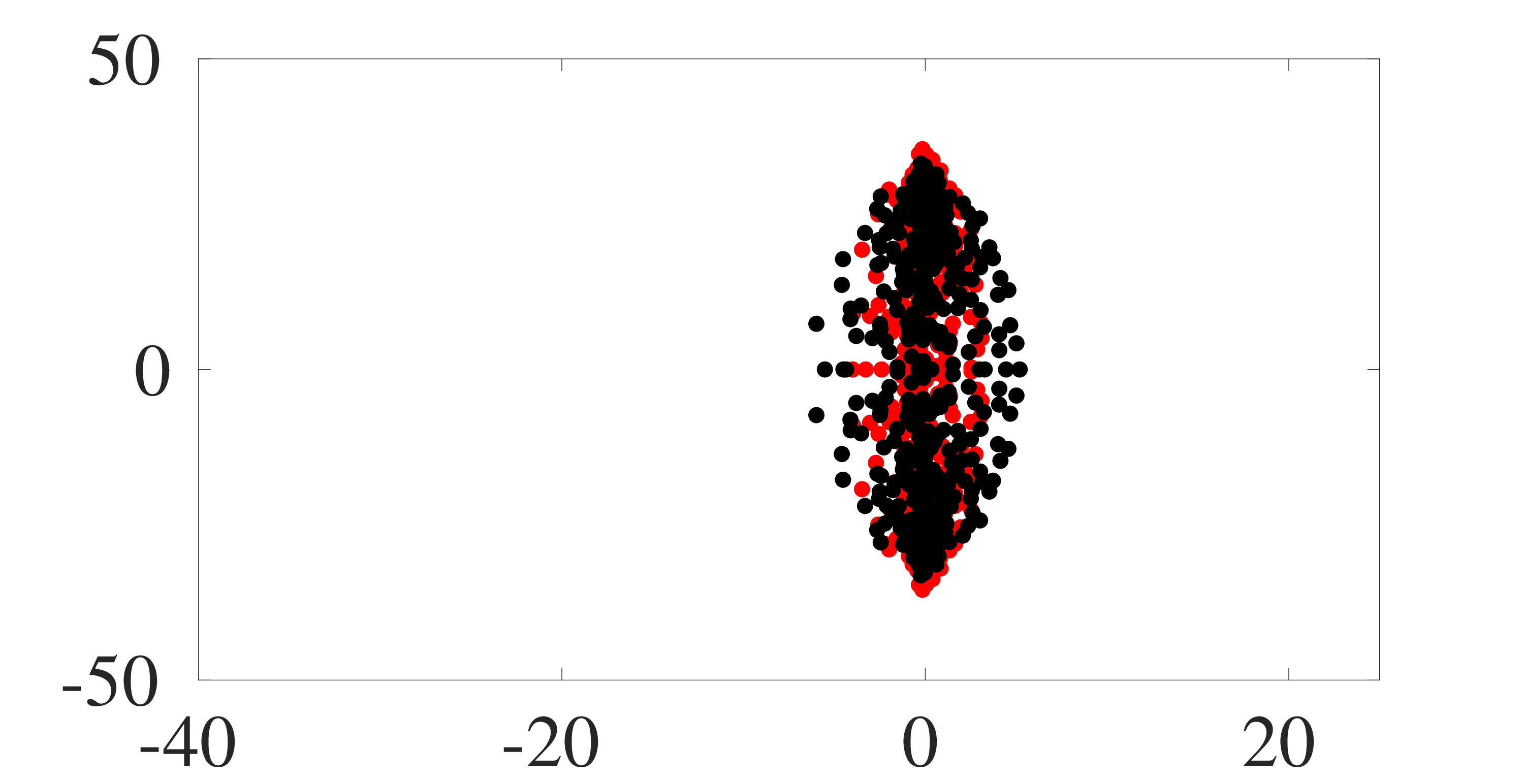}}
    \put(9,0){\includegraphics[width=0.24\linewidth]{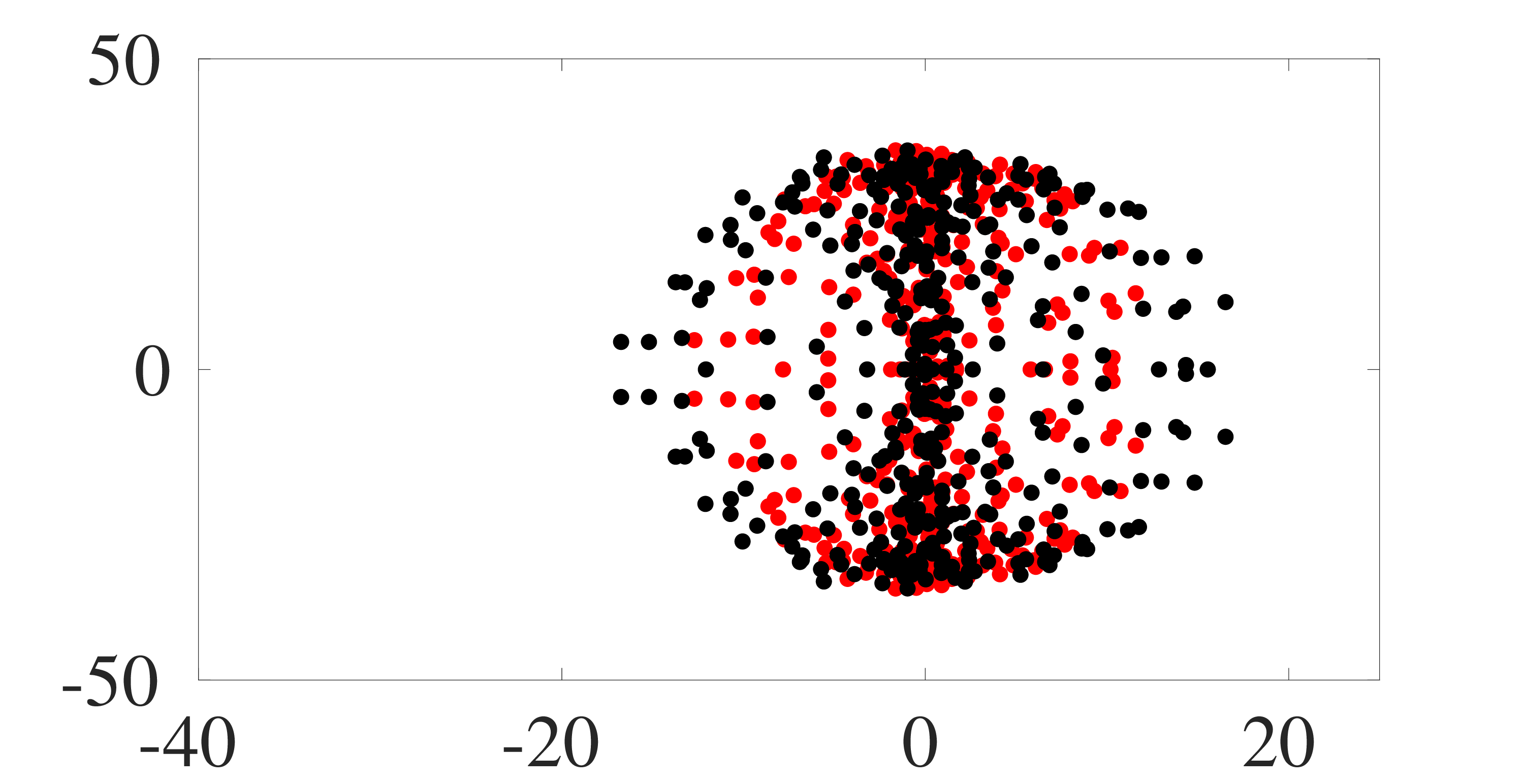}}
    \put(13.5,0){\includegraphics[width=0.24\linewidth]{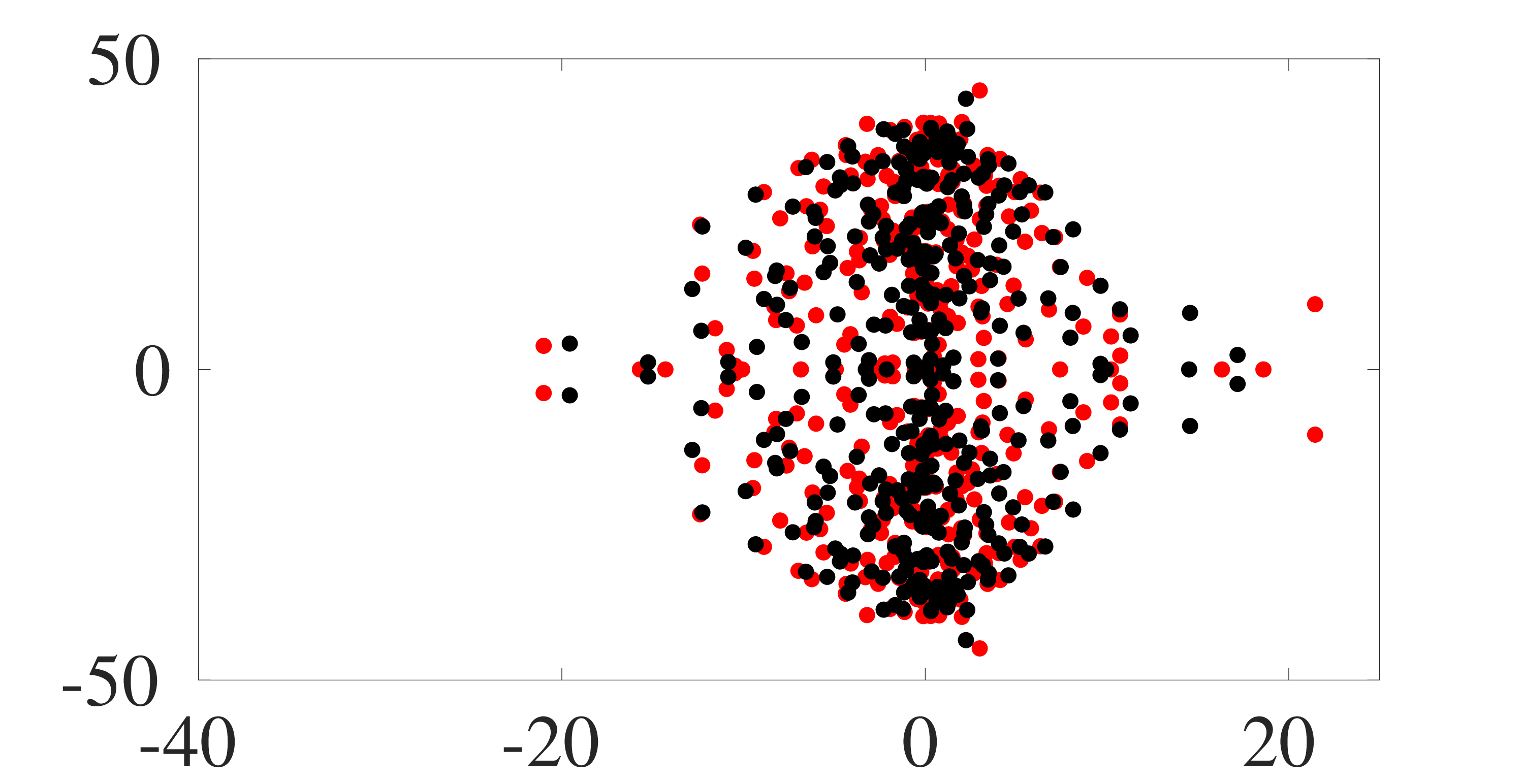}}
    \put(1.9,-0.45){$\Re\{\lambda\}$}
    \put(6.25,-0.45){$\Re\{\lambda\}$}
    \put(10.9,-0.45){$\Re\{\lambda\}$}
     \put(15.25,-0.45){$\Re\{\lambda\}$}
     \put(-0.2,0.8){\rotatebox{90}{$\Im\{\lambda\}$}}
     \put(4.3,0.8){\rotatebox{90}{$\Im\{\lambda\}$}}
     \put(8.8,0.8){\rotatebox{90}{$\Im\{\lambda\}$}}
     \put(13.2,0.8){\rotatebox{90}{$\Im\{\lambda\}$}}
    
    \end{picture}
    \end{center}
    \caption{Stability of gradient operator for RBF-FD approximations on a typical grid with $441$ pointsRed dots correspond to Gaussian RBF, black dots to those produced by linear combination of RBFs as detailed in \S\ref{subsec::Optimisation_of_coefficients}. From left to right: $m=0$, $\lvert\mathcal{N}_i\rvert=10$; $m=1$, $\lvert\mathcal{N}_i\rvert$=15; $m=2$, $\lvert\mathcal{N}_i\rvert=20$; $m=3$, $\lvert\mathcal{N}_i\rvert=25$.}
    \label{fig:rbf_grad_stab}
\end{figure}

\begin{figure}[t]
    \begin{center}
    \setlength{\unitlength}{1cm}
   \begin{picture}(18,3)(0,0)
    \put(0,0){\includegraphics[width=0.32\linewidth]{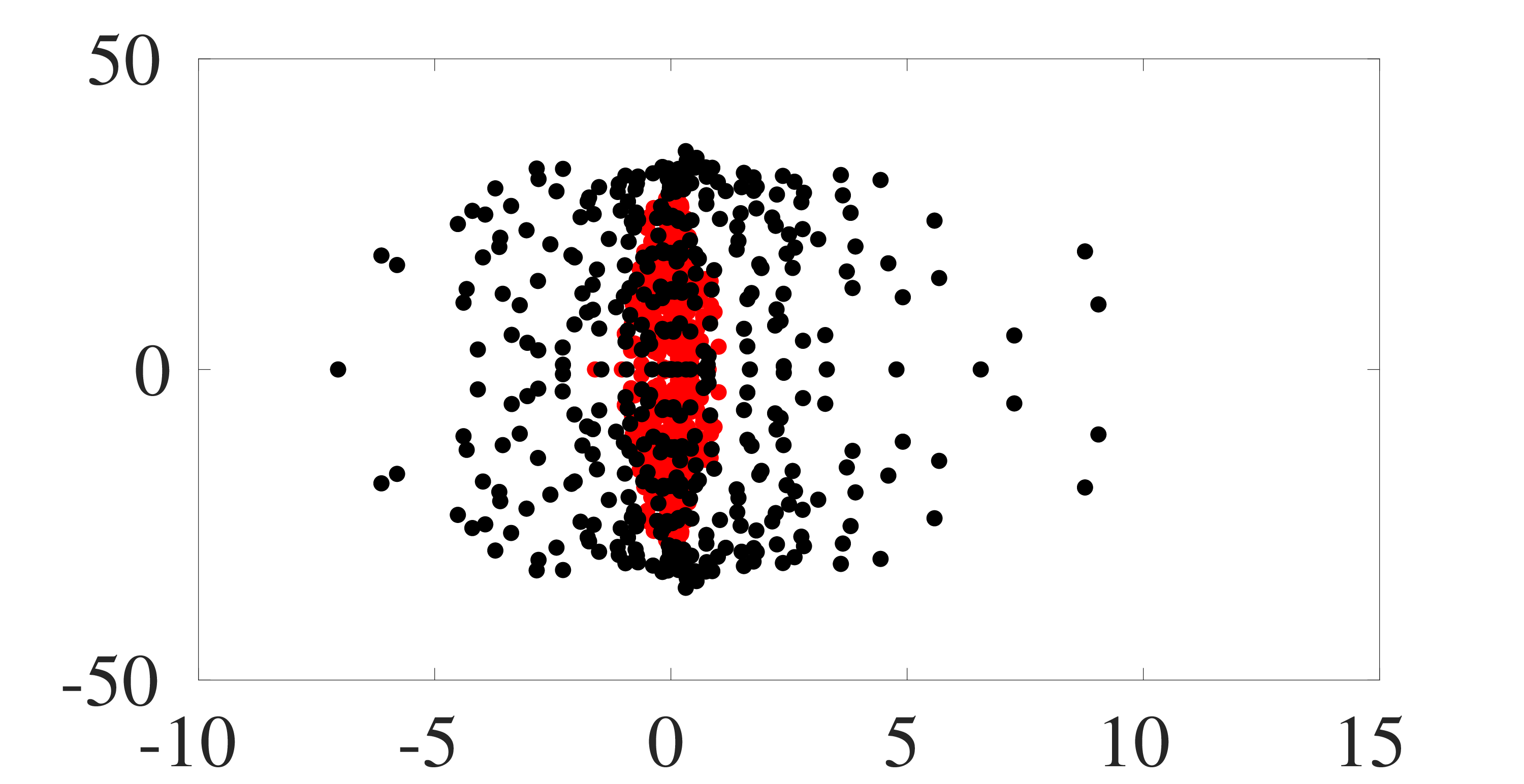}}
    \put(6,0){\includegraphics[width=0.32\linewidth]{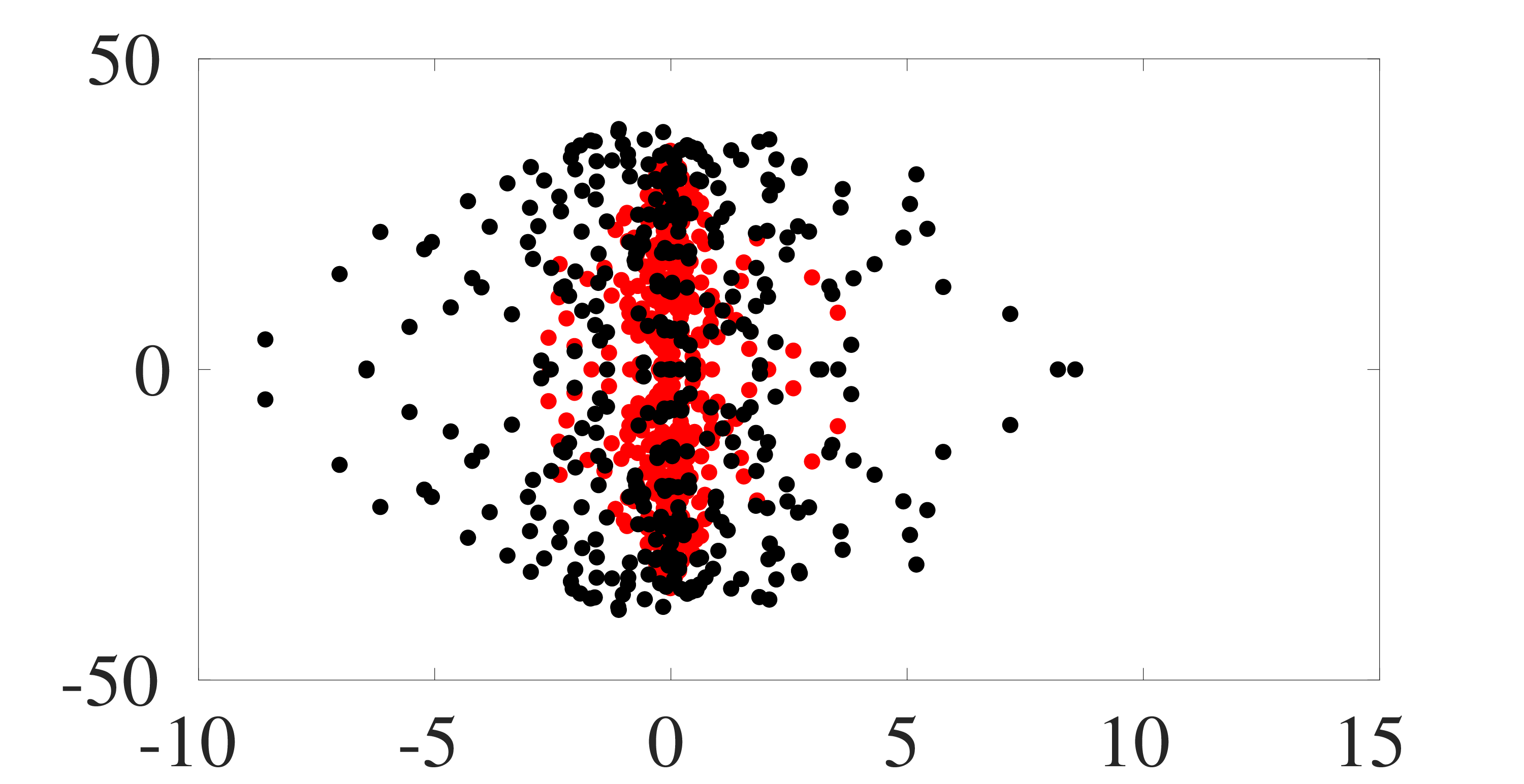}}
    \put(12,0){\includegraphics[width=0.32\linewidth]{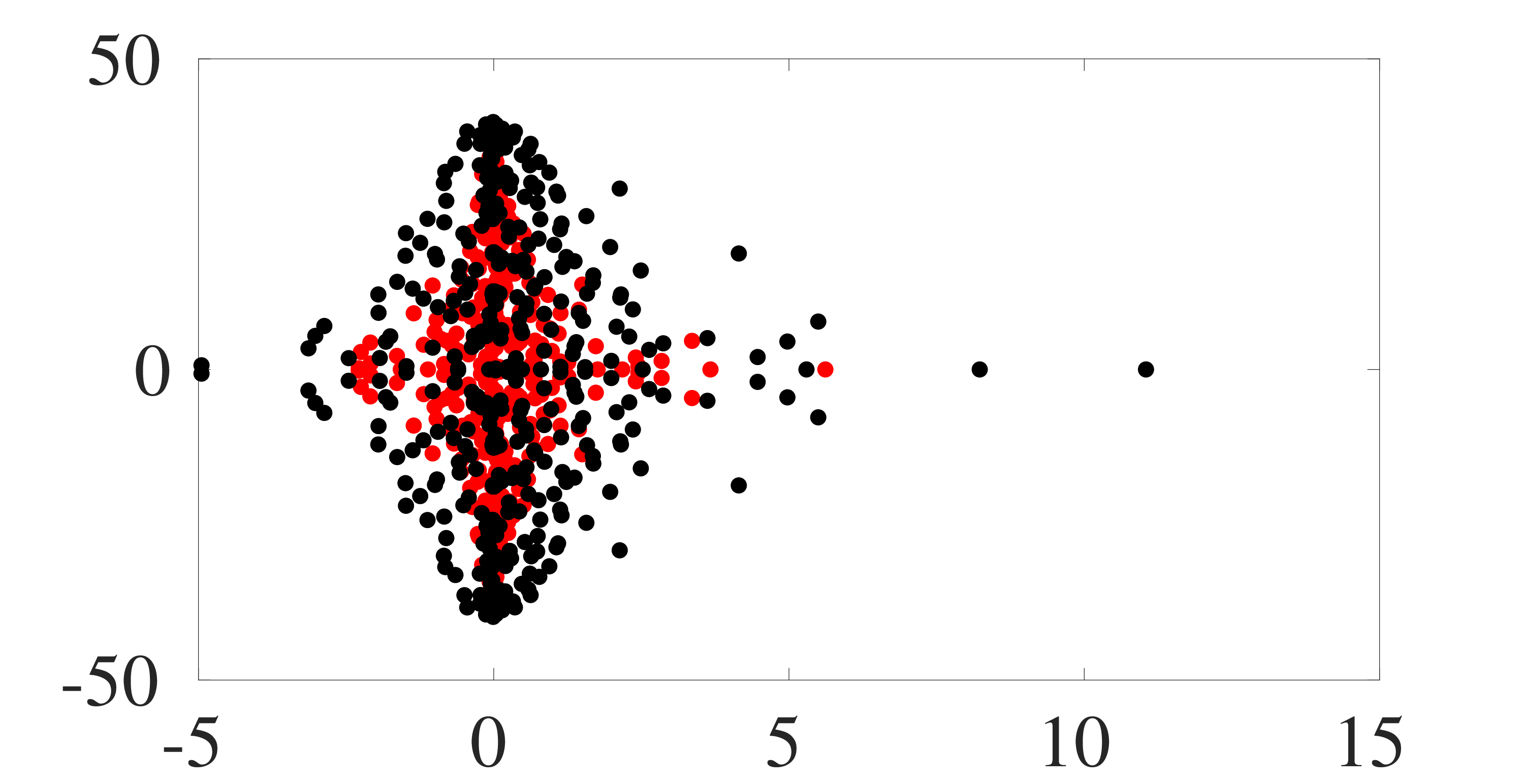}}
     \put(2.5,-0.45){$\Re\{\lambda\}$}
    \put(8.5,-0.45){$\Re\{\lambda\}$}
    \put(14.5,-0.45){$\Re\{\lambda\}$}
     \put(-0.25,1.15){\rotatebox{90}{$\Im\{\lambda\}$}}
     \put(5.75,1.15){\rotatebox{90}{$\Im\{\lambda\}$}}
     \put(11.75,1.15){\rotatebox{90}{$\Im\{\lambda\}$}}
    \end{picture}

    \end{center}
    \caption{Stability of gradient operator for LABFM approximations on a typical grid with $441$ points. Red dots correspond to weights produced by Hermite polynomials with Wendland C2 kernel, black dots to those produced by multi-basis formulation given in \S\ref{subsec::Optimisation_of_coefficients}. From left to right: $m=4$,  $m=6$, $m=8$.}
    \label{fig:labfm_grad_stab}
\end{figure}
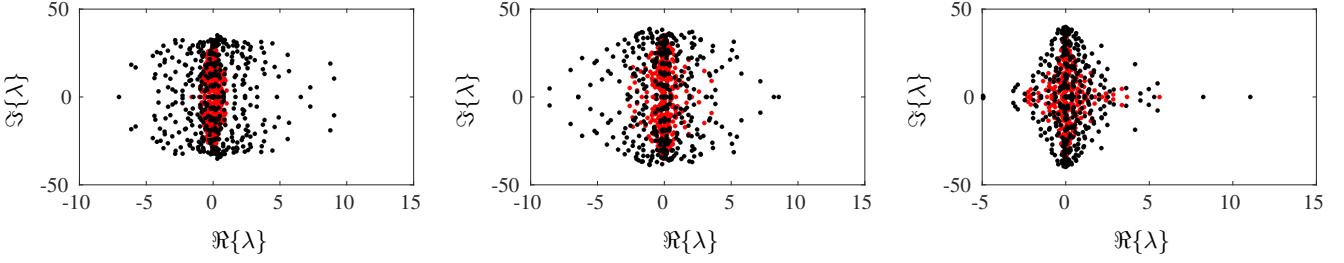

Figs. \ref{fig:sph_grad_stab}, \ref{fig:rbf_grad_stab}, \& \ref{fig:labfm_grad_stab} show eigenvalues for approximations to a gradient operator using the three meshless methods, comparing the usual single-kernel formulations to the multi-kernel formulations of \S\ref{subsec::Optimisation_of_coefficients}. In most cases we observe taking the approach of \S\ref{subsec::Optimisation_of_coefficients} slightly increases growth rate of the most unstable mode, though it always remains the same order of magnitude. For all cases (both SK and MK) we have $\max(\Re\{\boldsymbol{\lambda}\})>0$, therefore purely convective problems are unstable. Such numerical methods are therefore unsuitable for hyperbolic problems (as there is no obvious way to undertake upwinding), and instead are typically used for hyperbolic-parabolic systems, where diffusive terms compensate for the instability of the convective operator. A small increase in the instability of the convective operator due to the MK approach is tolerable because of the increased accuracy and the stabilisation provided by diffusive terms in the systems simulated.

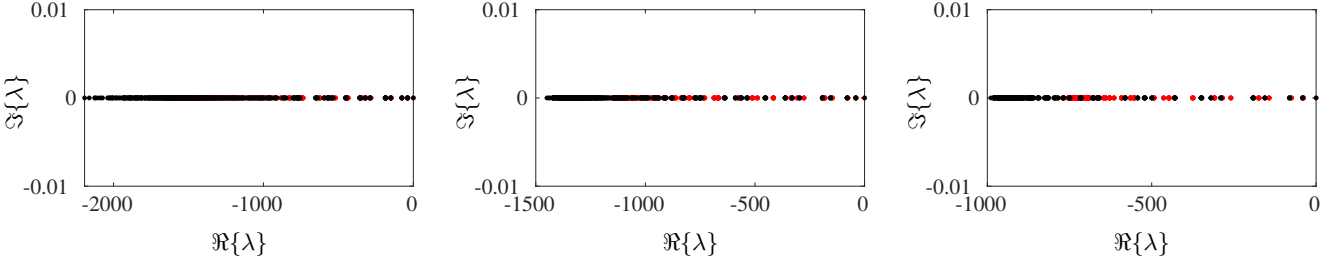
\begin{figure}[t]
    \begin{center}
    \setlength{\unitlength}{1cm}
    \begin{picture}(18,3)(0,0)
    \put(0,0){\includegraphics[width=0.32\linewidth]{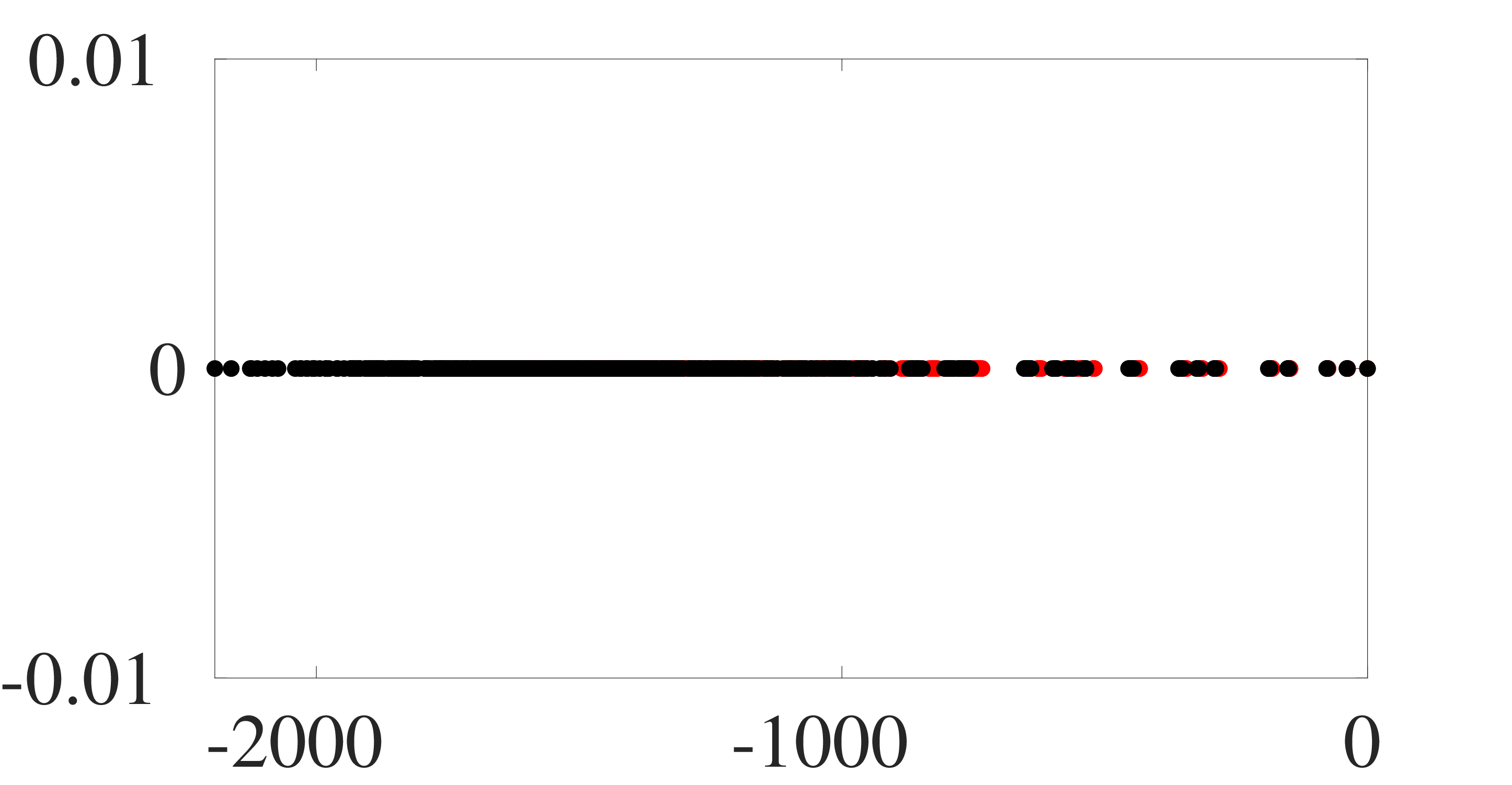}}
    \put(6,0){\includegraphics[width=0.32\linewidth]{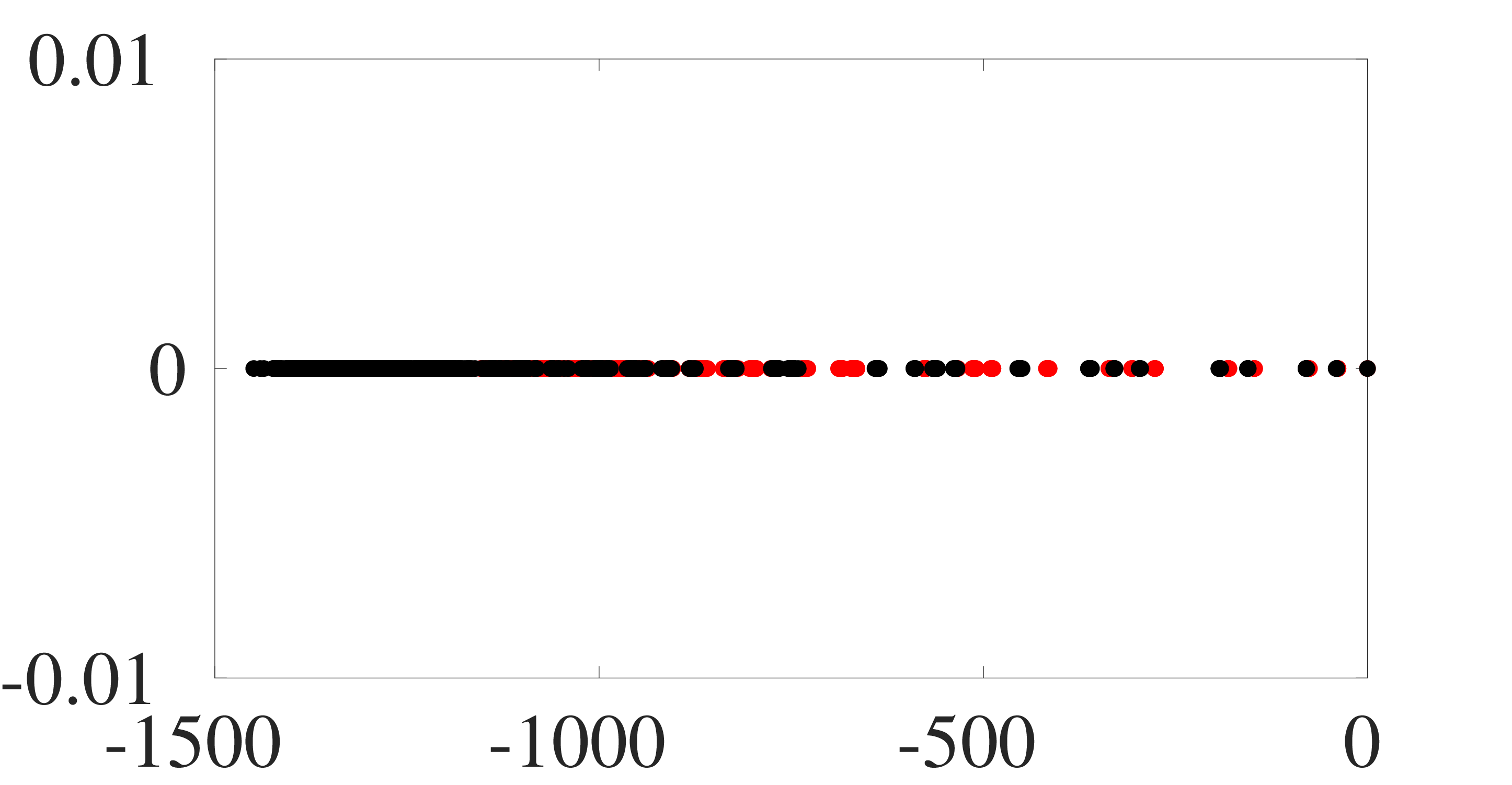}}
    \put(12,0){\includegraphics[width=0.32\linewidth]{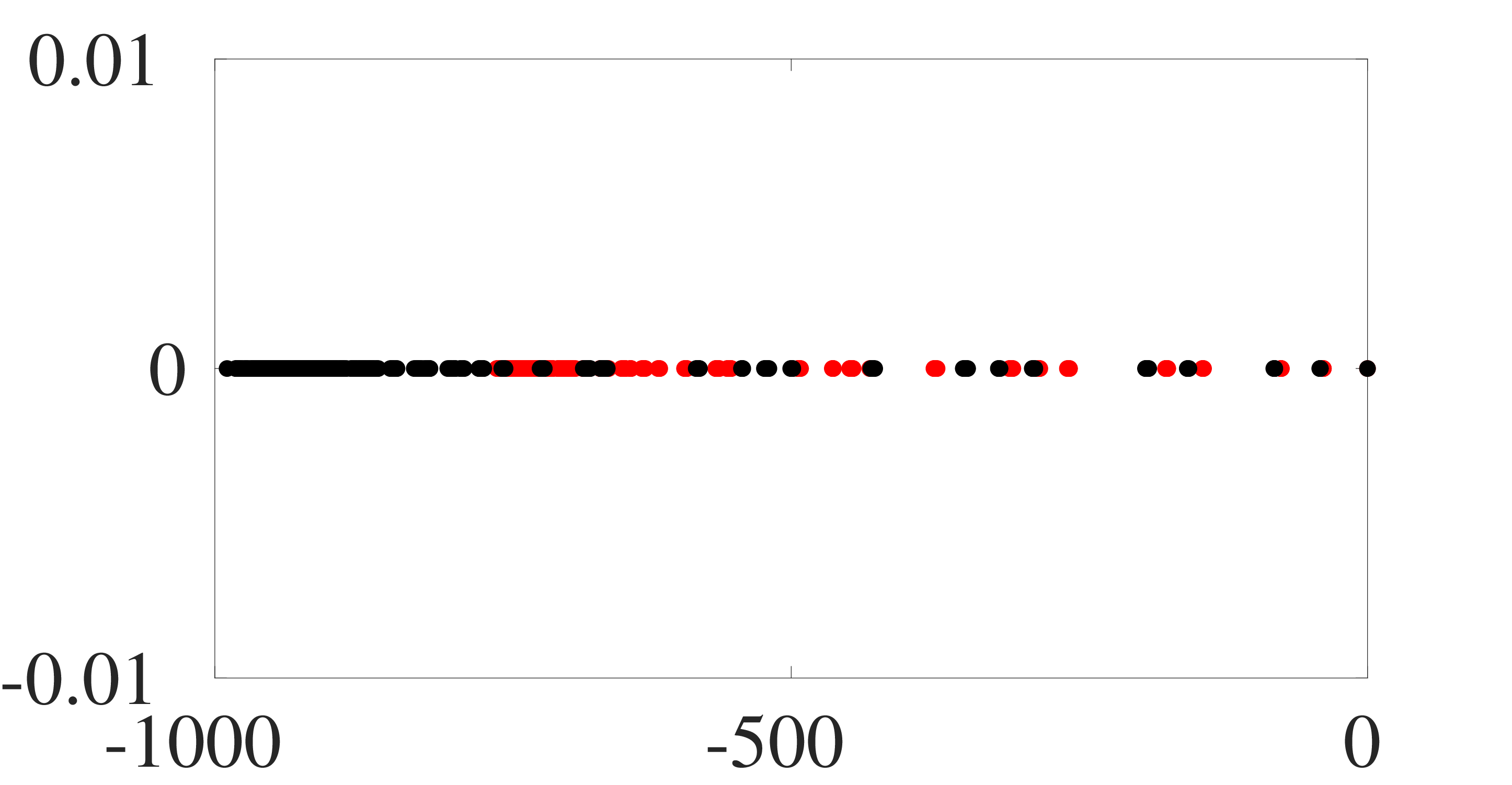}}
     \put(2.5,-0.45){$\Re\{\lambda\}$}
    \put(8.5,-0.45){$\Re\{\lambda\}$}
    \put(14.5,-0.45){$\Re\{\lambda\}$}
     \put(-0.25,1.15){\rotatebox{90}{$\Im\{\lambda\}$}}
     \put(5.75,1.15){\rotatebox{90}{$\Im\{\lambda\}$}}
     \put(11.75,1.15){\rotatebox{90}{$\Im\{\lambda\}$}}
    \end{picture}
    \end{center}
    \caption{Stability of Laplacian operator for SPH approximations on a typical node distribution with $441$ points. Red dots correspond to Wendland C2 kernel, black dots to those produced combination of two kernels as outlined in \S\ref{subsec::Optimisation_of_coefficients}. Left panel: $h/s=1.3$; Centre Panel: $h/s=1.6$; Right Panel: $h/s=2.0$.}
    \label{fig:sph_lap_stab}
\end{figure}

\begin{figure}[t]
    \begin{center}
    \setlength{\unitlength}{1cm}
     \begin{picture}(18,2)(0,0)
    \put(0,0){\includegraphics[width=0.24\linewidth]{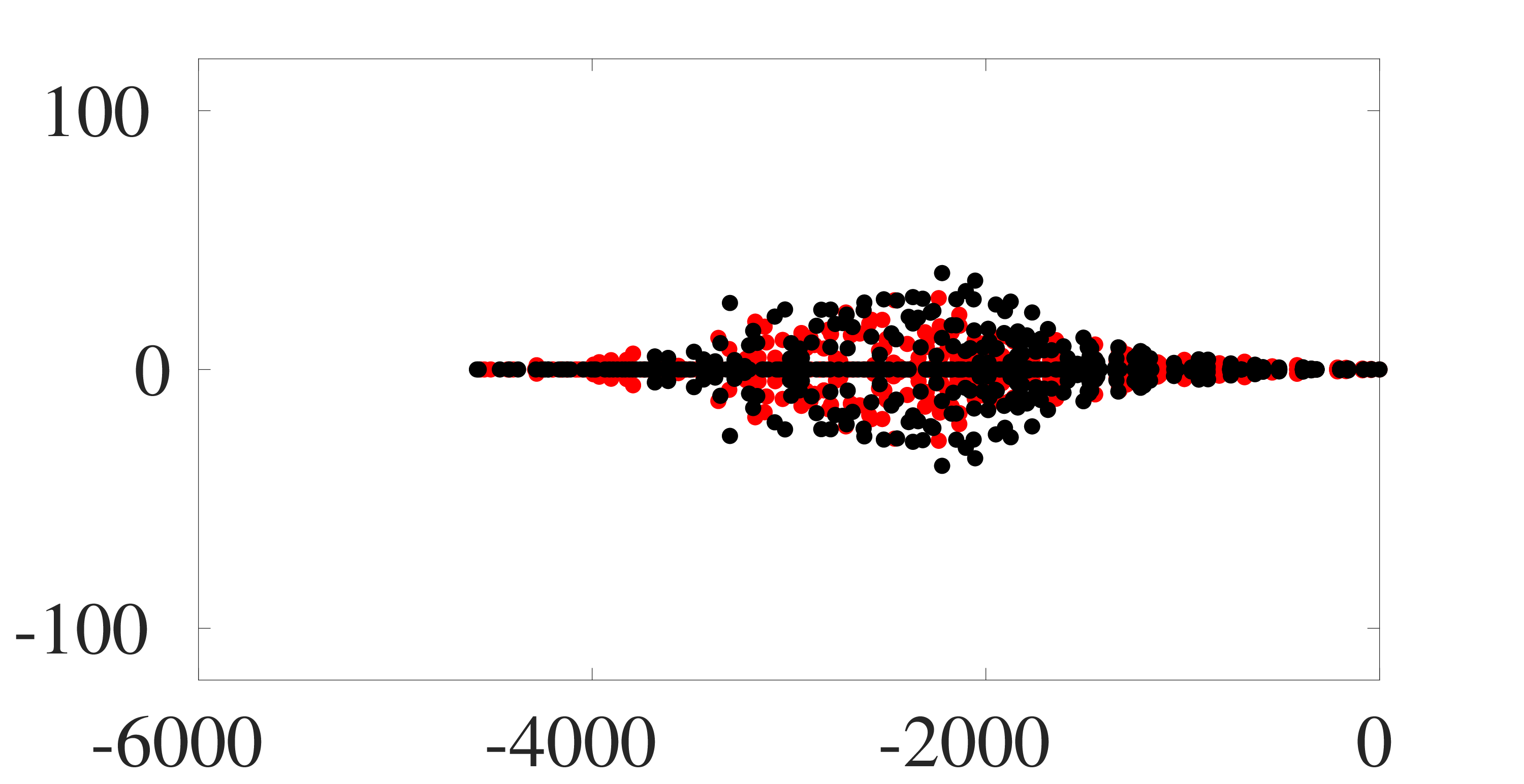}}
    \put(4.5,0){\includegraphics[width=0.24\linewidth]{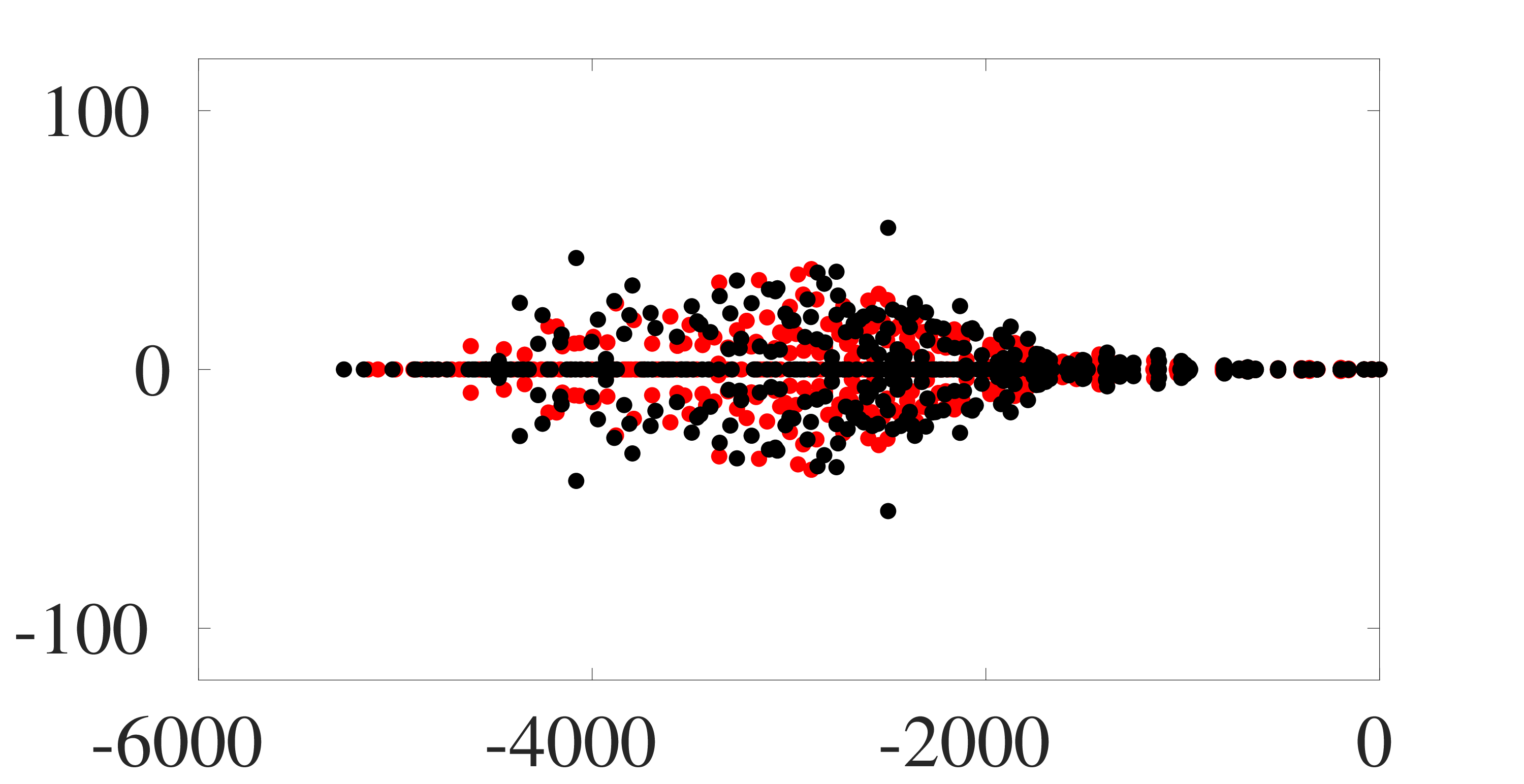}}
    \put(9,0){\includegraphics[width=0.24\linewidth]{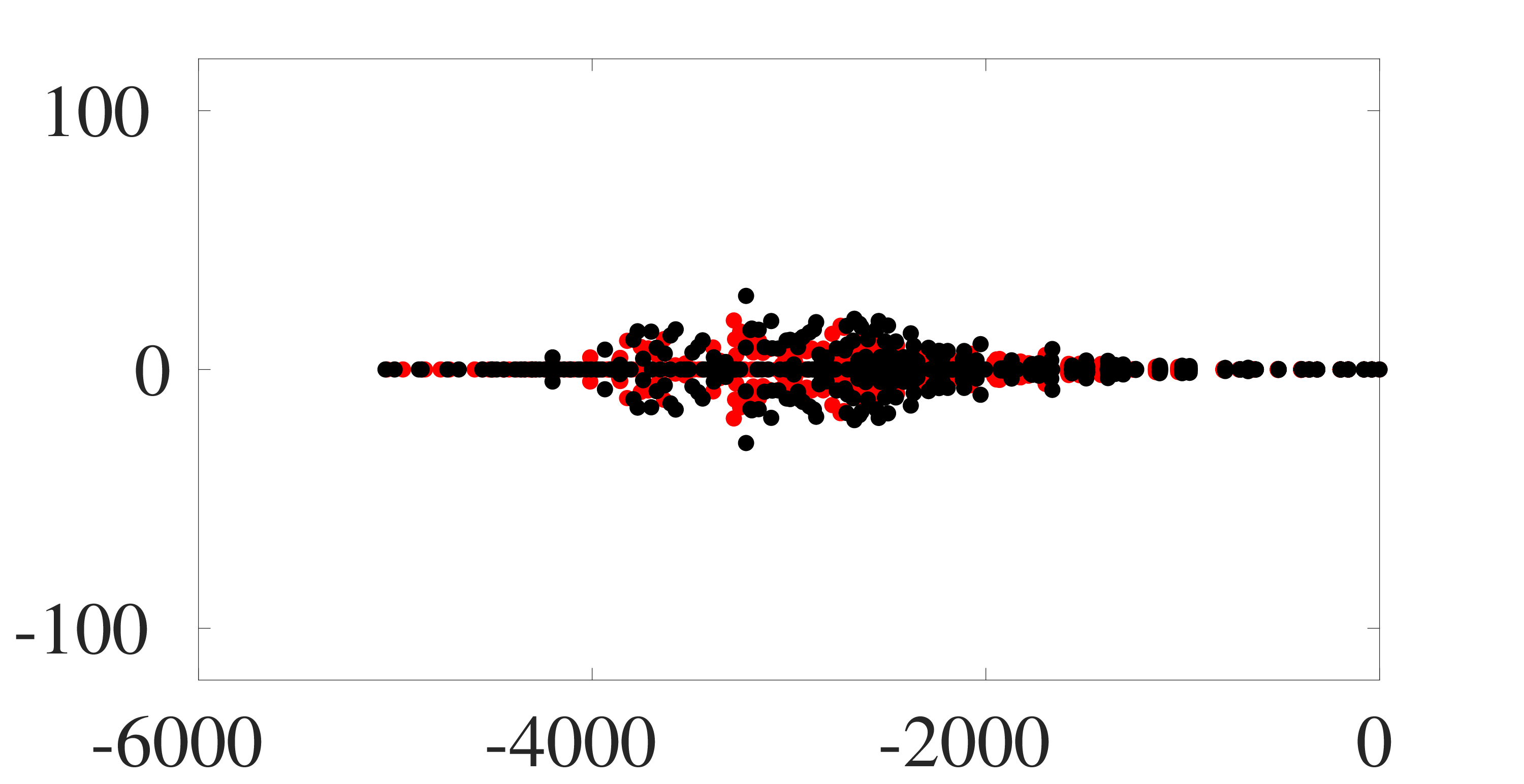}}
    \put(13.5,0){\includegraphics[width=0.24\linewidth]{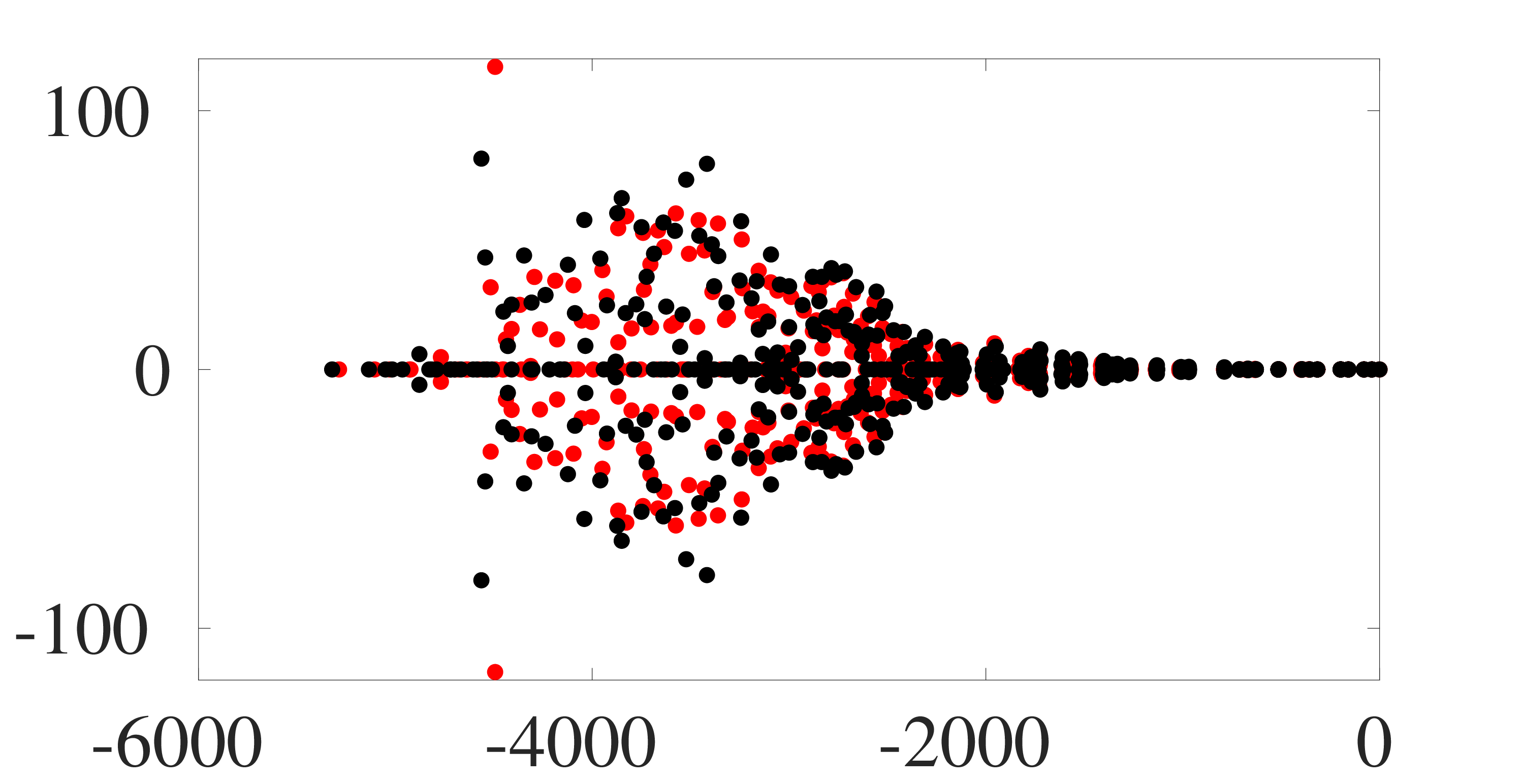}}
     \put(1.9,-0.45){$\Re\{\lambda\}$}
    \put(6.25,-0.45){$\Re\{\lambda\}$}
    \put(10.9,-0.45){$\Re\{\lambda\}$}
     \put(15.25,-0.45){$\Re\{\lambda\}$}
     \put(-0.2,0.8){\rotatebox{90}{$\Im\{\lambda\}$}}
     \put(4.3,0.8){\rotatebox{90}{$\Im\{\lambda\}$}}
     \put(8.8,0.8){\rotatebox{90}{$\Im\{\lambda\}$}}
     \put(13.2,0.8){\rotatebox{90}{$\Im\{\lambda\}$}}
    \end{picture}
    \end{center}
    \caption{Stability of Laplacian operator for RBF-FD approximations on a typical grid with $441$ points. Red dots correspond to Gaussian RBF, black dots to those produced by linear combination of RBFs as detailed in \S\ref{subsec::Optimisation_of_coefficients}. From left to right: $m=0$,  $\lvert\mathcal{N}_i\rvert=10$; $m=1$, $\lvert\mathcal{N}_i\rvert$=15; $m=2$, $\lvert\mathcal{N}_i\rvert=20$; $m=3$, $\lvert\mathcal{N}_i\rvert=25$.}
    \label{fig:rbf_lap_stab}
\end{figure}

\begin{figure}[t]
    \begin{center}
    \setlength{\unitlength}{1cm}
    \begin{picture}(18,3)(0,0)
    \put(0,0){\includegraphics[width=0.32\linewidth]{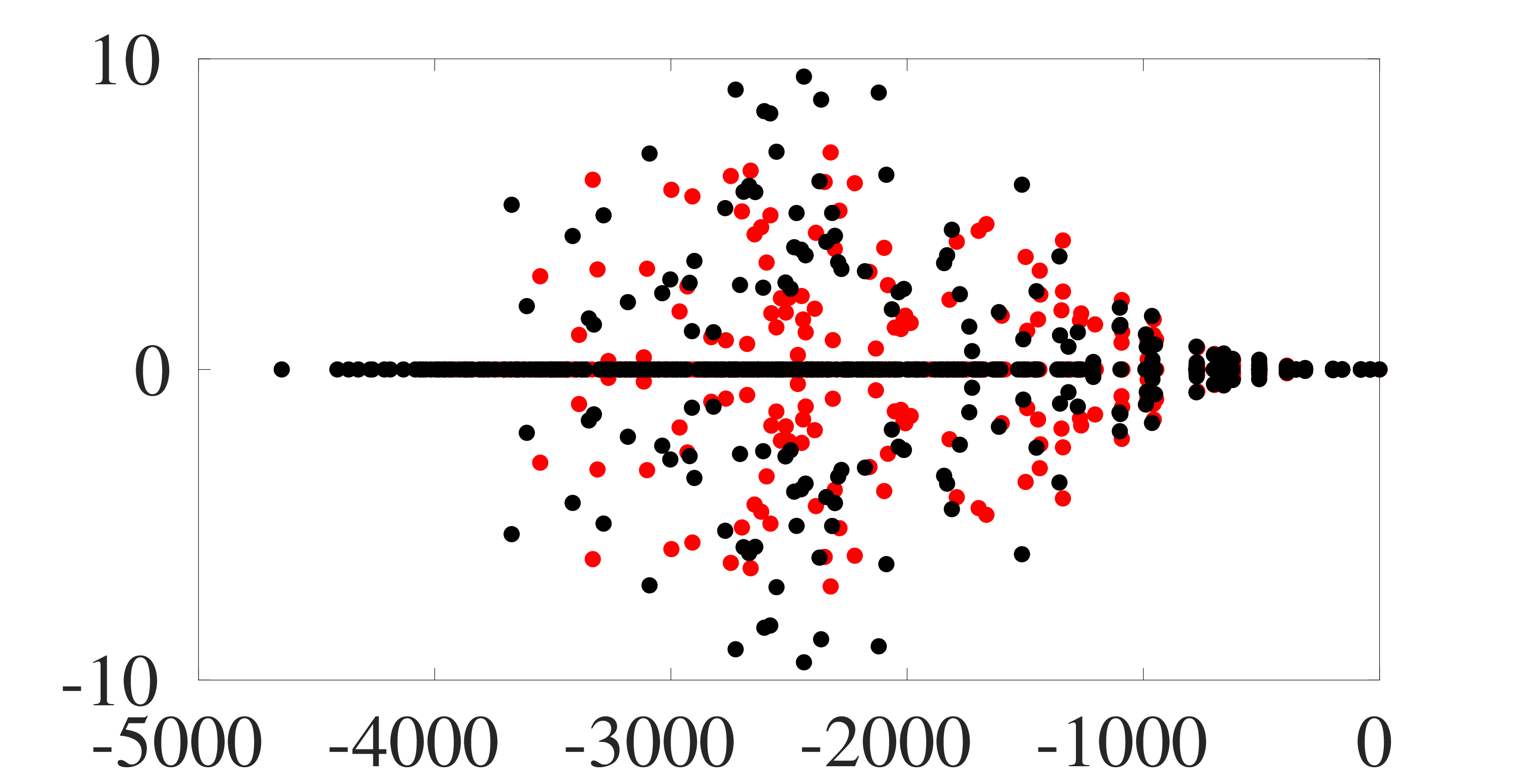}}
    \put(6,0){\includegraphics[width=0.32\linewidth]{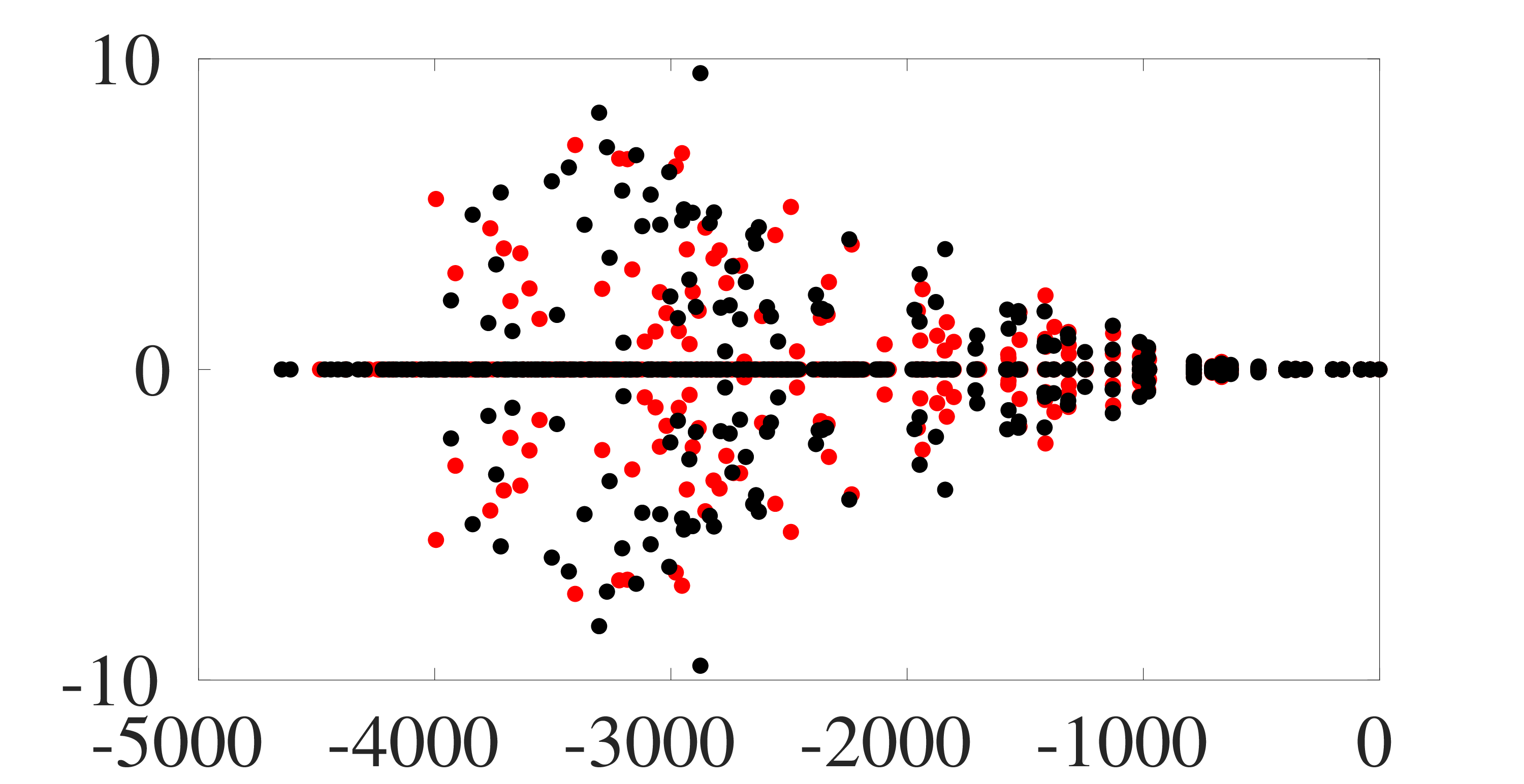}}
    \put(12,0){\includegraphics[width=0.32\linewidth]{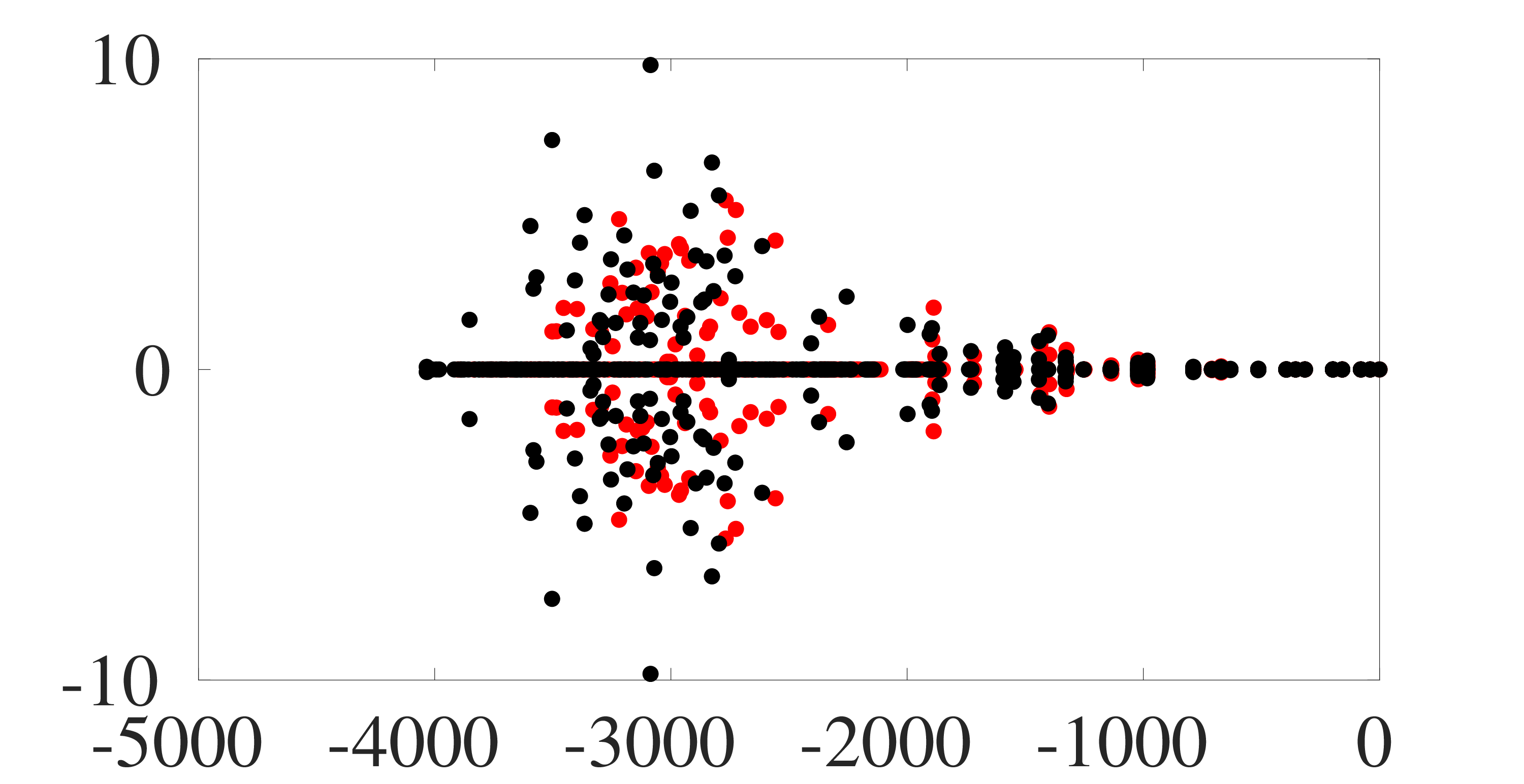}}
     \put(2.5,-0.45){$\Re\{\lambda\}$}
    \put(8.5,-0.45){$\Re\{\lambda\}$}
    \put(14.5,-0.45){$\Re\{\lambda\}$}
     \put(-0.25,1.15){\rotatebox{90}{$\Im\{\lambda\}$}}
     \put(5.75,1.15){\rotatebox{90}{$\Im\{\lambda\}$}}
     \put(11.75,1.15){\rotatebox{90}{$\Im\{\lambda\}$}}
    \end{picture}
    \end{center}
    \caption{Stability of Laplacian operator for LABFM approximations on a typical grid with $441$ points.Red dots correspond to weights produced by Hermite polynomials with Wendland C2 kernel, black dots to those produced by multi-basis formulation given in \S\ref{subsec::Optimisation_of_coefficients}. From left to right: $m=4$,  $m=6$, $m=8$.}
    \label{fig:labfm_lap_stab}
\end{figure}

We can also analyse the eigenvalues of the global derivative matrix for the Laplacian operator. In this case the ideal scenario would be $\Re\{\lambda_i\}\leq0,\ \Im\{\lambda_i\}=0, \forall i$, though the latter condition is not a requirement for stability. Figs. \ref{fig:sph_lap_stab}, \ref{fig:rbf_lap_stab}, \& \ref{fig:labfm_lap_stab} show example eigenvalue plots for approximation to the Laplacian operator, again comparing the SK approach of each meshless methods to the MK approach of \S\ref{subsec::Optimisation_of_coefficients}. In all cases there are no modes that grow in time, though the new formulations increases the phase speed error. Given the typically large values of $\Re\{\lambda\}$ (hence rapid decay in time), this behaviour is unlikely to present issues in time-dependent simulations.

\section{Application to PDEs}
\label{sec:PDE_solns}
We now test the impact of using a multi-kernel expansion to improve numerical solutions to systems of PDEs. The majority of our test cases will be in simple domains rather than complex geometries to which meshless methods are most suited due to the lack of analytic solutions to PDEs in complex domains. As in \S\ref{sec::Convergance_Studies}, we shall measure errors using the relative $L_2$ norm as defined in \eqref{eq::l2_norm_defn}.

In \S\ref{sec::Convergance_Studies} we demonstrated that the new formulation results in improved resolving power for three different numerical methods. For brevity of exposition and clarity we shall now focus solely on one method, LABFM, though similar improvements were found for both RBF-FDs and SPH operators.

\subsection{Poisson's Equation}
The first PDE system we consider is Poisson's equation
\begin{equation}
    \label{eq::Poisson_equation}
    \nabla^2 \phi=f(x,y)
\end{equation}
in domain $\Omega$.

We begin by considering a periodic test case in the square domain $(x,y)\in[0,1]\times[0,1]$. We set
\begin{equation}
    f(x,y)=-8\pi^2\sin(2\pi x)\sin(2\pi y)
\end{equation}
for which the analytic solution to \eqref{eq::Poisson_equation} is
\begin{equation}
    \phi(x,y)=\sin(2\pi x)\sin(2\pi y).
\end{equation}
To solve Poisson's equation we are required to solve a global (sparse) linear system which is constructed in the same manner as \cite{king_2020}. This procedure is performed using an in-house solver using the Bi-Conjugate Gradient Stabilised (BiCGStab) method with Jacobi preconditioner.

\begin{figure}[t]
    \begin{center}
    \setlength{\unitlength}{1cm}
    \begin{picture}(18,5)(0,0)
    \put(0,0){\includegraphics[width=0.49\linewidth]{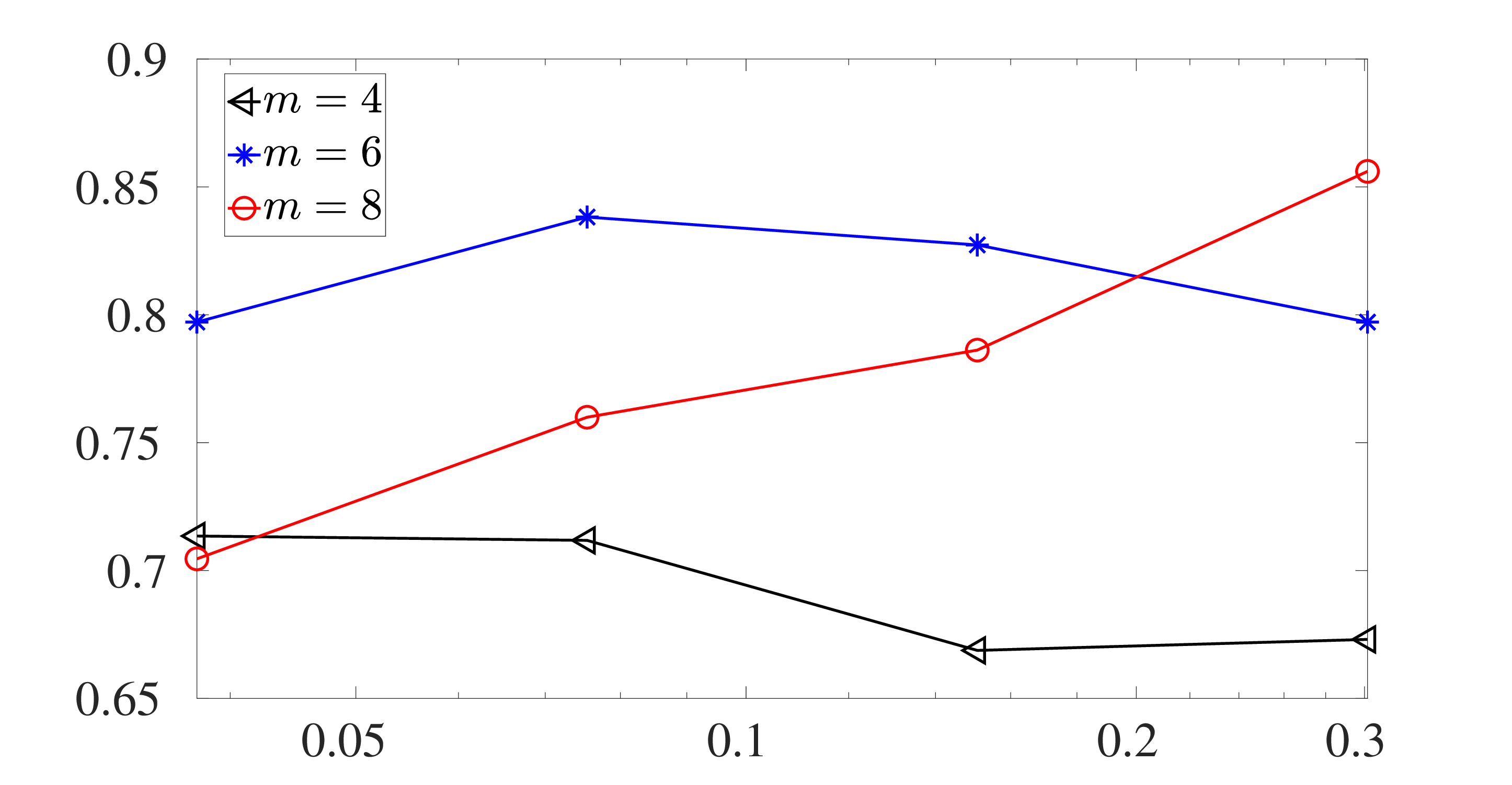}}
    \put(9,0){\includegraphics[width=0.49\linewidth]{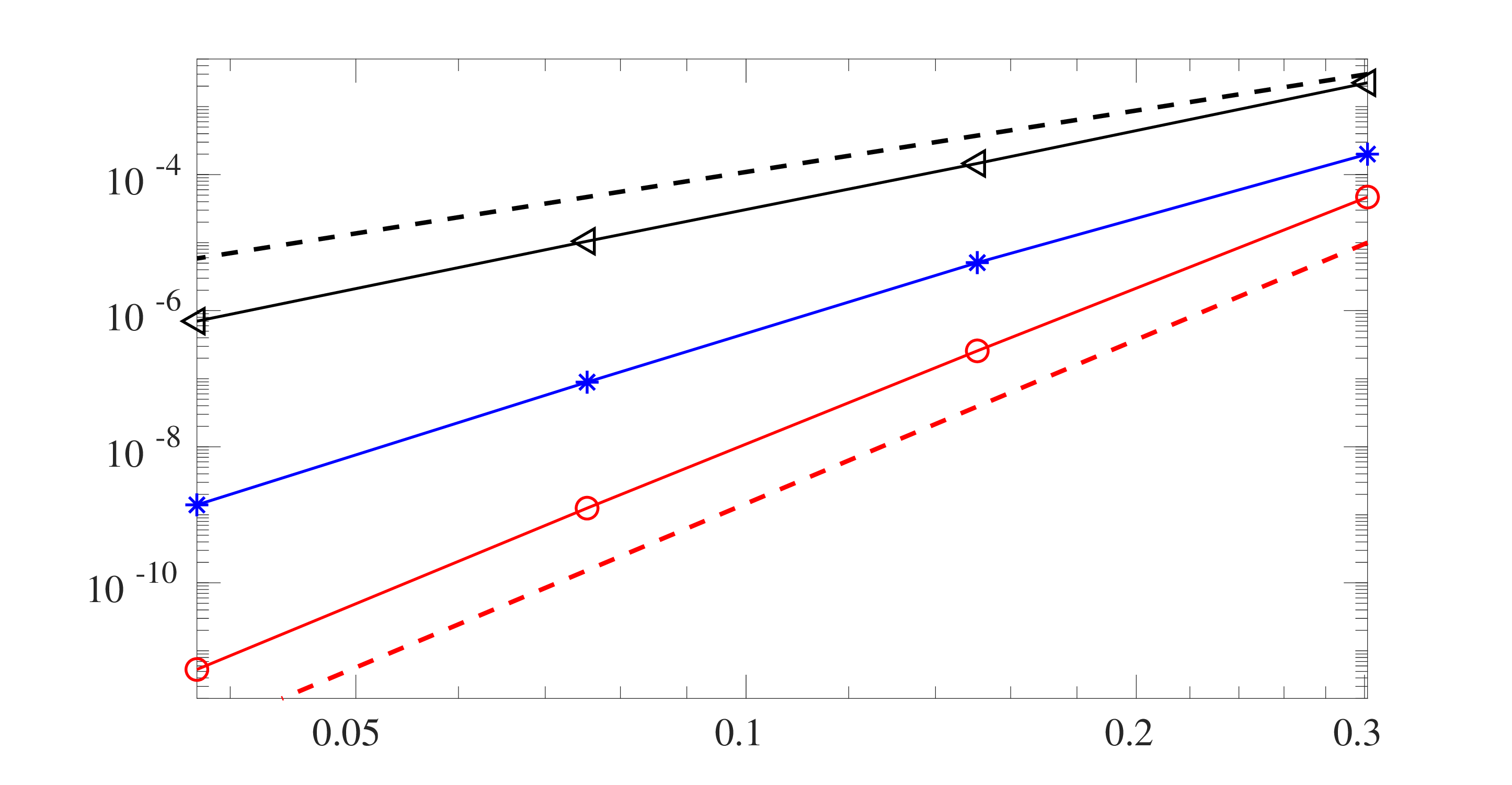}}
    \put(13.0,-0.25){$s$}
    \put(-0.1,2.25){\rotatebox{90}{$\mathcal{R}$}}
    \put(8.9,1.75){\rotatebox{90}{$L_2-$norm}}
   
    \put(4.0,-0.25){$s$}
    
    \end{picture}
    \end{center}
    \caption{Solution to periodic Poisson equation using LABFM. Left panel: comparison of ratio $\mathcal{R}$ of solutions. Right panel: convergence of new approach (dashed lines show convergence rates of 3 (black) and 7 (red)). Black $\triangleleft$ $m=4$; Blue $*$ $m=6$; Red $\circ$ $m=8$.}
    \label{fig:periodic_poisson}
\end{figure}

The left panel of Fig. \ref{fig:periodic_poisson} compares the ratio $\mathcal{R}$ of the $L_2-$ norms produced by the new approach and standard implementation for polynomial consistencies of LABFM. Similar to the Laplacian convergence study of $\S$\ref{sec::Convergance_Studies} we observe a consistent improvement in the accuracy of solutions across all resolutions and orders. Although the error reduction here is somewhat modest (in the range of $15-35\%$,) we note that the wavenumbers present in $\phi$ are small even on the coarsest meshes simulated. The right panel demonstrates the convergence of the new approach, with the solutions showing the same convergence rates as expected.

A more challenging test case is now examined. A key advantage of meshless methods is the ability to efficiently discretise non-trivial domains. We consider Poisson's equation over a circle of radius 1 with Dirichlet boundary conditions
\begin{gather}
    \phi= g_1(x,y)\qquad \mathrm{on}\quad \ \Gamma,
\end{gather}
where $\Gamma$ is the boundary of the domain. Taking $(r,\theta)$ to be the usual polar coordinates we set
\begin{equation}
    f(r,\theta)=\left[\left( 15r^2-\frac{\pi^2}{4}r^4\right)\cos\left(\frac{\pi r}{2}\right)-\frac{9\pi^2}{2}r^3\sin\left(\frac{\pi r}{2}\right)\right]\cos\left(\theta\right)
\label{eq:poisson_disc1}
\end{equation}
and set the Dirichlet boundary condition
\begin{equation}
    \phi(1,\theta)=0.
    \label{eq:poisson_disc2}
\end{equation}
This problem has the solution
\begin{equation}
    \phi(r,\theta)=r^4\cos\left(\frac{\pi r}{2}\right)\cos\left(\theta\right).
    \label{eq:poisson_disc3}
\end{equation}
We note that whilst we have expressed~\eqref{eq:poisson_disc1} to~\eqref{eq:poisson_disc3} in polar coordinates, our numerical calculations are performed in a Cartesian coordinate system. Due to the size of stencils required at the boundary (necessitated by incomplete support --- similar issues occur in RBF-FDs \cite{Bayona_2017}) we consider LABFM approximations of only fourth and sixth order polynomial consistencies.

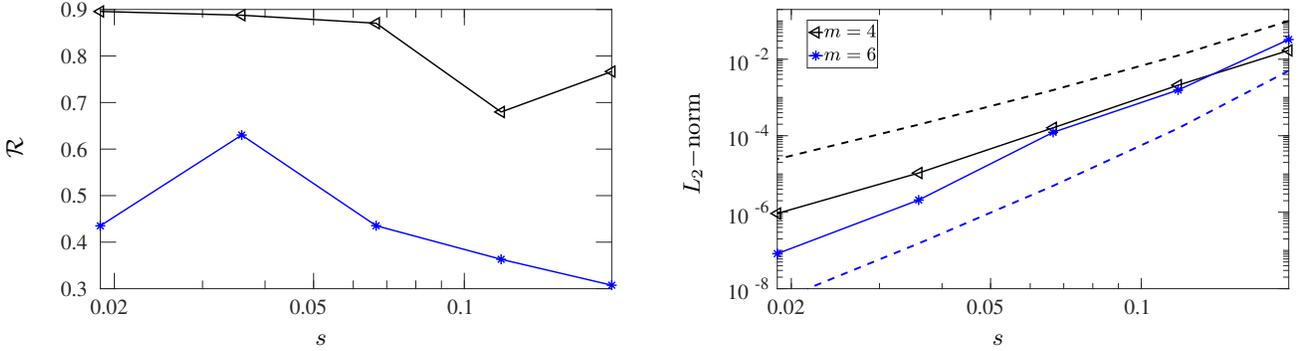
\begin{figure}[t]
    \begin{center}
    \setlength{\unitlength}{1cm}
    \begin{picture}(18,5)(0,0)
    \put(0,0){\includegraphics[width=0.49\linewidth]{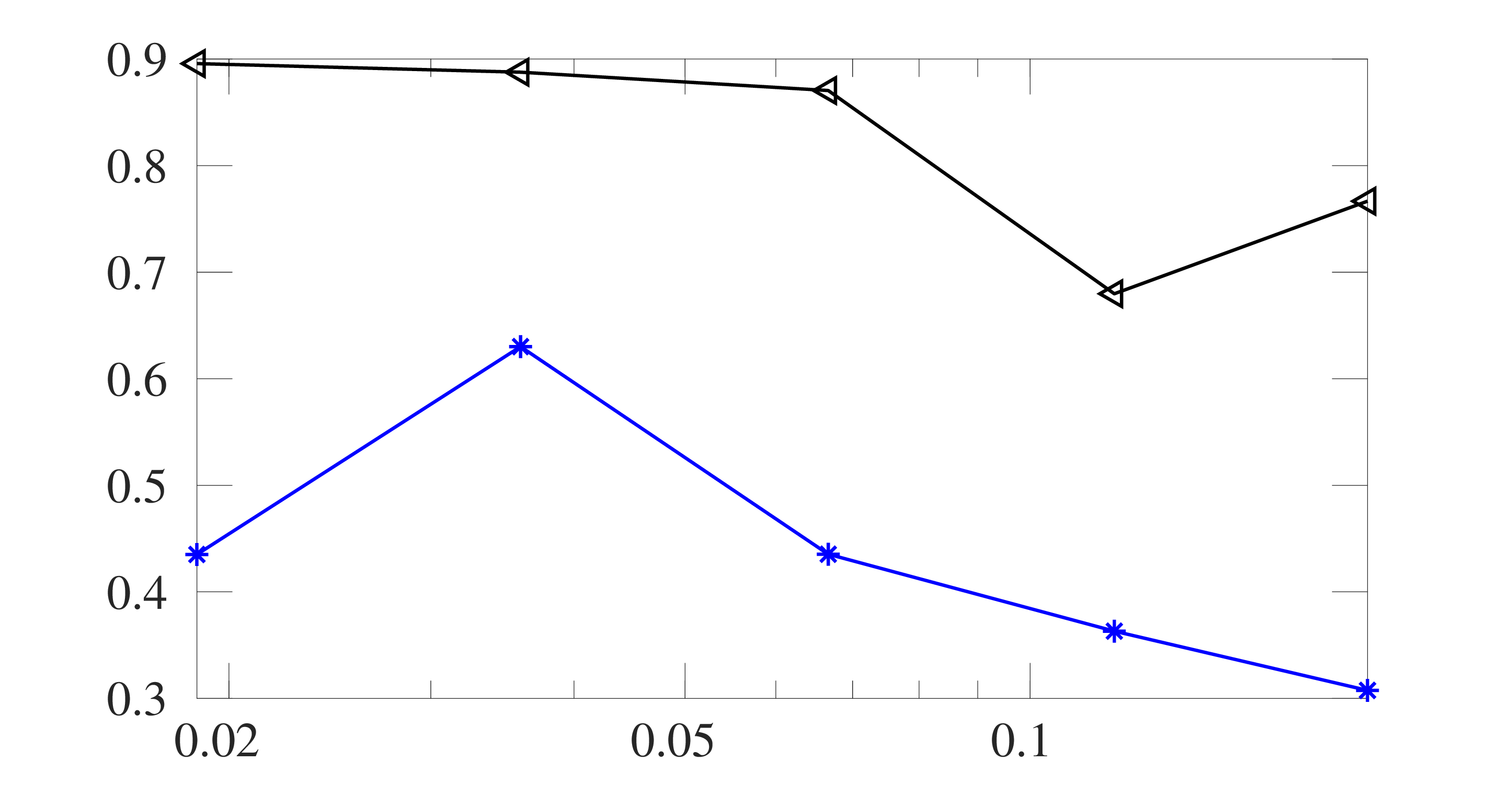}}
    \put(9,0){\includegraphics[width=0.49\linewidth]{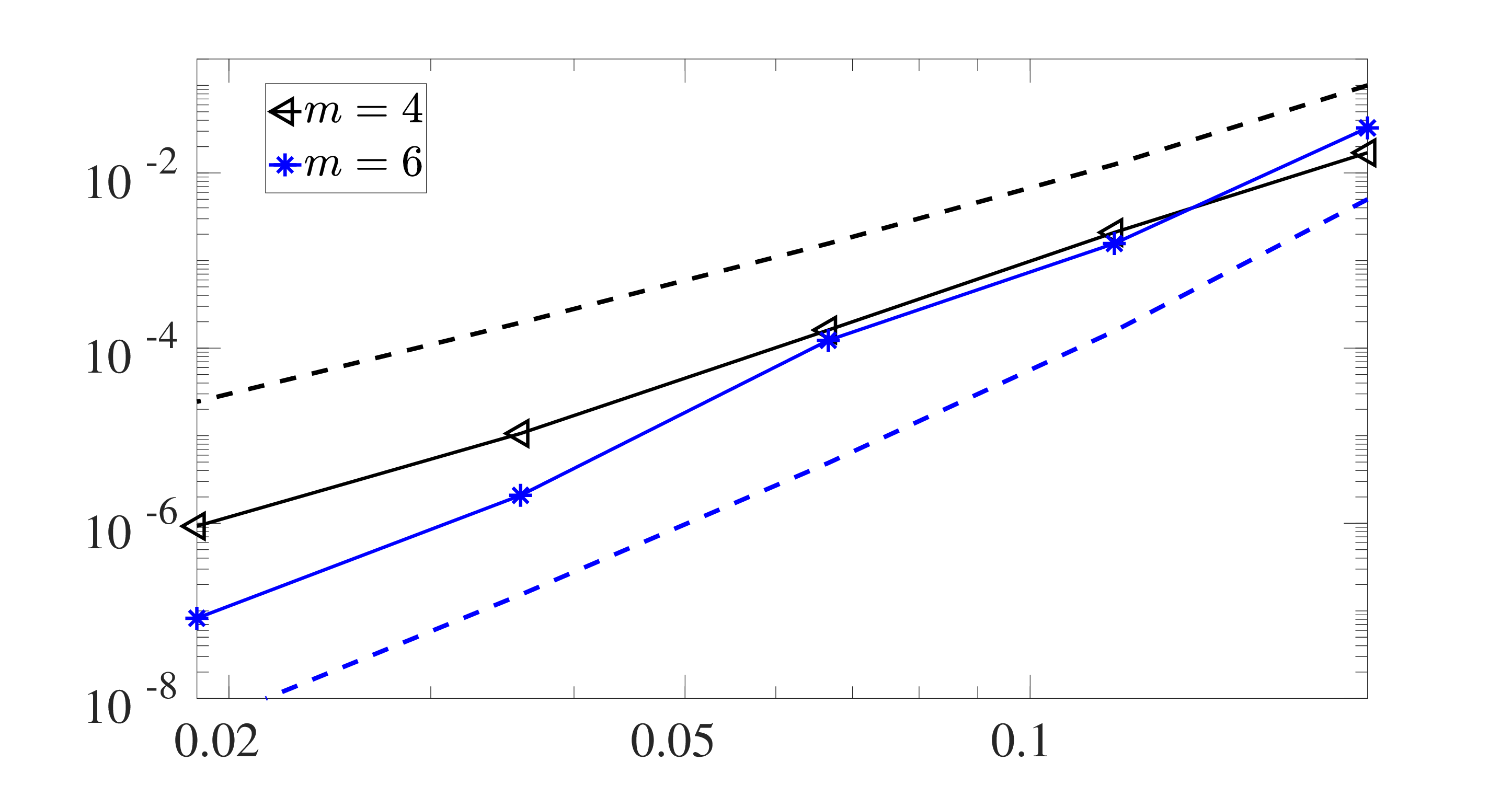}}
    \put(13.0,-0.25){$s$}
    \put(-0.1,2.25){\rotatebox{90}{$\mathcal{R}$}}
    \put(8.9,1.75){\rotatebox{90}{$L_2-$norm}}
   
    \put(4.0,-0.25){$s$}
    
    \end{picture}
    \end{center}
    \caption{Solution to Poisson's equation within a circular geometry using LABFM. Left panel: comparison of ratio $\mathcal{R}$ of solutions. Right panel: convergence of new approach (dashed lines show convergence rates of 3 (black) and 5 (blue)). Black $\triangleleft$ $m=4$; Blue $*$ $m=6$.}
    \label{fig:circular_poisson}
\end{figure}

The left panel in Fig. \ref{fig:circular_poisson} compares the $L_2-$norm of the new approach to that obtained by a standard implementation. The multi-kernel approach leads to a consistent improvement across all resolutions, with the higher order approximation seeing significant gains. This is likely due to increased accuracy near the boundary. The right panel shows the convergence of the solution. Convergence rates appear slower due to incomplete support of stencils near the boundary (see \cite{king_2020} for details). Development of high-order one-sided meshless operators is an ongoing area of research within our group.

\subsection{Advection-Diffusion Equation}
We now test our numerical methods for time-dependent systems of PDEs. Time integration is performed explicitly using the classical 4th order Runge-Kutta (RK4) scheme in all cases.

We begin by considering an advection-diffusion system
\begin{gather}
\label{eq::Advec_diff_eqn}
    \frac{\partial u}{\partial t}+\boldsymbol{a}\cdot\nabla u=\frac{1}{Re} \nabla^2u
\end{gather}
where $u(x,y)$ is a scalar, $\boldsymbol{a}$ is a vector of constants (potentially dependent on $(x,y)$) and $Re$ is a constant. Herein we take $\boldsymbol{a}=\left(\begin{smallmatrix}1\\0\end{smallmatrix}\right)$ to be constant throughout the domain in order to isolate resolving power effects (i.e. errors due to numerical dispersion) from effects of numerical anisotropy. We solve in the doubly-periodic domain $(x,y)\in [0,H]\times [0,H]$ and set the time step through
\begin{equation}
\label{eq::time_step}
    \delta t=\min \left( \frac{0.1 s}{\mathrm{max}(|u|)} ,\frac{0.05 s^2}{Re}\right).
\end{equation}
 To illustrate the dependence of resolving power of operators on the direction of waves in wavenumber space we shall consider two different initial conditions in our simulations. For Case 1 we set $H=1$ and use the initial condition
 \begin{equation}
     \label{eq::case_IC}
     u(x,y)=\sum_{m=1}^M\sin\left(2m\pi x\right)
 \end{equation}
for which \eqref{eq::Advec_diff_eqn} has the analytic solution
\begin{equation}
    \label{eq::ad_diff_soln1}
    u(x,t)=\sum_{m=1}^Me^{-\frac{4m^2\pi^2t}{Re}}\sin(2\pi m(x-at)).
\end{equation}
In the second case we have $H=\sqrt{2}$ and initially set
\begin{equation}
    \label{eq::IC_2}
    u(x,y)=\sum_{m=1}^M \sin\left(\sqrt{2}m\pi(x+y)\right)
\end{equation}
so that our numerical solution can be compared to the analytic expression
\begin{equation}
    \label{eq::ad_diff_soln2}
    u(x,y,t)=\sum_{m=1}^Me^{-\frac{4m^2\pi^2t}{Re}}\sin\left(\sqrt{2}m\pi(x+y-at)\right).
\end{equation}
When comparing solutions on different domains we ensure that the nodal spacing $s$ (and therefore the Nyquist wavenumber $k_{ny}$) remains constant to enable like-for-like comparison. We set the number of modes present in the solution $M=3$ in all cases noting on the coarsest discretisation the largest wavenumber of the solution is $\frac{3}{5}k_{ny}$. In this subsection we focus on approximations which have a polynomial consistency of $m=6$, noting that similar results were found for other orders of approximation.

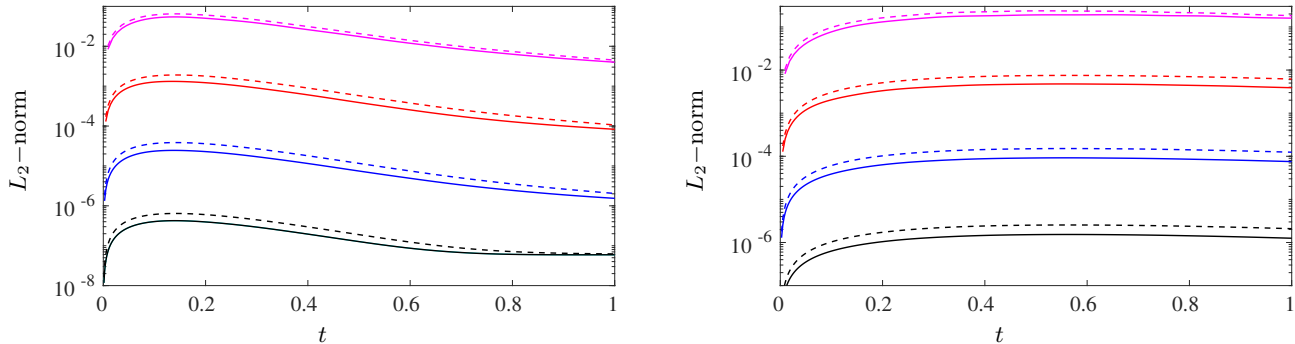
\begin{figure}
    \begin{center}
    \setlength{\unitlength}{1cm}
    \begin{picture}(18,5)(0,0)
    \put(0,0){\includegraphics[width=0.49\linewidth]{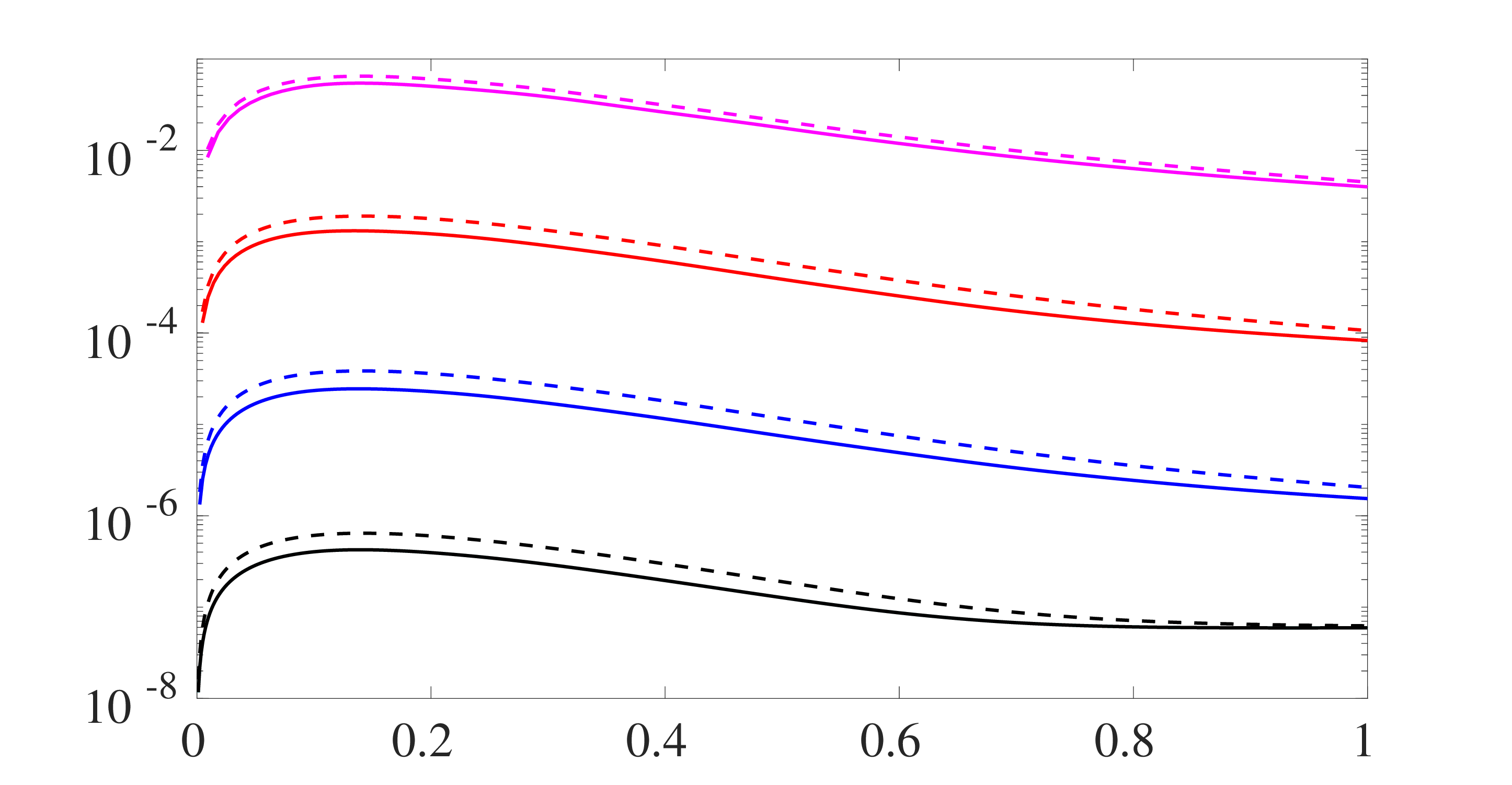}}
    \put(9,0){\includegraphics[width=0.49\linewidth]{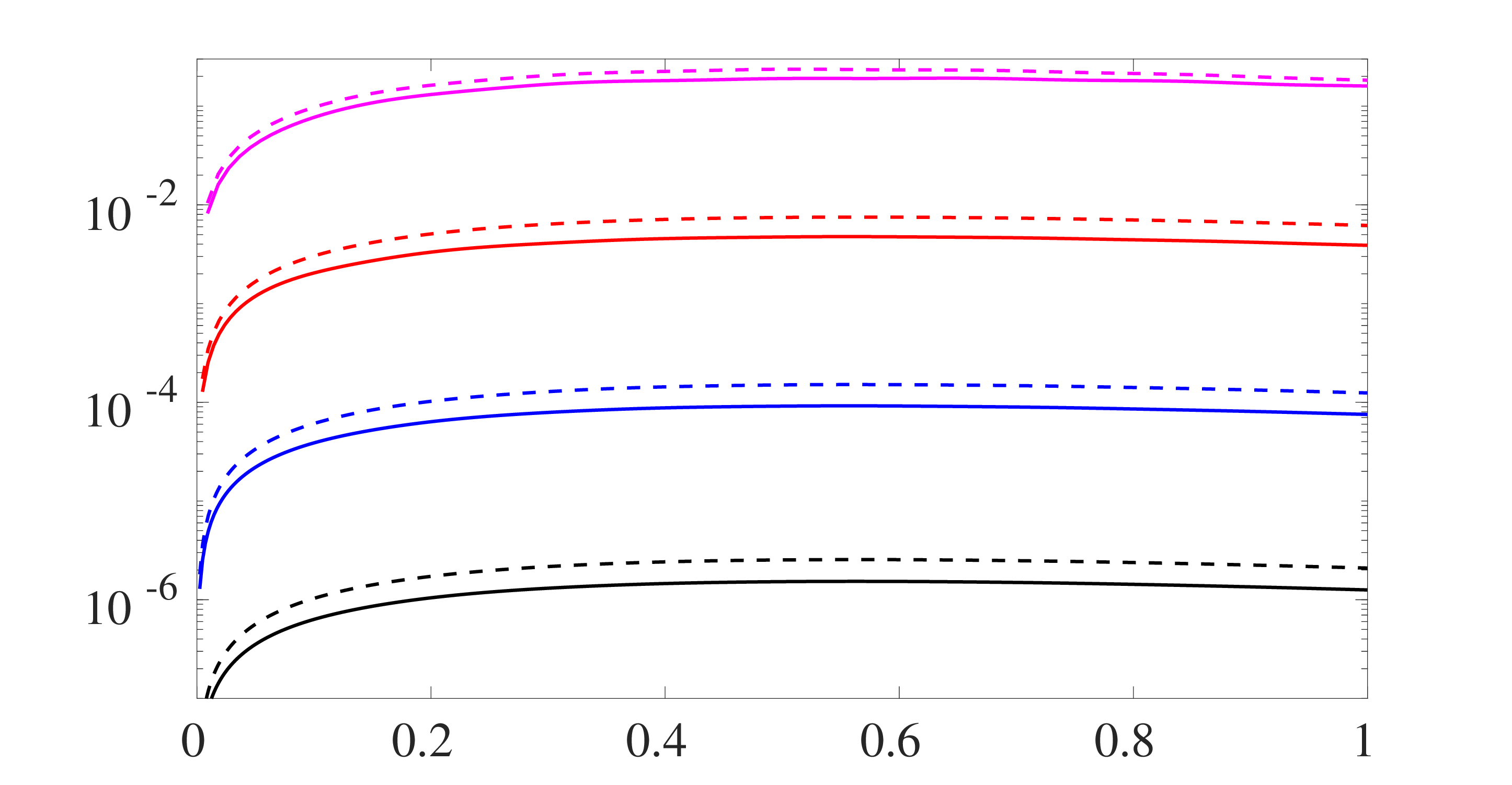}}
    \put(13.0,-0.25){$t$}
    \put(-0.1,1.75){\rotatebox{90}{$L_2-$norm}}
    \put(8.9,1.75){\rotatebox{90}{$L_2-$norm}}
   
    \put(4.0,-0.25){$t$}
    
    \end{picture}
    \end{center}
    \caption{Time evolution of $L_2-$norm of the error of the numerical solution to advection-diffusion equation for Case 1. Solid lines denote MK approach, dashed lines SK approach. Left panel: $Re=50$. Right panel: $Re=200$. Solid lines Pink $s=1/10$; red $s=1/20$; blue $s=1/40$; black $s=1/80$. }
    \label{fig:AD_case1}
\end{figure}

We begin by considering Case 1. Fig. \ref{fig:AD_case1} shows how the $L_2-$norm of the solution evolves in time for different choices of $Re$. In all cases the multi-kernel approach reduces error in the numerical simulation. For larger $Re$ a large improvement is typically observed, due to the (physical and numerical) solution decaying on a longer time scale. At short times, the MK approach typically leads to a $20-40\%$ reduction in error magnitude at little to no extra computational cost per time-step.  

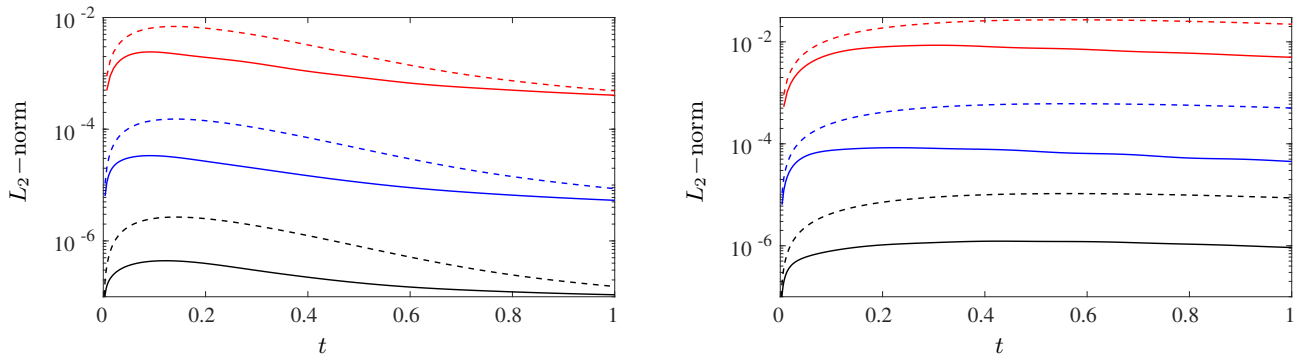
\begin{figure}
    \begin{center}
    \setlength{\unitlength}{1cm}
    \begin{picture}(18,5)(0,0)
    \put(0,0){\includegraphics[width=0.49\linewidth]{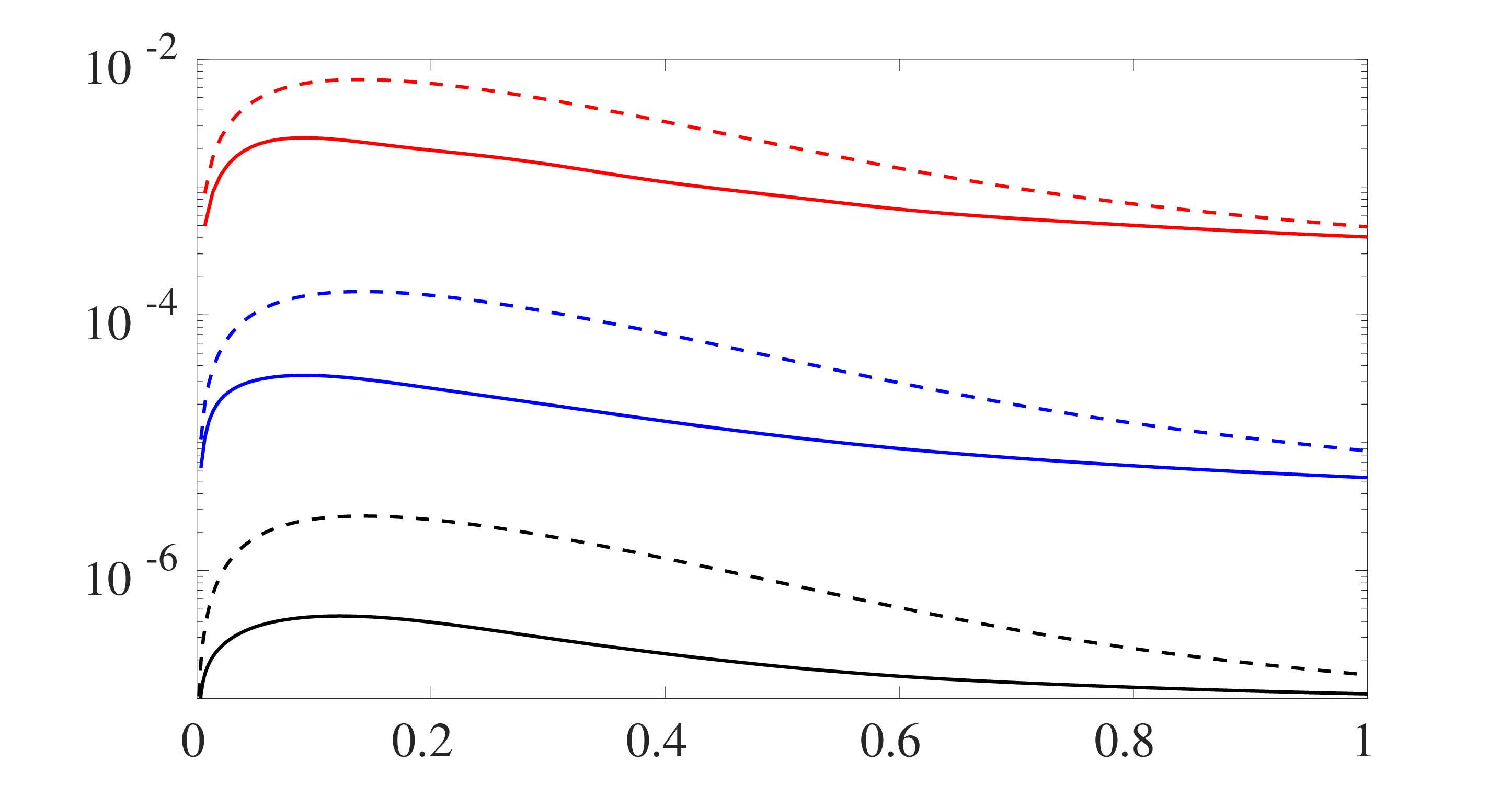}}
    \put(9,0){\includegraphics[width=0.49\linewidth]{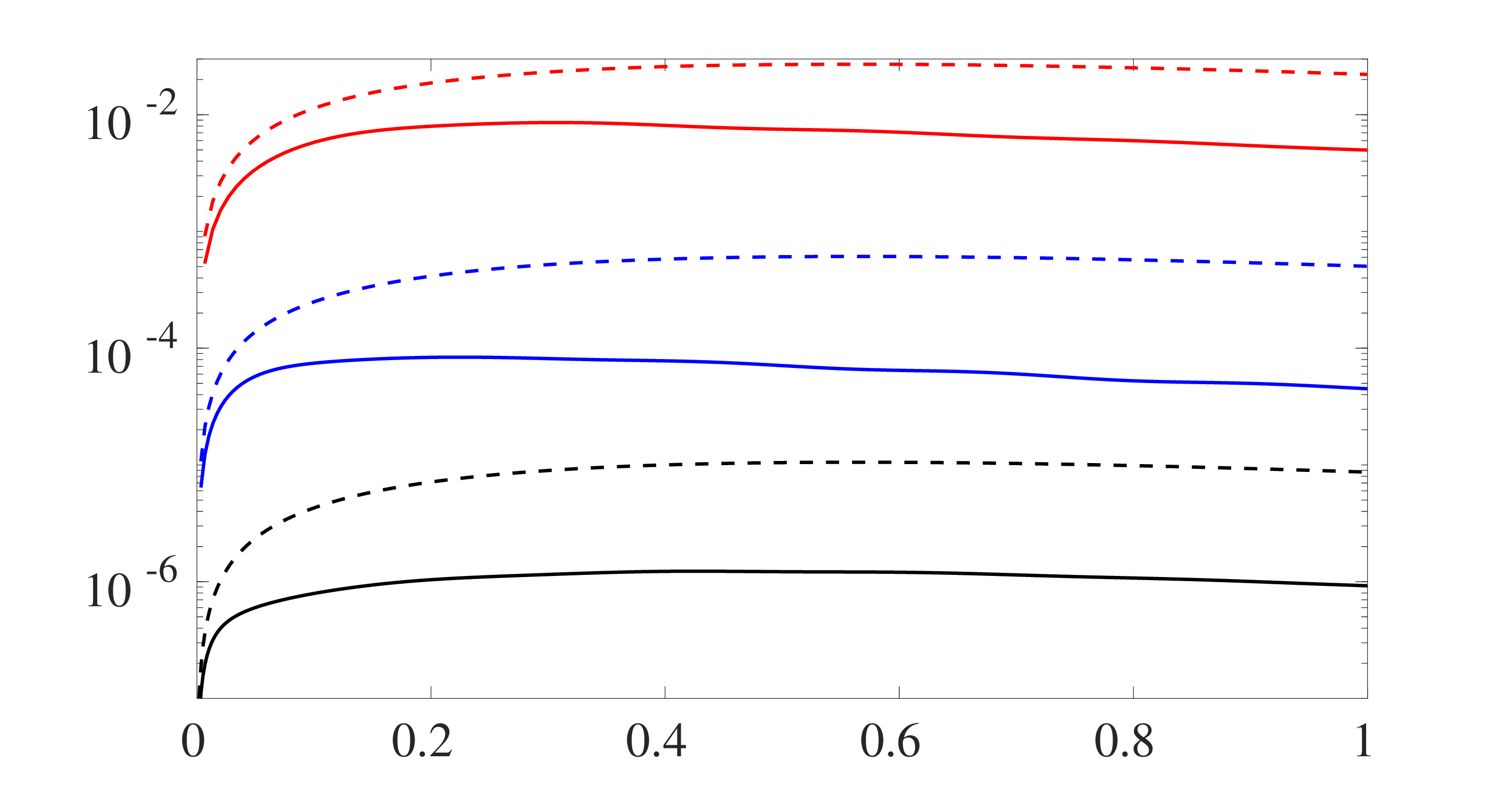}}
    \put(13.0,-0.25){$t$}
    \put(-0.1,1.75){\rotatebox{90}{$L_2-$norm}}
    \put(8.9,1.75){\rotatebox{90}{$L_2-$norm}}
   
    \put(4.0,-0.25){$t$}
    
    \end{picture}
    \end{center}
    \caption{Evolution of $L_2-$norm of solution to advection-diffusion equation for Case 2. Solid lines denote MK approach, dashed lines SK. Left panel: $Re=50$. Right panel: $Re=200$. Solid lines Red $s=1/20$; blue $s=1/40$; black $s=1/80$. }
    \label{fig:AD_case2}
\end{figure}

We now turn our attention to Case 2. In $\S$\ref{subsec::RP_results} greater gains in resolving power of operators were seen along the line $k_x=k_y$ than $k_y=0$, therefore we expect a larger percentage error reduction than in Case 1. Fig. \ref{fig:AD_case2} shows the time evolution of the $L_2-$norm for the same $Re$ as Fig. \ref{fig:AD_case1}. There is a noticeably larger improvement in the present case for both $Re$. Again, simulations at the larger $Re$ value saw larger error reductions, up to an order of magnitude on the finest resolutions. We note that in isotropic homogeneous turbulence there is typically a preponderance of modes with $k_x=k_y$, therefore the MK approach is likely to yield large improvement in many flow fields.

\subsection{Viscous Burgers Equation}
An important feature of many systems of PDEs is nonlinearity. To assess how the new formulation for improving resolving power performs for systems of nonlinear PDEs we consider viscous Burgers equation, a prototypical equation for the development, and subsequent dissipation, of shock waves. It takes the form
\begin{equation}
    \frac{\partial \boldsymbol{u}}{\partial t}+\boldsymbol{u}\cdot\nabla \boldsymbol{u}=\frac{1}{Re} \nabla^2 \boldsymbol{u}.
\end{equation}
We define the time-step as in \eqref{eq::time_step} and solve in a periodic domain $(x,y)\in [0,1]\times [0,1]$. We set the initial conditions to be 
\begin{gather}
    u(x,y,0)=\sin (2\pi x),\qquad v(x,y,0)\equiv 0
\end{gather}
for which the solution (see \cite{cole}) is
\begin{equation}
    \label{BE_soln}
    u(x,t)=\frac{4\pi}{Re}\frac{\sum_{n=1}^\infty nA_n\sin\left(2n\pi x\right)e^{-\frac{n^2\pi^2t}{Re}}}{A_0+\sum_{n=1}^\infty A_n\cos\left(2n\pi x\right)e^{-\frac{n^2\pi^2t}{Re}}}
\end{equation}
where
\begin{equation}
    A_0=e^{-\frac{Re}{4\pi}}I_0\left(\frac{Re}{4\pi}\right),\qquad A_n=2e^{-\frac{Re}{4\pi}}I_n\left(\frac{Re}{4\pi}\right)
\end{equation}
and $I_n$ are modified Bessel functions of the first kind. We use only the first $40$ terms due the decay of the $A_n$ coefficients, which is sufficient to obtain the analytic solution converged to machine precision. Here we consider the case of $Re=10$ as representative example.

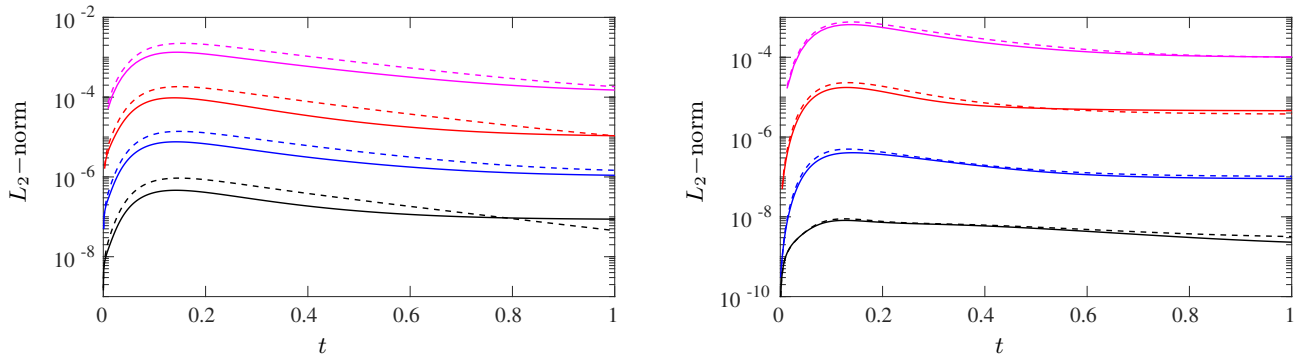
\begin{figure}
    \begin{center}
    \setlength{\unitlength}{1cm}
    \begin{picture}(18,5)(0,0)
    \put(0,0){\includegraphics[width=0.49\linewidth]{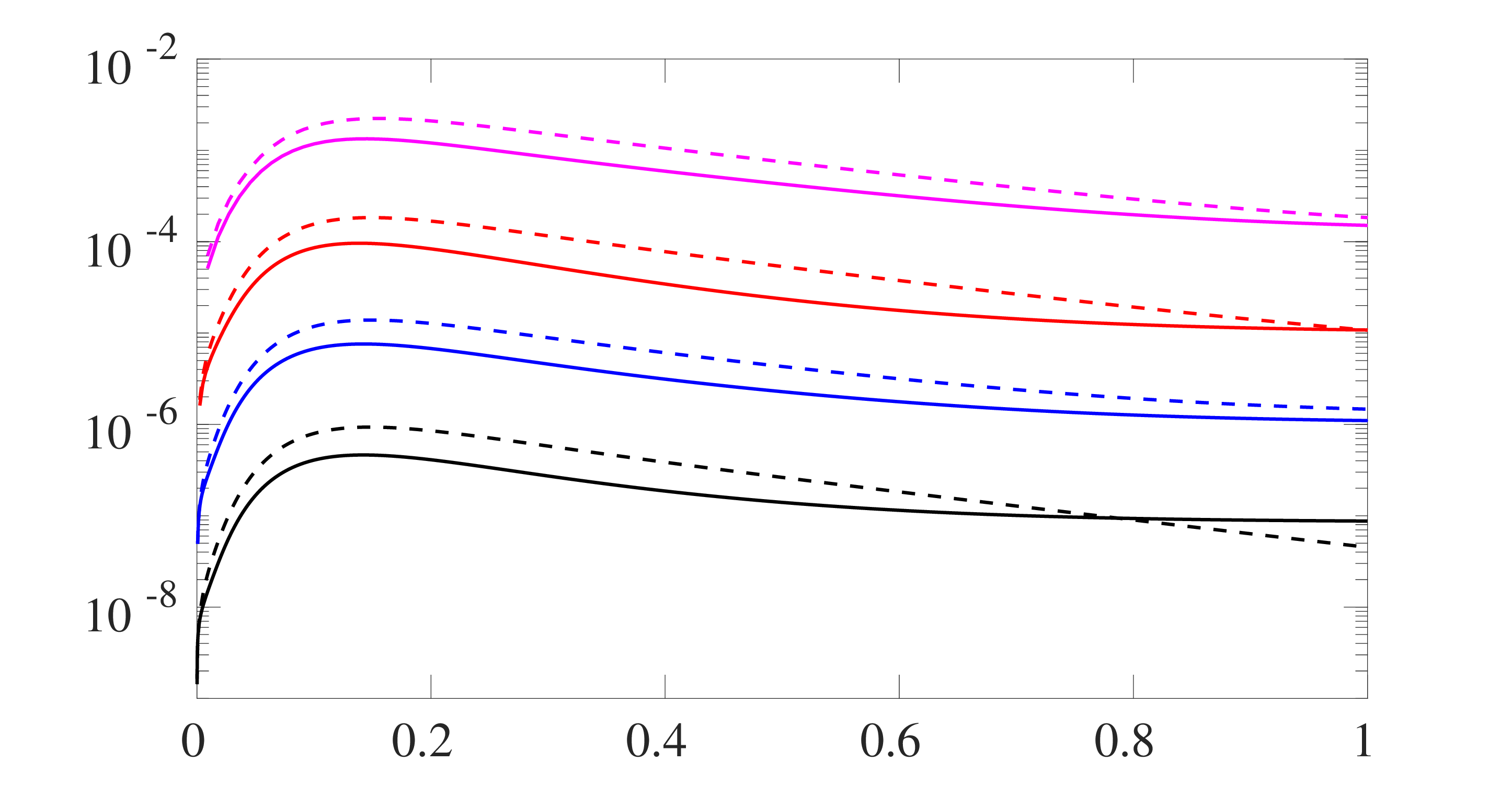}}
    \put(9,0){\includegraphics[width=0.49\linewidth]{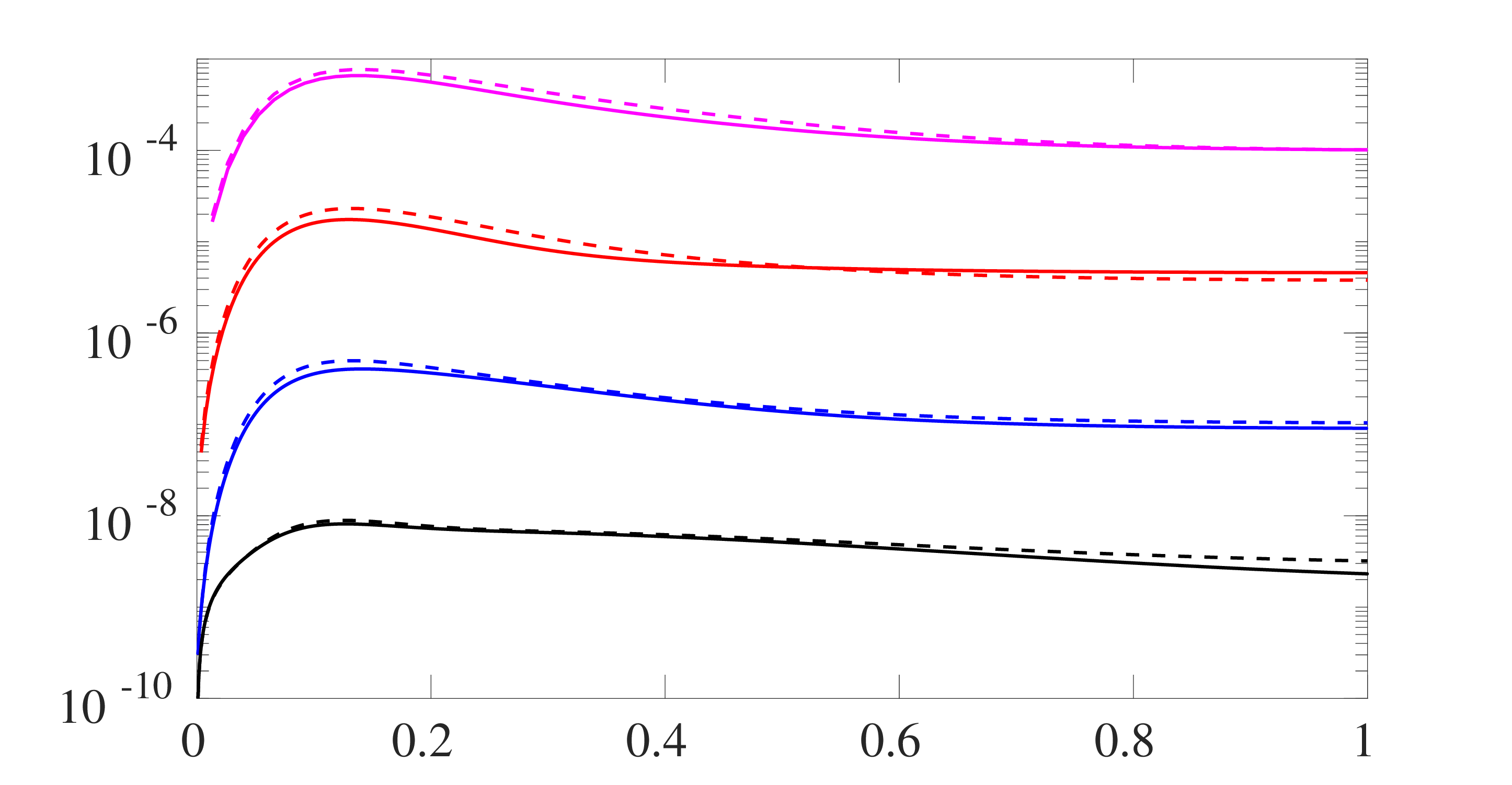}}
    \put(13.0,-0.25){$t$}
    \put(-0.1,1.75){\rotatebox{90}{$L_2-$norm}}
    \put(8.9,1.75){\rotatebox{90}{$L_2-$norm}}
   
    \put(4.0,-0.25){$t$}
    
    \end{picture}
    \end{center}
    \caption{Time evolution of $L_2-$norm of solution to Burgers equation with $Re=10$. Solid lines denote new approach, dashed lines standard approach. Left panel: polynomial consistency of $m=4$. Right panel: polynomial consistency of $m=6$. Solid lines Pink $s=1/10$; red $s=1/20$; blue $s=1/40$; black $s=1/80$. }
    \label{fig:Re10}
\end{figure}

Fig. \ref{fig:Re10} compares time evolution of the $L_2-$norm of the error for both the MK (solid lines) and SK (dashed lines) implementations of $m=4$ (left) and $m=6$ (right) approximations at various resolutions. In all cases, the MK approach presented in $\S$\ref{subsec::Optimisation_of_coefficients} performs better at early times with $20-40\%$ error reduction observed. Maximising resolving power allows the scheme to better capture the high wavenumber components present in the developing shock wave. For computations on a finer mesh the improvement due to maximising resolving power of higher wavenumbers is less pronounced as expected.

However, at later times, when the magnitude of the solution is small, there is no guarantee that the multi-kernel approach provides a more accurate solution. In fact, our numerical experiments suggest that the most accurate scheme for long time behaviour on a given discretisation is independent of resolving power. This is because the long time behaviour is not determined by resolving power, as conservation of properties analogous to linear momentum and kinetic energy are more important \cite{Price_2012}. Like most collocated methods, LABFM is non-conservative. The long-time behaviour of the solution to Burgers equation was found to be correlated to how well material properties were conserved (in the absence of dissipation) by the numerical scheme, which is independent of the current approach of maximising resolving power, although it is linked to the order of polynomial consistency.

We have included the present example in order to highlight the fact that other important factors play a role in determining the accuracy of numerical solutions to PDEs, especially over long time periods. We are actively pursuing a formulation where multiple sets of weights can be used to yield a conservative meshless collocated scheme.

\subsection{Turbulent Flow Test Case: Taylor--Green Vortices}
We now test the improved resolving power provided by the MK formuation using the benchmark test case of three-dimensional Taylor--Green vortices. This problem has often been used to test numerical codes \cite{CFD_methods_2012}, and spectrally accurate reference data is available online \cite{cfd_methods_webpage}. We solve the compressible Navier--Stokes equations in a periodic (in all dimensions) cubic domain with (non-dimensional) side length $2\pi$. The nodal distribution is generated in the same manner as previously in the $(x,y)$ plane, and is uniform in the $z-$direction thus derivatives with respect to $z$ are approximated using $8^{th}$ order accurate finite-differences. The non-dimensional initial conditions are given by
\begin{equation}
    \label{eq:TGV_ICs}
    \begin{split}u(x,y,z,0)=&\sin(x)\cos(y)\cos(z),\\
    v(x,y,z,0)=&-\cos(x)\sin(y)\cos(z),\\
    w(x,y,z,0)=&0,\\
    p(x,y,z,0)=&\frac{\rho_0}{16}\left(\cos(2x)+\cos(2y)\right)\left(2+\cos(2z)\right)
    \end{split}
\end{equation}
To enable comparison to literature we set $Re=1600$, $Ma=0.1$, $Pr=0.71$. The flow quickly transitions to turbulence, with maximum of the dissipation of kinetic energy at non-dimensional time $t\approx 8$. We solve the problem across a range of resolutions from $N=64^3$ to $N=512^3$, comparing results of simulations using: (i) $m=4$ SK approach; (ii) $m=8$ SK approach; (iii) $m=4$ MK approach. As mentioned previously, similar to other collocated methods LABFM admits high wavenumber instabilities, therefore filtering is necessary in order to stabilise the scheme. To enable like for like comparison in the present case we use the same filter (of polynomial consistency $m=8$) for all SK test cases, and implement a MK filter for the MK case. More details on the filtering procedure implemented in LABFM can be found in \cite{king_2022}, and we note the MK filter allows us to set the filter response at a second wavenumber (see Appendix \ref{appendix:MK_filter} for details).

Fig. \ref{fig:TGV_enstrophy} shows the time evolution of the volume averaged enstrophy for all the present simulations with spectrally accurate reference data from \cite{cfd_methods_webpage} included as a baseline comparison. For all but the most resolved simulation improvements resulting from the MK approach compared to the SK approach of both polynomial consistencies are clear graphically. When $N=512^3$, the improvement is more limited though the MK approach still outperforms the SK implementation (see inset).  This suggests that the MK approach with $m=4$ yields better resolution characteristics than the SK approach with $m=8$, and due to the reduced number of neighbours required to approximate operators the computational cost is reduced by approximately $30\%$. These results demonstrate the potential of the MK approach to significantly improve the accuracy of simulations relative to numerical schemes of the same polynomial consistency, and even produce more accurate results than schemes of higher polynomial consistency for substantially smaller computational cost.

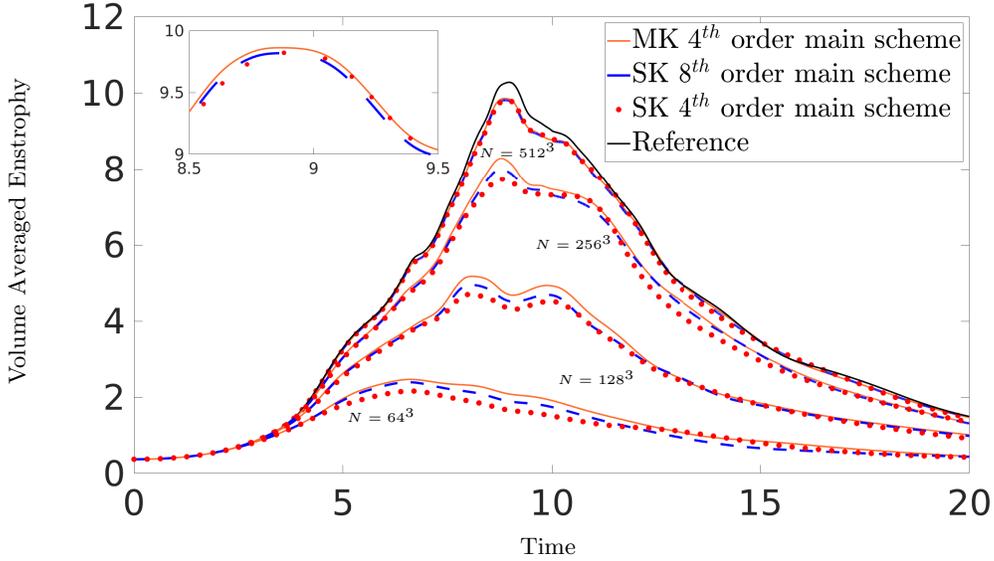
\begin{figure}[t]
    \begin{center}
    \setlength{\unitlength}{1cm}
    \begin{picture}(18,7)(0,0)
    \put(1.8,0){\includegraphics[width=0.8\linewidth]{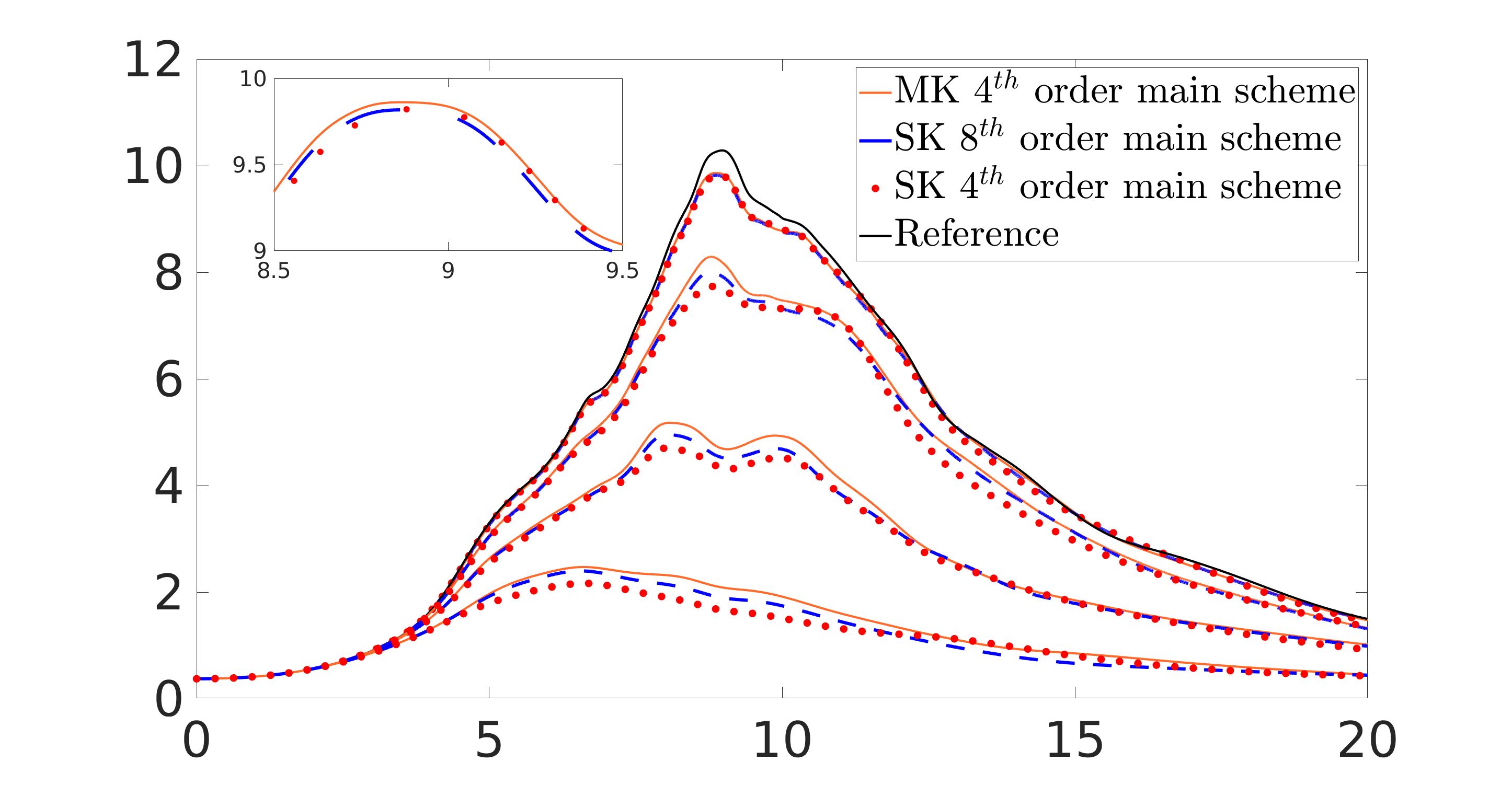}}
    \put(2,2){\rotatebox{90}{Volume Averaged Enstrophy}}
   
    \put(8.8,-0.25){Time}
    \put(8.25,5.03){\tiny{{$N=512^3$}}}
   \put(9,3.8){\tiny{{$N=256^3$}}}
   \put(9.3,2){\tiny{{$N=128^3$}}}
   \put(6.5,1.5){\tiny{{$N=64^3$}}}
    
    \end{picture}
    \end{center}
    \caption{Time evolution of volume averaged enstrophy for three-dimensional Taylor--Green Vortices at $Re=1600, Ma=0.1$ for MK approach $m=4$ (orange solid lines), SK approach $m=4$ (dotted red lines), SK approach $m=8$ (blue dot-dash lines), spectrally accurate reference data of \cite{cfd_methods_webpage}(black line). Inset shows zoom in of $N=512^3$ case near the global maximum (Time $\in[8.5, 9.5]$). }
    \label{fig:TGV_enstrophy}
\end{figure}

\section{Conclusion}
\label{sec:conclusion}
Due to the asymmetry of node distributions, the errors in derivative approximations in meshless methods are dependent on both the magnitude and orientation of the wavenumber vector of the field being approximated. An analysis of the resolving power of meshless methods needs to account for this anisotropy in wavenumber space, and in this work we have presented a framework for such analysis. We have exploited the non-uniqueness of weights used to approximate operators in meshless methods by using multiple bases/kernels simultaneously in order to reduce error. To do so, we have optimised the resolving power of each operator. By maximising resolving power, the new approach provides significant gains in accuracy when approximating operators for three different numerical methods. This approach is applicable to a wide class of collocated numerical methods, far more than can be detailed in one manuscript. Collocated numerical methods which are suitable to improvement in this manner will have the key features:
\begin{itemize}
\item The domain is discretised with a scattered nodal distribution.
    \item Operators can be expressed as a weighted sum of node properties within a computational stencil.
    \item For a given stencil and desired polynomial consistency, the weights are not unique.
\end{itemize}

When canonical PDEs were simulated, the new approach provides consistent reductions in error for regimes where resolving high wavenumber components is important. In time-dependent cases where the magnitude of the solution exponentially decays, other properties such as conservation may be more important. However, there is clear potential for use of a multi-kernel approach in simulations of turbulent flows performed using meshless methods, where resolving the smallest scales admitted by the discretisation is vital.

Our focus herein has been on explicit, Eulerian methods. The same approach is possible in a Lagrangian framework, however, the computational cost of calculating an extra set of weights and then determining an optimal linear combination is not negligible. Both Radial Basis Functions Finite-Differences (RBF-FDs) and the local anisotropic basis function method (LABFM) are typically used in an Eulerian framework so this is not an issue for these methods. Smoothed Particle Hydrodynamics (SPH) on the other hand is widely used in an explicit Lagrangian framework, although calculation of an extra set of weights would not be a major bottleneck for such simulations. We believe that the benefits of the new approach would outweigh costs when simulating turbulent flows, though this likely depends on the flow scenario in question. Comparison with more expensive higher order corrections would also be of interest in such flows.

It is of course possible to formulate meshless methods as implicit schemes, and the polynomial consistency of RBF-FDs and LABFM means there is a natural extension to compact schemes. However, without a concept of resolving power in meshless methods it is difficult to know how to choose the implicit coefficients. By providing the resolving power of meshless methods with a firm mathematical basis it will now be possible to optimise the choice of implicit coefficients in compact meshless schemes. Development of a compact formulation of LABFM optimised for resolving power is ongoing, with details reserved for a future study.

\begin{acknowledgements}
This work was partially funded by the Engineering and Physical Sciences Research Council (EPSRC) grant EP/W005247/2. JK is supported by the Royal Society via a University Research Fellowship (URF\textbackslash R1\textbackslash 221290). We thank EPSRC for computational time made available on the UK supercomputing facility ARCHER2 via the UK Turbulence Consortium (EP/X035484/1). We are grateful for assistance given by Research IT and the use of the Computational Shared Facility at the University of Manchester. 
\end{acknowledgements}

\appendix

\section{Numerical Implementations}
\label{sec::Numerical_implementations}

We briefly outline the numerical methods used in this manuscript chosen due to their desirable convergence properties or widespread use in the meshless community. Other collocated meshless methods of course exist, and as long as the conditions described in the conclusion are satisfied the analysis presented herein remains valid.

\subsection{SPH}
To approximate operators using SPH, we must choose a kernel and a stencil radius. For MK approximations we use a combination of the Wendland C2 kernel and Gaussian kernel with the same stencil size of $2h$, both with a constant stencil size throughout the domain. We use the antisymmetric gradient operator which guarantees zeroth order consistency \cite{quinlan_2006,king_2020}, and upon using the notation $(\cdot)_{ji}=(\cdot)_j-(\cdot)_i$ the gradient approximation of a test function $\phi$ is 
\begin{equation}
    \nabla \phi\Big|_i=\sum_{j\in\mathcal{N_i}}\phi_{ji}\nabla W_{ij} V_j
\end{equation}
where $\mathcal{N}_i$ is the set of neighbours of $i$ within a stencil of size $2h_i$, $\nabla W$ is the kernel gradient, and $V_j$ is the volume associated with the point $j$. Note that throughout this paper $h_i$ and $V_j$ are taken to be constant (when considering SPH).

To approximate the Laplacian we use the Morris operator \cite{Morris}
\begin{equation}
   \nabla^2 \phi_i = 2\sum_j \phi_{ji}\boldsymbol{e}_{ji}\cdot \nabla W_{ij} V_j/r_{ji}
\end{equation}
where $\boldsymbol{e}_{ji}$ is the unit vector pointing from $j$ to $i$ and $r_{ji}$ the distance between $j$ and $i$. Formally, this operator is divergent, with leading order error terms $\mathcal{O}\left(h_{0}\right)$, although for a range of resolutions higher order error terms dominate, and convergence is observed.

Throughout this work, for ease of computation we fix the particle positions and use an Eulerian SPH scheme. Doing so negates some of the advantages of SPH compared to other methods, however the goal of this paper is to demonstrate how meshless methods can be improved by exploiting the non-uniqueness of weights, not to compare the numerical methods themselves. Therefore an Eulerian SPH scheme suffices for the present proof of concept study.

\subsection{RBF-FDs}
RBF-FDs approximations to an operator $L$ acting on a function $\phi$ with polynomial consistency $m$ take the form
\begin{equation}
    L(\phi)\Big|_i=\sum_{j\in\mathcal{N}_i}\phi_{ji}w_{ji}^L+\mathcal{O}(s^{m+1})
\end{equation}
where the $w_{ji}^L$ are weights which are implicitly dependent on the radial basis function chosen $\psi(||\boldsymbol{x}_{ji}||)$, and are found by solving (in 2D) an $(\lvert\mathcal{N}_i\rvert+1+(m^2+3m)/2)\times (\lvert\mathcal{N}_i\rvert+1+(m^2+3m)/2)$ linear system.
Herein we use the Gaussian ($\psi_1$) RBF for the SK implementation and the inverse multiquadric ($\psi_2$) RBF for the MK approach (see \cite{Fornberg&Flyer2015} for details). When calculating weights we choose to use a constant number of nearest neighbours $\lvert\mathcal{N}_i\rvert$ rather than defining a stencil size. This number increases for higher order approximations. There is a complex interplay between the choice of $\lvert\mathcal{N}_i\rvert$, the stability of the operator, and the accuracy of the approximation, acknowledged in various studies e.g. \cite{Tominec_2021}. Our choices are made to ensure the stability of RBF-FD approximations is comparable to the other methods considered herein (larger stencils are typically more unstable \cite{Fornberg&Flyer2015} requiring hyperviscosity to stabilise time-dependent simulations), and are similar to \cite{Fornberg2011} and \cite{Kolar_RBF-FD_stencils}.

For each RBF we must also choose a flatness parameter $\varepsilon$. We allow this to vary on a stencil by stencil basis, instead requiring at each node $i$
\begin{equation}
    \mathrm{min}\{ \psi_1(||\boldsymbol{x}-\boldsymbol{x}_i||)\}=\frac{1}{2},\qquad \mathrm{max }\{ \psi_2(||\boldsymbol{x}-\boldsymbol{x}_i||)\}=\frac{5}{4}
\end{equation}
which ensures that the linear system remains well-conditioned and the operators remain relatively stable. We note that there is no definitive consensus in the RBF-FD community for choices of $\lvert N_i\rvert$ and $\varepsilon$ despite their influence on both stability and accuracy. Although tangential to the main topic of this paper, we note that the resolving power analysis in \S\ref{sec::resolving_power_formulation} gives a metric by which different choices of $\lvert\mathcal{N}_i\rvert$ and $\varepsilon$ may be compared.

\subsection{LABFM}
LABFM operators of polynomial consistency $m$ take the form
\begin{equation}
    L(\phi)\Big|_i= \sum_{j\in\mathcal{N}_i}\phi_{ji}w_{ji}^L+\mathcal{O}(s^{m+1})
\end{equation}
where $w_{ji}^L$ are weights constructed from a linear sum of anisotropic basis functions. To obtain these weights a linear system of size $(m^2+3m)/2 \times (m^2+3m)/2$ must be solved. Herein we use Hermite ABFs with either a Wendland C2 kernel (SK implementation) and/or a Gaussian kernel (second kernel in MK approach). Complete details are given in \cite{king_2020,king_2022}. Increasing $m$ requires more points to be used in the calculation of the weights, hence a larger stencil must be used. Stencil size is determined using the procedure outlined in \cite{king_2022}, with the stencil size is systematically reduced from a stable value whilst ensuring the calculation of weights remains well-conditioned.

\section{Multi-Kernel Filter Operator}
\label{appendix:MK_filter}
Collocated numerical methods often admit high wavenumber instabilities, and therefore require filtering to prevent spurious modes from destroying the numerical solution. Here we briefly discuss how a multi-kernel approach can be used to improve filter operators in meshless methods. As in \cite{king_2022}, time-dependent variables are filtered at each time-step through
\begin{equation}
    \Tilde{\phi}_i=(1+\kappa_{m,i}\nabla^m)\phi \lvert_i
\end{equation}
where $\kappa_{m,i}$ is a constant to be chosen (though it can differ at each point $i$), and $m$ is the order of the filter operator. Approximating the filter operator in a MK framework gives
\begin{equation}
    \nabla^m \phi\lvert_i=\sum_{j\in \mathcal{N}_i}\phi_{ji}\left(\hat{c}_i\hat{w}_{ji}^{\mathrm{hyp}}+(1-\hat{c}_i)\overline{w}_{ji}^{\mathrm{hyp}}\right).
\end{equation}
In \cite{king_2022}, $\kappa_{m,i}$ was chosen in order to set the response of the filter operator at a given point in wavenumber space (response of 2/3 at a wavenumber $k_x=k_y=(2/3) k_{Ny}$). In the MK framework we can now choose the response of the filter operator at two different points in wavenumber space. To illustrate this let's assume we wish to set the filter operator at wavenumbers $\boldsymbol{k}_1, \boldsymbol{k}_2$ to the values $\lambda_1, \lambda_2$. We require
\begin{subequations}
    \begin{gather}
        \kappa_{m,i}\sum_{j\in \mathcal{N}_i}\left(1-\cos(\boldsymbol{k}_1\cdot \boldsymbol{x}_{ji})\right)\left(\hat{c}_i\hat{w}_{ji}^{\mathrm{hyp}}+(1-\hat{c}_i)\overline{w}_{ji}^{\mathrm{hyp}}\right)=\lambda_1,\\
        \kappa_{m,i}\sum_{j\in \mathcal{N}_i}\left(1-\cos(\boldsymbol{k}_2\cdot \boldsymbol{x}_{ji})\right)\left(\hat{c}_i\hat{w}_{ji}^{\mathrm{hyp}}+(1-\hat{c}_i)\overline{w}_{ji}^{\mathrm{hyp}}\right)=\lambda_2.
    \end{gather}
\end{subequations}
These are two equations for two unknowns ($\kappa_{m,i}$, $\hat{c}_i$) which can be solved trivially to give the desired filter responses. Choices of $\boldsymbol{k}_1, \boldsymbol{k}_2,\lambda_1,\lambda_2$ are made to ensure sufficient filtering at high wavenumbers whilst requiring the response remains as small as possible at low wavenumbers. For all simulations with filtering performed herein we have taken $m=8$.

  \bibliographystyle{elsarticle-num} 
  \bibliography{HB_bib}

\begin{thebibliography}{10}
\expandafter\ifx\csname url\endcsname\relax
  \def\url#1{\texttt{#1}}\fi
\expandafter\ifx\csname urlprefix\endcsname\relax\def\urlprefix{URL }\fi
\expandafter\ifx\csname href\endcsname\relax
  \def\href#1#2{#2} \def\path#1{#1}\fi

\bibitem{Lele}
S.~K. Lele, Compact {F}inite {D}ifference {S}chemes with {S}pectral-like {R}esolution, J. Comp. Phys. 103 (1992) 16--42.

\bibitem{Li_Review}
S.~Li, W.~K. Liu, Meshfree and particle methods and their applications, Appl. Mech. Rev. 55 (2002) 1--34.

\bibitem{Garg_Review}
S.~Garg, M.~Pant, Meshfree methods: a comprehensive review of applications, Int. J. Comput. Methods 15 (2018).

\bibitem{gingold}
R.~A. Gingold, J.~J. Monaghan, Smoothed particle hydrodynamics: theory and application to non-spherical stars, Mon. Not. R. Astron. Soc 181 (1977) 375--389.

\bibitem{Lucy}
L.~B. Lucy, Numerical approach to the testing of the fission hypothesis, Astron. J. 82 (1977) 1013--1024.

\bibitem{Monaghan_Review}
J.~J. Monaghan, Smoothed {P}article {H}ydrodynamics and its diverse applications, Ann. Rev. Fluid Mech. 44 (2012) 323--346.

\bibitem{Bonet&Lok}
J.~Bonet, T.-S. Lok, Variational and momentum preservation aspects of smooth particle hydrodynamics formulations, Comput. Meth. Appl. Mech. Eng. 180 (1999) 97--115.

\bibitem{Sibilla}
S.~Sibilla, An algorithm to improve consistency in smoothed particle hydrodynamics, Comput. Fluids 118 (2015) 148--158.

\bibitem{Cabezon&Garcia}
R.~M. Cabez{\'{o}}n, D.~Garc{\'{i}a}-Senz, Mixing {S}inc kernels to improve interpolations in smoothed particle hydrodynamics without pairing instability, Mon. Not. R. Astron. Soc 528 (2024) 3782–3796.

\bibitem{MK_Arxiv}
R.~Wissing, T.~Quinn, B.~Keller, J.~Wadsley, S.~Shen, Optimized smoothing kernels for {SPH} (2025).
\newblock \href {http://arxiv.org/abs/2508.08471} {\path{arXiv:2508.08471}}.

\bibitem{TraskGMLS}
N.~Trask, M.~Perego, P.~Bochev, A high-order staggered meshless method for elliptic problems, SIAM J. Sci. Comput. 39 (2017) A479–A502.

\bibitem{TraskCMLS}
M.~Trask, N.~Maxey, X.~Hu, Compact moving least squares: An optimization framework for generating high-order compact meshless discretisations, J. Comp. Phys. 326 (2016) 596--611.

\bibitem{Weinan}
E.~Weinan, J.-G. Liu, Essentially compact schemes for unsteady viscous incompressible flow, J. Comp. Phys. 126 (1996) 122--138.

\bibitem{RBF_Cartography}
R.~L. Hardy, Multiquadric equations of topography and other irregular surfaces, J. Geophys. Rev. 76 (1971) 1905--1915.

\bibitem{shu}
C.~Shu, H.~Ding, K.~Yeo, Local radial basis function-based differential quadrature ann its application to solve two-dimensional incompressible {N}avier--{S}tokes equations, Comput. Methods Appl. Mech. Eng. 192 (2003) 941--954.

\bibitem{Cecil}
T.~Cecil, J.~Qian, S.~Osher, Numerical methods for high dimensional {H}amilton--{J}acobi equations using radial basis functions, J. Comp. Phys. 196 (2004) 327--347.

\bibitem{Wright}
G.~B. Wright, Radial {B}asis {F}unction {I}nterpolation: {N}umerical and {A}nalytic {D}evelopments, Ph.D. thesis, Boulder (2004).

\bibitem{Fornberg&Flyer2015}
B.~Fornberg, N.~Flyer, Solving pdes with radial basis functions, Acta Numer. 24 (2015) 215--258.

\bibitem{Hybridkernels_RBFFD}
P.~K. Mishra, G.~E. Fausshauer, M.~K. Sen, L.~Ling, A stabilized radial basis-finite difference {(RBF-FD)} method with hybrid kernels, Comput. Math. Appl. 77 (2019) 2354--2368.

\bibitem{king_2020}
J.~R.~C. King, S.~J. Lind, A.~M.~A. Nasar, High order difference schemes using the local anisotropic basis function method, Journal of Computational Physics 415 (2020) 109549.

\bibitem{king_2022}
J.~R.~C. King, S.~J. Lind, High-order simulations of isothermal flows using the local anisotropic basis function method ({LABFM}), Journal of Computational Physics 449 (2022) 110760.

\bibitem{king_2024_ve}
J.~R.~C. King, S.~J. Lind, A mesh-free framework for high-order simulations of viscoelastic flows in complex geometries, Journal of Non-Newtonian Fluid Mechanics 330 (2024) 105278.

\bibitem{king_2024_Combustion}
J.~R.~C. King, A mesh-free framework for high-order direct numerical simulations of combustion in complex geometries, Comp. Meth. Appl. Mech. Eng. 421 (2024) 116762.

\bibitem{Roberts&Weiss}
K.~V. Robest, N.~O. Weiss, Convective {D}ifference {S}chemes, Math. Comput 20 (1966) 272--299.

\bibitem{Vichnevetsky&Bowles}
R.~Vichnevetsky, J.~B. Bowles, Fourier {A}nalysis of {N}umerical {A}pproximations of {H}yperbolic {E}quastion, SIAN, 1982.

\bibitem{Pencil_code}
A.~Brandenburg, Computational aspects of astrophysical {MHD} and turbulence, in: A.~Ferriz-Mas, M.~N{\'u}{\~n}ez (Eds.), Advances in nonlinear dynamos, Taylor \& Francis, London and New York, 2003, p. 269–334.

\bibitem{Fornberg&Flyer_disc}
B.~Fornberg, N.~Flyer, Fast generation of 2{D} node distributions for mesh-free {PDE} discretizations., Comput. Math. Appl. 69 (2015) 531--544.

\bibitem{Bayona_2017}
V.~Bayona, N.~Flyer, B.~Fornberg, G.~A. Barnett, On the role of polynomials in {RBF}-{FD} approximations: {II}. {N}umerical solution of elliptic {PDE}s, J. Comp. Phys. 332 (2017) 257--273.

\bibitem{cole}
J.~D. Cole, On a quasi-linear parabolic equation occurring in aerodynamics, Q. Appl. Math. 9 (1951) 225--236.

\bibitem{Price_2012}
D.~J. Price, Smoothed particle hydrodynamics and magnetohydrodynamics, J. Comp. Phys. 231 (2012) 759--794.

\bibitem{CFD_methods_2012}
Z.~Wang, K.~Fidowski, R.~Abgrall, F.~Bassi, D.~Caraeni, A.~Cary, H.~Deconinck, R.~Hartmann, K.~Hillewaert, H.~Huynh, K.~N., G.~May, P.-O. Persson, B.~van Leer, l.~M. Visba, High-order {CFD} methods: current status and perspective, Int. J. Numer. Methods Fluids 72 (2012) 811--845.

\bibitem{cfd_methods_webpage}
Z.~J. Wang, R.~Abgrall, F.~Bassi, D.~Caraeni, A.~Cary, H.~Deconinck, C.~Fidkowski, R.~Hartmann, K.~Hillewaert, H.~T. Huynh, N.~Kroll, G.~May, P.-O. Persson, B.~van Leer, M.~Visbal, 1st international workshop on high-order {CFD} methods, \url{https://cfd.ku.edu/hiocfd.html}, accessed: 2025-10-03 (2012).

\bibitem{quinlan_2006}
N.~J. Quinlan, M.~Basa, M.~Lastiwka, Truncation error in mesh-free particle methods, Int. J. Num. Meth. Engng 66 (2006) 2064–2085.

\bibitem{Morris}
J.~P. Morris, P.~J. Fox, Y.~Zhu, Modeling {L}ow {R}eynolds {N}umber {I}ncompressible {F}lows {U}sing {SPH}, J. Comp. Phys 136 (1997) 214--226.

\bibitem{Tominec_2021}
I.~Tominec, E.~Larsson, H.~A., A least squares radial basis function finite difference method with improved stability properties, SIAM J. Sci. Comput. 43 (2021).

\bibitem{Fornberg2011}
B.~Fornberg, E.~Lehto, Stabilization of {RBF}-generated finite difference methods for convective {PDE}s, J. Comp. Phys 230 (2011) 2270--2285.

\bibitem{Kolar_RBF-FD_stencils}
A.~Kolar-Po{\v{z}}an, M.~Jan{\v{c}}i{\v{c}}, M.~Rot, G.~Kosec, Some observations regarding the {RBF-FD} approximation accuracy dependence on stencil size, J. Comp. Sci. 79 (2024).

\end{thebibliography}

\end{document}